\documentclass[pdflatex,sn-vancouver]{sn-jnl}

\usepackage{graphicx}%
\usepackage{rotating}
\usepackage{adjustbox}
\usepackage{multirow}%
\usepackage{upgreek}
\usepackage{amsmath,amssymb,amsfonts}%
\usepackage{amsthm}%
\usepackage{mathrsfs}%
\usepackage[title]{appendix}%
\usepackage{xcolor}%
\usepackage{textcomp}%
\usepackage{manyfoot}%
\usepackage{booktabs}%
\usepackage{algorithm}%
\usepackage{algorithmicx}%
\usepackage{algpseudocode}%
\usepackage{listings}%
\usepackage{subcaption}
\usepackage{url}

\theoremstyle{thmstyleone}%
\theoremstyle{thmstyletwo}%

\theoremstyle{thmstylethree}%

\begin{document}

\title[Article Title]{AI Driven Laser Parameter Search: Inverse Design of Photonic Surfaces using Greedy Surrogate-based Optimization}


\author*[1]{\fnm{Luka} \sur{Grbcic}}\email{lgrbcic@lbl.gov}

\author[2]{\fnm{Minok} \sur{Park}}\email{minokpark@lbl.gov}

\author[3]{\fnm{Juliane} \sur{M\"uller}}\email{juliane.mueller@nrel.gov}
\author[2,4]{\fnm{Vassilia} \sur{Zorba}}\email{vzorba@lbl.gov}
\author*[1]{\fnm{Wibe Albert} \sur{de Jong}}\email{wadejong@lbl.gov}

\affil*[1]{\orgdiv{Applied Mathematics and Computational Research}, \orgname{Lawrence Berkeley National Laboratory}, \orgaddress{\street{1 Cyclotron Rd}, \city{Berkeley}, \postcode{94720}, \state{California}, \country{USA}}}

\affil[2]{\orgdiv{Energy Technologies Area}, \orgname{Lawrence Berkeley National Laboratory}, \orgaddress{\street{1 Cyclotron Rd}, \city{Berkeley}, \postcode{94720}, \state{California}, \country{USA}}}

\affil[3]{\orgdiv{Computational Science Center}, \orgname{National Renewable Energy Laboratory}, \orgaddress{\street{15013 Denver West Parkway}, \city{Golden}, \postcode{80401}, \state{Colorado}, \country{USA}}}

\affil[4]{\orgdiv{Department of Mechanical Engineering}, \orgname{University of California at Berkeley}, \orgaddress{\street{6141 Etcheverry Hall}, \city{Berkeley}, \postcode{94709}, \state{California}, \country{USA}}}

\abstract{Photonic surfaces designed with specific optical characteristics are becoming increasingly important for use in in various energy harvesting and storage systems. , In this study, we develop a surrogate-based optimization approach for designing such surfaces. The surrogate-based optimization framework employs the Random Forest algorithm and uses a greedy, prediction-based exploration strategy to identify the laser fabrication parameters that minimize the discrepancy relative to a user-defined target optical characteristics. We demonstrate the approach on two synthetic benchmarks  and  two specific cases of photonic surface inverse design targets. It exhibits superior performance when compared to other optimization algorithms across all benchmarks. Additionally, we demonstrate a technique of inverse design warm starting for changed  target optical characteristics which  enhances the performance of the introduced approach.}


\maketitle

\section{Introduction}\label{sec:intro}

The photonic surface is a type of material that is excellent at absorbing light and emitting thermal radiation. Its efficiency depends on its spectral absorptivity and emissivity, more specifically, how well it emits energy across different wavelengths when at a stable temperature. These quantities measure the energy emitted at each wavelength compared to that of an ideal emitter (\citet{brewster1992thermal}, \citet{howell2020thermal}). Photonic surfaces are increasingly used for energy applications like harvesting and storage. They are used in Thermophotovoltaic (TPV) systems (\citet{fan2020near}, \citet{lapotin2022thermophotovoltaic}), radiative cooling systems (\citet{raman2014passive}, \citet{heo2020janus}), solar-based water desalination systems (\citet{menon2020enhanced}, \citet{ni2016steam}), and concentrated solar power systems (\citet{weinstein2015concentrating}, \citet{he2020perspective}). Designing photonic surfaces to meet specific target spectral emissivity values is therefore a key optimization and inverse design task. 

Modern approaches used for photonic materials inverse design can be sorted into two main categories: (i) Deep Learning (DL) (\citet{wiecha2021deep}, \citet{ma2021deep}), and (ii) Optimization-based methods (\citet{mao2021inverse}, \citet{wang2021intelligent}). DL-based approaches are increasingly used and they include architectures like Tandem Neural Networks (\citet{xu2021improved}, \citet{park2024tnn}), Generative Adversarial Networks (GAN) (\citet{ma2022benchmarking}, \citet{jiang2020simulator}), and Autoencoders (AE) (\citet{wiecha2021deep}, \citet{kudyshev2020machinea}). These techniques are popular since they handle unstructured data that are common for parametrization of photonic materials (\citet{liu2020topological}, \citet{liu2021tackling}, \citet{kudyshev2020machinea}). DL methods are beneficial as they can be reused for inverse design, provided that the same design space is considered. However, these methods require large amounts of data to achieve sufficient accuracy and  this ``cost'' also needs to be taken into account (\citet{habibi2023actually}). Moreover,  if the design space of the target value drastically changes compared to that of the training data, DL approaches can struggle to accurately extrapolate.

Optimization-based methods for photonic inverse design can  be further divided into two distinct categories. The first category is photonic inverse design using the adjoint optimization method (\citet{zhu2023inverse}, \citet{gershnabel2022reparameterization}, \citet{hughes2018adjoint}, \citet{minkov2020inverse}, \citet{lalau2013adjoint}, \citet{wang2020inverse}). The adjoint optimization method is a specific type of the gradient-based optimization method that is generally  computationally more efficient for inverse design (\citet{molesky2018inverse}). The main drawback is in its larger implementation complexity, as well as its dependency on simulations, i.e., it is not suitable for inverse design based on experimental data (\citet{ma2021deep}).

The other  category of optimization-based approaches are hybrid machine learning (ML) and optimization methods. A forward ML model is trained, also known as a surrogate model, that can quickly and efficiently predict a solution given on a design vector. The next step is to formulate the problem as an inverse design optimization problem that is solved by an optimization algorithm which employs the surrogate model to assess each design (\citet{deng2022hybrid}, \citet{ma2021deep}). Another option is to use an ML surrogate to infer initial designs used for further simulation-optimization steps which has shown to be generally beneficial for inverse design tasks (\citet{habibi2023actually}). Hybrid approaches include combinations like AE coupled with Differential Evolution (DE) (\citet{kudyshev2020machineb}), fully connected DNNs with DE (\citet{hegde2019photonics}), adjoint optimization with AEs (\citet{kudyshev2020machinea}), Long Short-term Memory networks with a gradient-free optimization algorithm (\citet{yao2023inverse}), adjoint optimization with GANs (\citet{kudyshev2021optimizing}), \citet{yeung2022deepadjoint}), adjoint optimization with Convolutional Neural Networks (\citet{yeung2022enhancing}), and Random Forests (RF) with DE (\citet{grbcic2024ensemble}). 

The overwhelming majority of previous literature are based on simulation data, and most have a unique design methodology tailored for a specific photonic material type and target property. Furthermore, design parameters often include complex topological information (\citet{liu2020topological}, \citet{liu2021tackling}, \citet{kudyshev2020machinea}), \citet{kudyshev2020machineb}). Finally, previous research of photonic materials inverse design is mostly not focused on minimizing the required resources needed to train highly accurate ML models, nor the minimization of inverse design (simulation or experimental) function evaluations with optimization algorithms. This is an increasingly important aspect of inverse design that needs to be considered, especially in the era of self-driving laboratories and autonomous experimentation (\citet{noack2023methods}, \citet{hase2019next}). 

The inverse design approach we introduce here builds upon our previous work (\citet{park2024tnn}, \citet{grbcic2024ensemble}) where we showed that it is possible to design photonic surfaces utilizing pulsed femtosecond laser ablation on surfaces of plain materials. We determine an inverse design mapping between laser fabrication parameters such as laser power, scanning speed, and spacing between consecutive scan lines, thereby bypassing the complex light-matter interaction physics, and complex topological data usually used for photonic inverse design. Furthermore, unlike the majority of previous inverse design methods, we utilize real experimental data originally presented in our previous work. The experimental data contains laser fabrication parameters and spectral emissivity curves for two different materials, namely stainless steel and Inconel. 

Given all these major benefits over other approaches, we introduce the AI Laser Parameter Search (ALPS) optimization framework used for efficient inverse design of photonic surfaces. 
ALPS is a surrogate-based optimization approach that utilizes the RF algorithm and a prediction-based exploration strategy (\citet{kochenderfer2019algorithms}) to approximate a target design. This approach relies on evaluating the RF surrogate at specific designs, and then determining the optimal design or a batch of optimal designs for experimental model evaluation based on its distance or discrepancy relative to the target design. This simple greedy sampling strategy is shown to work efficiently on a set of benchmark problems (\citet{paria2022greedy}, \citet{wang2014evaluation}). The RF algorithm is used as it is shown to perform accurately and efficiently as a forward model for photonic surfaces in previous work (\citet{grbcic2024ensemble}, \citet{elzouka2020interpretable}). 

Furthermore, one of the major benefits of ALPS is that the RF algorithm is trained to model the forward relationship of the photonic surface design (laser fabrication parameters and spectral emissivity) during the optimization process, and not the relationship between the input designs and their discrepancy relative to the target which is usually the case for most model-based optimization algorithms for inverse design. This feature of ALPS implies re-usability through warm starting when the inverse design target is changed. It should be noted that this could also be achieved through Bayesian Optimization (BO) methods that utilize Gaussian Processes (GP) as they are extremely efficient, however, they require more complex implementations (\citet{liu2018remarks}) if we want to accurately capture the multi-input and multi-output relationships present in the photonic surface design, whereas in ALPS we can simply utilize an out-of-the-box RF algorithm.

Finally, in order to showcase the benefits of using ALPS for general and photonic surfaces inverse design optimization, we compare it with an array of other established optimization algorithms (including BO) on two synthetic benchmarks and four different photonic surfaces benchmark targets varied by both the used plain surface material (Stainless steel and Inconel) and the shape of the the target spectral emissivity curves. Furthermore, we show the re-usability of the ALPS approach for a cross-target and even cross-material inverse design process.

The remainder of this article is organized as follows. The photonic surface inverse design problem is introduced in Sec. \ref{sec:InverseDesign}. This section includes the formal mathematical definition of inverse design, the photonic surface design loop, the experimental models, the photonic surface design target benchmarks, and most importantly, the ALPS algorithm. The results of ALPS and other optimization algorithms for the synthetic and photonic inverse design benchmarks (with and without warm starting) are presented in Sec. \ref{sec:results}. Details on the experimental data, machine learning experimental models and validation, synthetic benchmarks, algorithms used for comparison, and detailed results can be found in App. \ref{app:experimental_model}, \ref{app:exp_setup}, and \ref{app:fullresults}, respectively.

\section{Photonic Surface Inverse Design}\label{sec:InverseDesign}

In this section, we introduce the mathematical definition of the general inverse design problem and notation, as well as the photonic surface inverse design problem. Furthermore, we define the photonic surface inverse design benchmarks and the algorithmic details of ALPS.

\subsection{Mathematical Definition of Inverse Design}\label{sec:inversemath}

The objective of inverse design is to obtain a set of design parameters that yield a known target value or property. The inverse design problem is mathematically defined as:

\begin{equation}
\mathbf{x} = f^{-1}(\mathbf{y})
\label{eqn:inverse_design_definition}
\end{equation}

In Eq. (\ref{eqn:inverse_design_definition}), \(\mathbf{x} \in \mathbb{R}^M\)  is the design vector, \(\mathbf{y} \in \mathbb{R}^N\)  is the target vector that is defined \textit{a priori}, and \(f: \mathbb{R}^M \to \mathbb{R}^N\) represents the objective function which in the case of inverse design is inverted, unlike in forward optimization problems. Inverse design problems are typically ill-posed, meaning that multiple values of $\mathbf{x}$ can yield similar values of $\mathbf{y}$.
As an optimization objective function, the inverse design problem is defined as:

\begin{equation}
\begin{aligned}
& \underset{\mathbf{x}}{\text{minimize}}
& & \epsilon(f(\mathbf{x}), \mathbf{y}) \\
& \text{s. t.}
& & \mathbf{x_{lb}} \leq \mathbf{x} \leq \mathbf{x_{ub}}
\label{eqn:inverse_optimization}
\end{aligned}
\end{equation}

In Eq. (\ref{eqn:inverse_optimization}), \(\epsilon: \mathbb{R}^N \times \mathbb{R}^N \to \mathbb{R}\) is a measure of discrepancy between the desired target vector $\mathbf{y}$ and the forward function $f$ evaluated design vector $\mathbf{x}$. The design vector $\mathbf{x}$ is defined as $\textbf{x} = [x_1,..., x_M]^T$ in the decision space $\mathbb{R}^M$, where $M$ is the dimension of vector. The objective here is to minimize $\epsilon$, and ideally find a perfect match between $f(\mathbf{x})$ and $\mathbf{y}$, i.e., an $\epsilon$ value of 0. Moreover, during the minimization process, the design vector must be within predefined lower and upper boundaries,  $\mathbf{x_{lb}}$ and  $\mathbf{x_{ub}}$, respectively.

For all inverse design benchmarks in this manuscript, the chosen discrepancy measure $\epsilon$ is equivalent to the Root Mean Square Error (RMSE), and mathematically, it is defined as:

\begin{equation}
\begin{aligned}
\text{RMSE} = \sqrt{\frac{1}{N} \sum_{i=1}^{N} (y_{i} - f(\mathbf{x})_{i})^2}
\label{eqn:rmse}
\end{aligned}
\end{equation}

The variables $y_{i}$ and $f(\mathbf{x})_{i}$ in Eq. (\ref{eqn:rmse}) are the i$^{th}$ components of the target vector $\mathbf{y}$ and  $f$ evaluated at the design vector $\mathbf{x}$, respectively. The value of $f(\mathbf{x})$ is a vector itself and has the same number of components $N$ as $\mathbf{y}$.

\subsection{Photonic Surfaces Inverse Design Computational Framework}\label{sec:photonicdesign}

\begin{figure}[!h]
\centering
\begin{subfigure}[b]{0.49\textwidth}
   \raisebox{0.9cm}{\includegraphics[width=\linewidth]{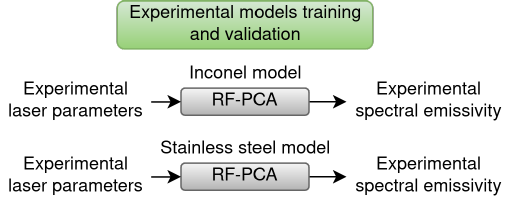}}     \caption{}
   \label{fig:Figure1a}
\end{subfigure}
\begin{subfigure}[b]{0.49\textwidth}
   \includegraphics[width=\linewidth]{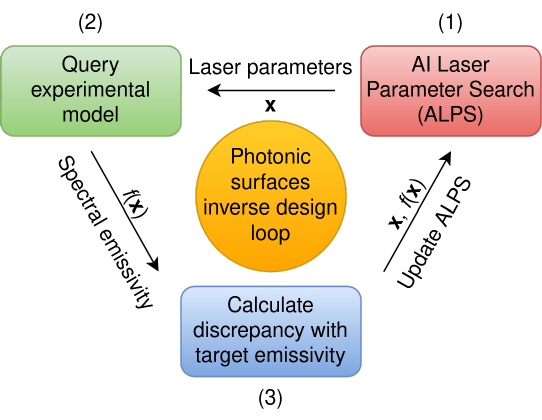}     
   \caption{}
   \label{fig:Figure1b} 
\end{subfigure}
\begin{subfigure}[b]{0.49\textwidth}
   \includegraphics[width=\linewidth]{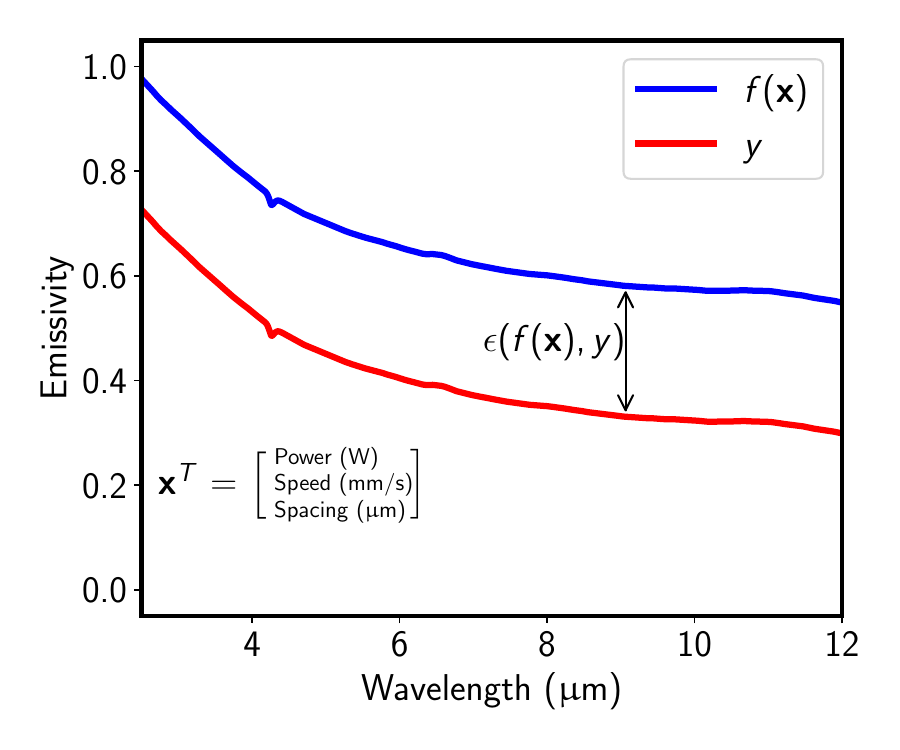}
   \caption{}
   \label{fig:Figure1c}
\end{subfigure}
\begin{subfigure}[b]{0.49\textwidth}
   \includegraphics[width=\linewidth]{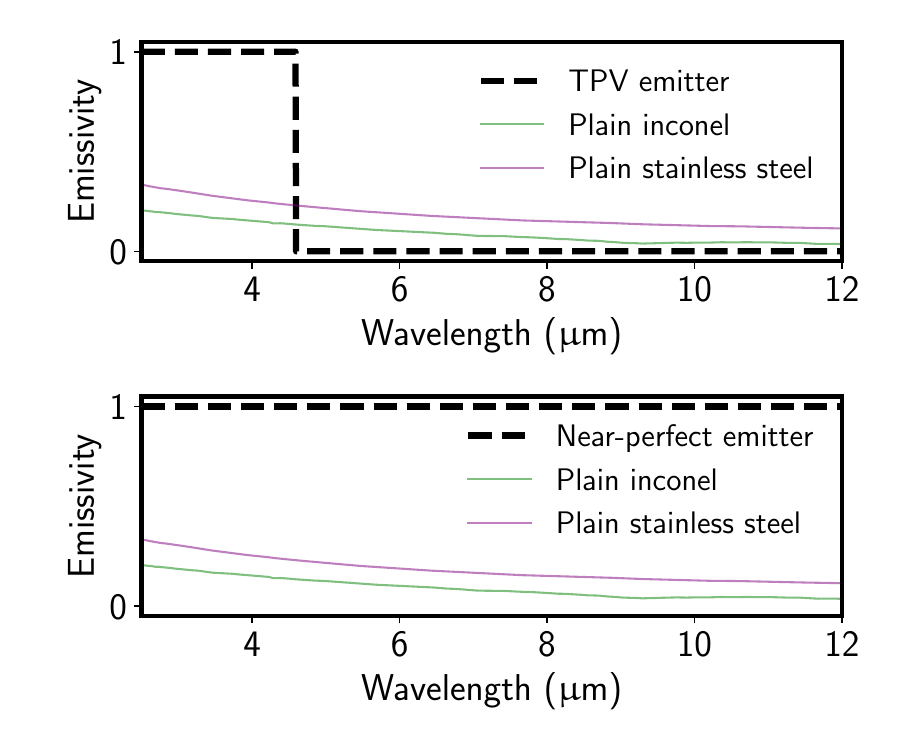}
   \caption{}
   \label{fig:Figure1d}
\end{subfigure}
\caption[]{Photonic surfaces inverse design segments and examples: (a) The ML experimental model pipeline developed to assess laser fabrication parameters. Each model is trained using distinct datasets and categorized according to the surface material employed for laser texturing. The RF algorithm predicts the PCA components, which are then transformed into spectral emissivity curves. (b) The inverse design loop comprises three phases: (1) Generation of laser parameters using the ALPS framework and assessment via an experimental model (defined in (a)) to produce a spectral emissivity curve (2), followed by (3) comparison of this curve against the target spectral emissivity and finally updating  ALPS with new data for iterative decision making. (c) Example of the discrepancy, $\epsilon$, as defined in Eq. (\ref{eqn:inverse_optimization}), between the user-defined target spectral emissivity curve $\mathbf{y}$ and the spectral emissivity curve derived from evaluating the design vector $\mathbf{x}$ with the model defined in (a). (d) Photonic surfaces inverse design benchmark targets: TPV emitter (top) and the near-perfect emitter (bottom). The TPV emitter target switches to 0 emissivity at a wavelength of 4.6 $\upmu$m. The plain spectral emissivity curves are measured from the material without laser texturing.}
\label{fig:Figure1}
\end{figure}

The relationship between laser fabrication parameters and spectral emissivity curves is utilized for the inverse design of photonic surfaces. The laser fabrication parameters are the laser power (W), scanning speed (mm/s), and spacing ($\upmu$m), denoted as $L_p$, $S_s$, and $S$, respectively. The spectral emissivity curves are defined as $N$ dimensional vectors ($N$ = 822) where each component represents an emissivity value (0 to 1) for each of the $N$ wavelength values (interval from 2.5 to 12.5 $\upmu$m). Instead of inferring the spectral emissivity curves from the laser fabrication parameters experimentally, we use ML algorithms to accurately model this relationship (denoted as the experimental model throughout the manuscript). As shown in Fig. \ref{fig:Figure1a}, we utilize a combined RF and Principle Component Analysis (PCA) algorithms (denoted as RF-PCA) as it was shown in our previous work (\citet{grbcic2024ensemble}) that it is accurate and robust. Two experimental models are trained using two distinct datasets, differentiated by the type of plain surface materials used for texturing: Inconel (\citet{grbcic2024ensemble}) and Stainless Steel (\citet{park2024tnn}). Full details of the datasets, RF-PCA hyperparameters, and the experimental model validation procedure can be found in  App. \ref{app:experimental_model}.

Moreover, the photonic surface inverse design loop is shown in Fig. \ref{fig:Figure1b}. The principal aim of this loop is to determine the optimal laser fabrication parameters that yield a specific target spectral emissivity with the fewest possible evaluations of the experimental model. The process commences by generating an initial set of laser fabrication parameters through Latin Hypercube Sampling (LHS), followed by an update of the ALPS framework. Secondly, the ALPS framework explores the laser fabrication parameter space and determines the design vector that should be evaluated by the experimental model to obtain the spectral emissivity value. The target and experimental model-inferred spectral emissivity curves are compared, and based on this discrepancy, the best laser parameters and corresponding spectral emissivity curves are selected to update the ALPS framework to improve its capacity for more precise parameter space exploration. The loop is terminated when the maximum number of experimental model evaluations is reached.

A visual example of the photonic surface inverse design process is shown in Fig. \ref{fig:Figure1c}. The goal is to minimize the discrepancy $\epsilon$ between the two curves as outlined in Eq. (\ref{eqn:inverse_optimization}). The two spectral emissivity curves are the user-defined target spectral emissivity denoted as $\mathbf{y}$, and the spectral emissivity $f(\mathbf{x})$ generated by evaluating the design vector $\mathbf{x}$ ($\textbf{x} = [L_p, S_s, S]^T$) with the experimental model $f$. Due to differences in the material datasets, the lower and upper boundaries $\mathbf{x_{lb}}$ and $\mathbf{x_{ub}}$ vary slightly for each parameter. For both materials, the boundaries are set within the same ranges for $L_p$, $0.2 \, (\text{W}) \leq L_p \leq 1.3 \, (\text{W}) $, and $S_s$, $10 \, (\text{mm/s}) \leq S_s \leq 700 \, (\text{mm/s})$. However, for $S$, the ranges differ: for Inconel, the range is $15 \, (\upmu\text{m}) \leq S \leq 28 \, (\upmu\text{m})$, while for stainless steel, the range is $1 \, (\upmu\text{m}) \leq S \leq 42 \, (\upmu\text{m})$.

Finally, Fig. \ref{fig:Figure1d} shows two spectral emissivity targets that serve as benchmarks to demonstrate the performance of the ALPS framework. The ideal step function spectral emissivity (top curve in Fig. \ref{fig:Figure1d}) represents the optimal emissivity profile that a photonic surface should exhibit for Thermophotovoltaic (TPV) applications, specifically TPV emitters as shown by \citet{park2024tnn}. The bottom curve in Fig. \ref{fig:Figure1d} is the near-perfect emitter where the goal is to determine the laser parameters that yield an emissivity profile that is equal to 1 at all wavelengths (used as a benchmark in \citet{grbcic2024ensemble}). Besides these two photonic surface inverse design benchmarks, we also utilize two synthetic benchmarks to show if our approach can perform well on different inverse design problems. The additional synthetic benchmarks, as well as the algorithms used for comparison with ALPS, and all the numerical experiment setup parameters needed to reproduce this study, are thoroughly described in App. \ref{app:exp_setup}.

\newpage
\subsection{ALPS: AI Laser Parameter Search}\label{subsec:bench}

ALPS plays a main role in the computational framework we have developed for photonic surface inverse design. Its task is to decide on the most informative laser parameters for experimental model evaluation to reconstruct the user-defined target vector through discrepancy minimization defined in Eq. 
(\ref{eqn:inverse_optimization}). The ALPS algorithm is defined in detail in Alg. \ref{alg:alps} with all of the required hyperparameters included.

Firstly, we consider the simpler variant of ALPS where we do not include warm starting (parameter $ws$ is set to False). To start the inverse design process based on the user-defined target vector, we sample the design vector space using LHS, evaluate the samples using the experimental model, and train an initial RF surrogate denoted as $f_{rf}$. The default scikit-learn 1.2.2. hyperparameters of the RF algorithm were used (\citet{pedregosa2011scikit}). Note that we  employed other ML algorithms as the surrogate model in ALPS, but none perform as well as RF. Selecting new design vectors to evaluate with the experimental model is done through a greedy strategy defined as prediction-based exploration by \citet{kochenderfer2019algorithms}, or minimization of an interpolating surface by \citet{jones2001taxonomy}. Moreover, we extend this infill criteria or strategy by selecting a batch of design vectors, instead of a single design vector for experimental model evaluation. 

The mathematical expressions defined in Eq. (\ref{eqn:epsilon}), (\ref{eqn:xbatch}) and (\ref{eqn:selection}), explain the procedure for selecting the best samples for experimental model evaluation. More specifically, given a matrix of samples \( \mathbf{X_s} \) generated using LHS, where each row is a vector $  \mathbf{x}   \in [\mathbf{x_{lb}}, \mathbf{x_{ub}
}] \subset  \mathbb{R}^M$, we compute the discrepancy \( \epsilon(f_{rf}(\mathbf{x}), \mathbf{y}) \) between the target vector $ \mathbf{y} \in \mathbb{R}^N $ and the RF surrogate function $f_{rf}$ evaluated design vector $\mathbf{x}$:
\begin{equation}
\epsilon_i = \epsilon(f_{rf}(\mathbf{x}_i), \mathbf{y}), \quad \forall \mathbf{x}_i \in \mathbf{X_s}.
\label{eqn:epsilon}
\end{equation}

To select the best \( n_{\text{batch}} \) vectors, we define \( \mathbf{X_b} \) as:
\begin{equation}
\mathbf{X_b} = \begin{bmatrix}
\mathbf{x}^{(1)} \\
\mathbf{x}^{(2)} \\
\vdots \\
\mathbf{x}^{(n_{batch})}
\end{bmatrix},
\label{eqn:xbatch}
\end{equation}
where \( \mathbf{x}^{(i)} \) are the top $n_{batch}$ vectors from \( \mathbf{X}_s \) that minimize \( \epsilon \):
\begin{equation}
\mathbf{x}^{(i)} \in \arg\min_{\mathbf{x} \in \mathbf{X}_s} \epsilon(f_{rf}(\mathbf{x}), \mathbf{y}), \quad i = 1, \ldots, n_{batch}.
\label{eqn:selection}
\end{equation}

The $\epsilon$ minimization procedure is done through a simple sorting algorithm. After the selection of $\mathbf{X_b}$, the RF surrogate is retrained with all of the previously sampled design vectors $\mathbf{x}$ as well as the new batch of design vectors $\mathbf{X_b}$, and their respective experimental model evaluations $f(\mathbf{x})$. The process is repeated until the maximum evaluation value $n_{max}$ is reached.

When the warm starting option of ALPS is enabled by setting the boolean parameter $ws$ to True, we utilize a pre-trained model in conjunction with LHS, $f_{ws}$, for generating initial samples, otherwise, we use LHS to generate design values within $\mathbf{x_{lb}}$ and $\mathbf{x_{ub}}$. The model $f_{ws}$ is trained and saved during a previous inverse design process, specifically for a photonic surface design aimed at a different target. We are training the surrogate model using the design vector $\mathbf{x}$ and the corresponding experimental model value  $f(\mathbf{x})$, rather than focusing on the error $\epsilon(f(\mathbf{x}), \mathbf{y}))$. Therefore, \(f_{ws}: \mathbb{R}^M \to \mathbb{R}^N \) is designed to capture the forward relationship between the design vectors and the target space, rather than the inverse design error landscape. This advantage opens up the potential for reusing the surrogate model after each inverse design process is completed to accelerate every subsequent process. If warm starting is enabled, the initial samples are determined with the process described in Eq. (\ref{eqn:epsilon}), (\ref{eqn:xbatch}) and (\ref{eqn:selection}), however, we use $f_{ws}$ instead of $f_{rf}$.

\begin{algorithm}[H]
\caption{AI Laser Parameter Search (ALPS)}
\label{alg:alps}
\begin{algorithmic}[1]
\Require target $\mathbf{y}$, experimental model $f$, batch size $n_{batch}$, initial sample size $n_{init}$,  surrogate sample size $n_s$, maximum evaluations $n_{max}$, lower boundary vector $\mathbf{x_{lb}}$, upper boundary vector $\mathbf{x_{ub}}$, warm start $ws$, warm start model $f_{ws}$
\If{$ws$ is True}  
    \State $\mathbf{X_s} \gets \text{LHS}(n_{s}, \mathbf{x_{lb}}, \mathbf{x_{ub}})$  \Comment{Generate $n_s$ design vectors using LHS bounded by $\mathbf{x_{lb}}$ and $\mathbf{x_{ub}}$ and store in matrix $\mathbf{X_s}$}
    \State $\hat{\mathbf{X}} \gets \mathbf{x}^{(i)} \in \arg\min_{\mathbf{x} \in \mathbf{X}_s} \epsilon(f_{ws} (\mathbf{x}), \mathbf{y}), \quad i = 1, \ldots, n_{init}$
    \Comment{Select $n_{init}$ rows from $\mathbf{X_s}$ based on the smallest $\epsilon$ values and store in $\hat{\mathbf{X}}$}
\Else
    \State $\hat{\mathbf{X}} \gets \text{LHS}(n_{init}, \mathbf{x_{lb}}, \mathbf{x_{ub}})$ \Comment{Generate $n_{init}$ initial design vectors using LHS bounded by $\mathbf{x_{lb}}$ and $\mathbf{x_{ub}}$}
\EndIf
\State $\hat{\mathbf{F}} \gets f(\mathbf{\hat{\mathbf{X}}})$ \Comment{Evaluate each row in matrix $\hat{\mathbf{X}}$ using the experimental model $f$, and store in matrix $\hat{\mathbf{F}}$}
\State $n \gets \text{rows}(\hat{\mathbf{F}})$ \Comment{Number of rows in the matrix $\hat{\mathbf{F}}$}
\State $f_{rf} \gets \text{RF}(\hat{\mathbf{X}}, \hat{\mathbf{F}})$  \Comment{Train the RF surrogate $f_{rf}$ using matrix $\hat{\mathbf{X}}$ and the responses $\hat{\mathbf{F}}$}
\While{$n < n_{max}$}
    \State $\mathbf{X_s} \gets \text{LHS}(n_{s}, \mathbf{x_{lb}}, \mathbf{x_{ub}})$  \Comment{Generate $n_s$ design vectors using LHS bounded by $\mathbf{x_{lb}}$ and $\mathbf{x_{ub}}$ and store in matrix $\mathbf{X_s}$}
    \State $\mathbf{X_b} \gets \mathbf{x}^{(i)} \in \arg\min_{\mathbf{x} \in \mathbf{X}_s} \epsilon(f_{rf} (\mathbf{x}), \mathbf{y}), \quad i = 1, \ldots, n_{batch}$
    \Comment{Select $n_{batch}$ rows from $\mathbf{X_s}$ based on the smallest $\epsilon$ values and store in $\mathbf{X_b}$}
    \State $f(\mathbf{X_b})$ \Comment{Evaluate design vectors in matrix $\mathbf{X_b}$ using the true experimental model $f$}
    \State $\hat{\mathbf{X}} \gets \mathbf{X_b}$ \Comment{Add each row of $\mathbf{X_b}$ into matrix $\hat{\mathbf{X}}$}
    \State $\hat{\mathbf{F}} \gets f(\mathbf{X_b})$ \Comment{Add each response of $f(\mathbf{X_b})$ into matrix $\hat{\mathbf{F}}$}
    \State $f_{rf} \gets \text{RF}(\hat{\mathbf{X}}, \hat{\mathbf{F}})$ \Comment{Retrain $f_{rf}$ with the updated $\hat{\mathbf{X}}$ and the true responses $\hat{\mathbf{F}}$}
    \State $n \gets \text{rows}(\hat{\mathbf{F}})$ \Comment{Update $n$}
\EndWhile
\State $\hat{\epsilon} \gets \epsilon(\mathbf{y}, \hat{\mathbf{F}})$
\Comment{Obtain discrepancy array $\hat{\epsilon}$ based on target $\mathbf{y}$ and all values in $\hat{\mathbf{F}}$ i.e., apply Eq. (\ref{eqn:inverse_optimization}) for each value in $\hat{\mathbf{F}}$}
\State $\mathbf{x_{best}} \gets \text{row corresponding to the smallest value of $\hat{\epsilon}$ in $\hat{\mathbf{X}}$}$ \Comment{Select the best row from $\hat{\mathbf{X}}$ based on the smallest $\hat{\epsilon}$ value}
\end{algorithmic}
\end{algorithm}

\section{Results}\label{sec:results}

In this section, we present the ALPS framework results for the two photonic surface inverse design benchmarks (Fig. \ref{fig:Figure1d}) and two synthetic benchmarks (Fig. \ref{fig:FigureB1b}). ALPS is compared to other optimization algorithms such as Particle Swarm Optimization (PSO), DE, Mesh Adaptive Direct Search (MADS), Nelder Mead (NM), Limited memory Broyden–Fletcher–Goldfarb–Shanno with Boundaries (LBFGSB) , and BO. We also include a simple random sampling algorithm for comparison (denoted as Random). The specifics and hyperparameters of these algorithms are detailed in App. \ref{subapp:algorithms} (Tab. \ref{tab:optimization_algorithms}).

We also showcase the major benefit of using ALPS with warm starting for inverse design by demonstrating its capability of cross-target and cross-material inverse design optimization. For all numerical experiments only 100 experimental model evaluations (denoted as Experiments) are considered, i.e., the $n_{max}=100$. Other ALPS-specific parameters used for all benchmarks and comparisons are $n_{batch}$ = 5, $n_{init}$ = 5, and $n_s$ = 600 (for hyperparameter details see Alg. \ref{alg:alps}). 

To rigorously assess the performance of optimization algorithms across all  benchmarks, we execute each algorithm for 100 repeated runs. Throughout each run, we monitor the error value ($\epsilon$), as defined in Eq. (\ref{eqn:inverse_optimization}), during each of the 100 experimental model evaluations. For these evaluations, we track the "best found so far" value of $\epsilon$, updating this record only if a lower $\epsilon$ value is discovered during subsequent evaluations. This progressive update of the minimum $\epsilon$ forms the basis of our convergence graphs, which illustrate the optimization progress over time. At the end of the trials, we analyze these results by calculating the mean, and the 10$^{th}$ and 90$^{th}$ percentiles, of these minimum values across the 100 repetitions, providing insights into the algorithms' efficiency and performance variability. 

In order to visualize the quality of the solutions generated by each optimization algorithm, we obtain the best design vectors found for each of the 100 repeated runs and we use them as inputs to our benchmark models to obtain the corresponding curves. We average these solutions and calculate the 10$^{th}$ and 90$^{th}$ percentiles and visualize them juxtaposed with the target. We denote these graphs as solution reconstruction graphs. 

\subsection{Synthetic Benchmarks Results}\label{subsec:benchresults1}

\begin{figure}[!h]
\centering
\begin{subfigure}[b]{0.49\textwidth}
   \includegraphics[width=\linewidth]{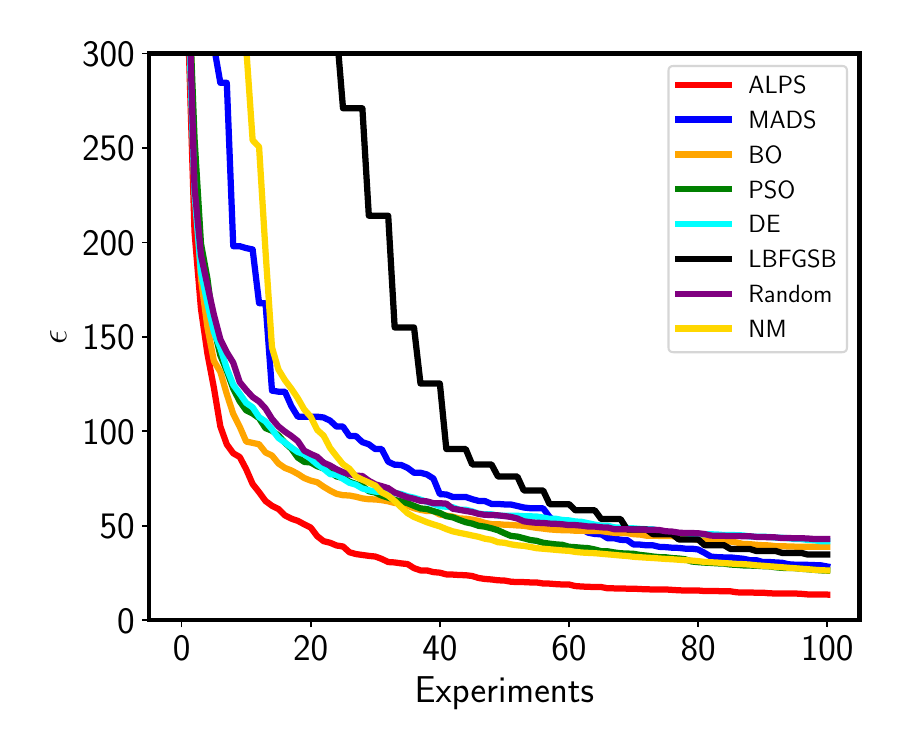}
   \caption{}
   \label{fig:Figure2a}
\end{subfigure}
\begin{subfigure}[b]{0.49\textwidth}
   \includegraphics[width=\linewidth]{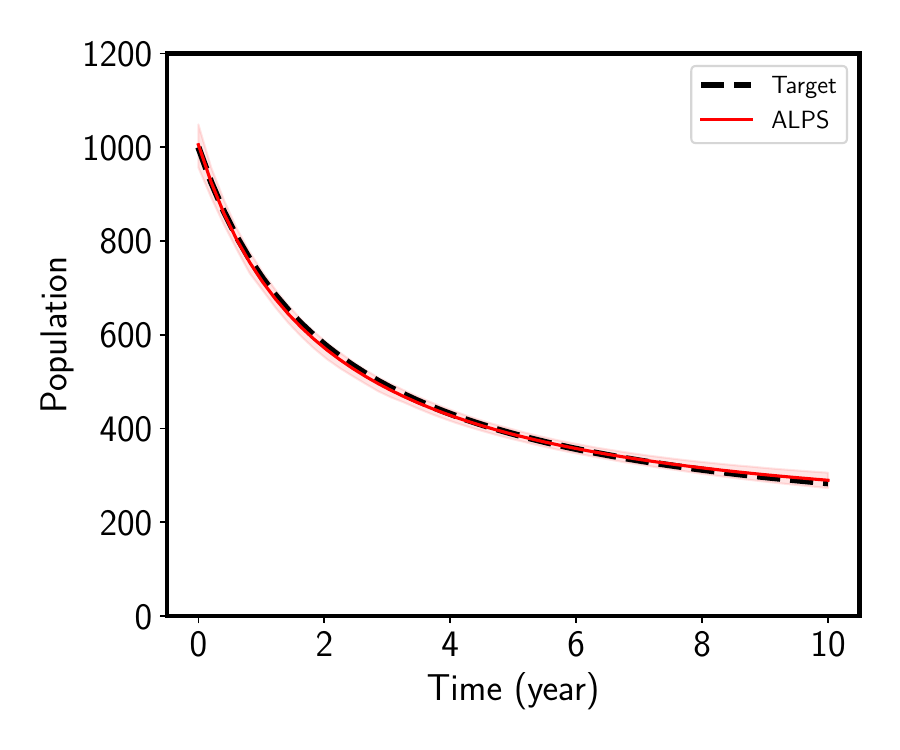}     
   \caption{}
   \label{fig:Figure2b} 
\end{subfigure}
\begin{subfigure}[b]{0.49\textwidth}
   \includegraphics[width=\linewidth]{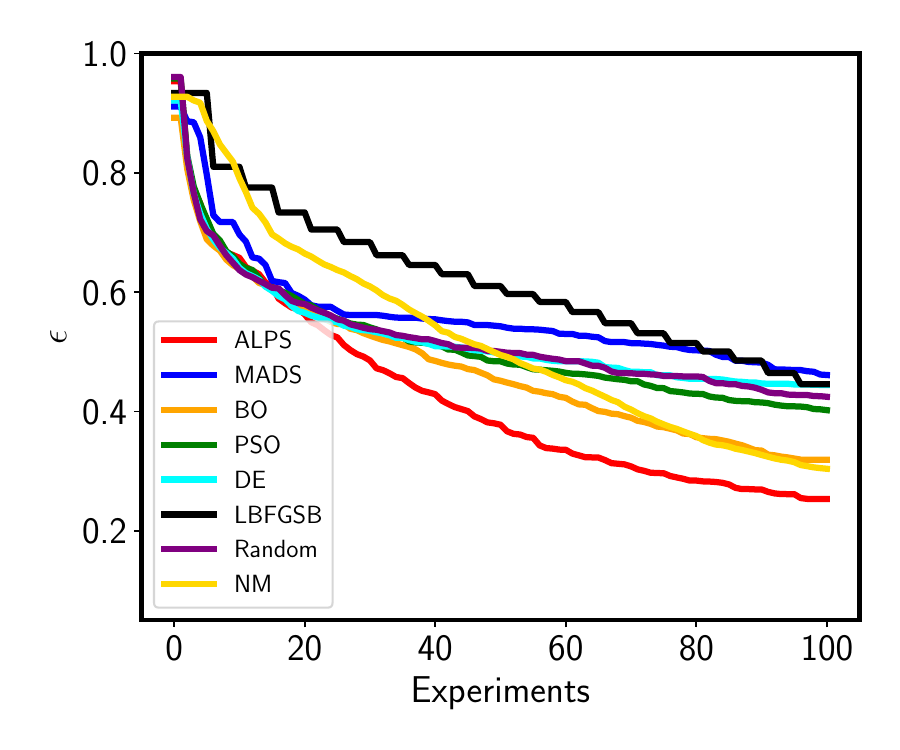}
   \caption{}
   \label{fig:Figure2c}
\end{subfigure}
\begin{subfigure}[b]{0.49\textwidth}
   \includegraphics[width=\linewidth]{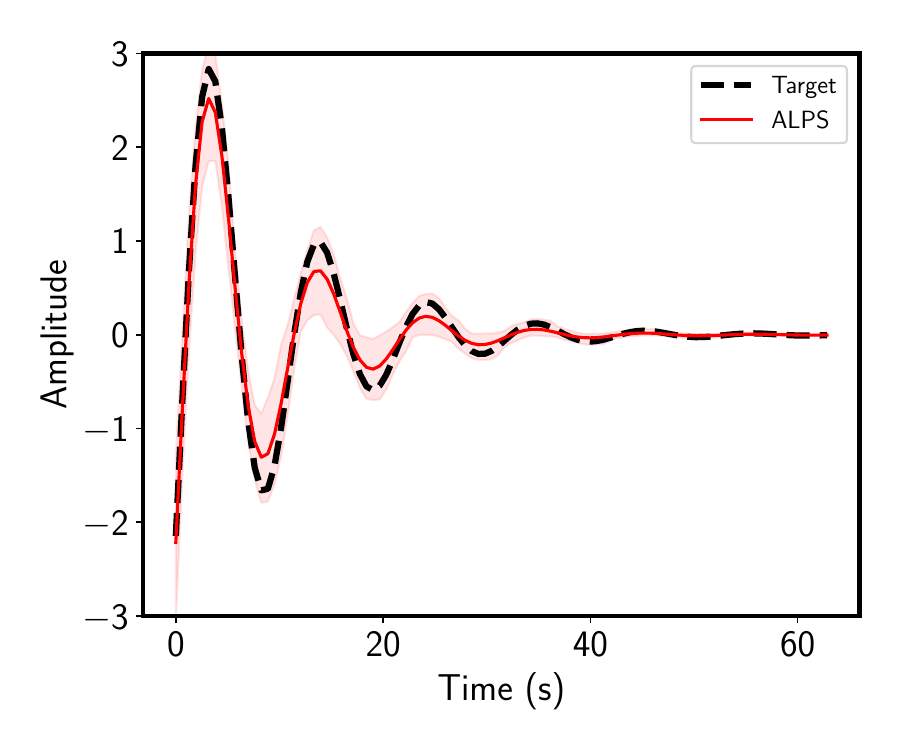}
   \caption{}
   \label{fig:Figure2d}
\end{subfigure}
\caption[]{The inverse design results for the logistic growth (top) and sinusoidal oscillation with damping (bottom) benchmarks: (a) Convergence graphs for all optimization algorithms for the logistic growth benchmark. (b) ALPS solution reconstruction graph for the logistic growth benchmark. (c) Convergence graphs for all optimization algorithms for the sinusoidal oscillation with damping benchmark. (d) ALPS solution reconstruction graph for the sinusoidal oscillation with damping benchmark.}
\label{fig:Figure2}
\end{figure}

In Fig. \ref{fig:Figure2}, we present the results of the three-dimensional logistic growth and four-dimensional sinusoidal oscillation with damping benchmarks. Specifically, Fig. \ref{fig:Figure2a} displays the convergence graphs for all evaluated optimization algorithms, while Fig. \ref{fig:Figure2b} showcases the logistic growth benchmark solution reconstruction graph. It is evident that for the logistic growth benchmark, ALPS solutions cluster tightly around the target value, reflecting the reduced uncertainty in the inverse design. The ALPS mean and standard deviation of the $\epsilon$ value for the logistic growth benchmark after 100 repeated runs are 13.52, and 9.03, respectively. The second best performing algorithm is PSO  with the mean of $\epsilon$ as 26.20  and standard deviation 16.14.

Fig. \ref{fig:Figure2c} and \ref{fig:Figure2d} show the convergence graphs for all algorithms and the ALPS-reconstructed solution graph for the sinusoidal oscillation with damping benchmark, respectively. The ALPS algorithm outperforms all other optimization algorithms; however, NM and BO are the second and third best, with significantly better accuracy than the rest. The reconstructed ALPS solution closely aligns with the target design, but the uncertainty remains slightly higher across 100 repeated runs compared to that of the logistic growth benchmark target. The ALPS mean and standard deviation for $\epsilon$ for the sinusoidal oscillation with damping benchmark are 0.25 and 0.11, while the second and third best, NM and BO, have the vales 0.30 and 0.26, and 0.31 and 0.13, respectively. 

Detailed statistics from all repeated runs are provided in Tab. \ref{tab:detailed_results_logistic} and \ref{tab:detailed_results_sinusoidal}, which cover logistic growth and sinusoidal oscillation with damping, respectively. For the logistic growth benchmark, ALPS demonstrates superior performance in terms of both accuracy and reliability, as it exhibits the best mean and standard deviation across the 100 repeated runs. However, the minimal error is achieved by NM, which, despite obtaining the lowest error, exhibits the largest standard deviation, indicating no reliability. In the case of sinusoidal oscillation with damping, ALPS also excels by showing the lowest mean error and one of the lowest standard deviations in the final model evaluation. Although NM finds the best solution, it proves to be the least reliable, heavily dependent on the initial random design point.

Further in-depth convergence graphs and solution reconstruction graphs of both benchmarks are available in App. \ref{subapp:synbencresults}. Figure \ref{fig:FigureC1} provides convergence graphs for all algorithms across both benchmarks, incorporating measures of uncertainty. Fig. \ref{fig:FigureC2} and \ref{fig:FigureC3} present the solution reconstruction graphs for the logistic growth and sinusoidal oscillation with damping benchmarks, respectively. 

\newpage
\subsection{Photonic Surfaces Design Results}\label{subsec:phodesignresults}

\begin{figure}[!h]
\centering
\begin{subfigure}[b]{0.49\textwidth}
   \includegraphics[width=\linewidth]{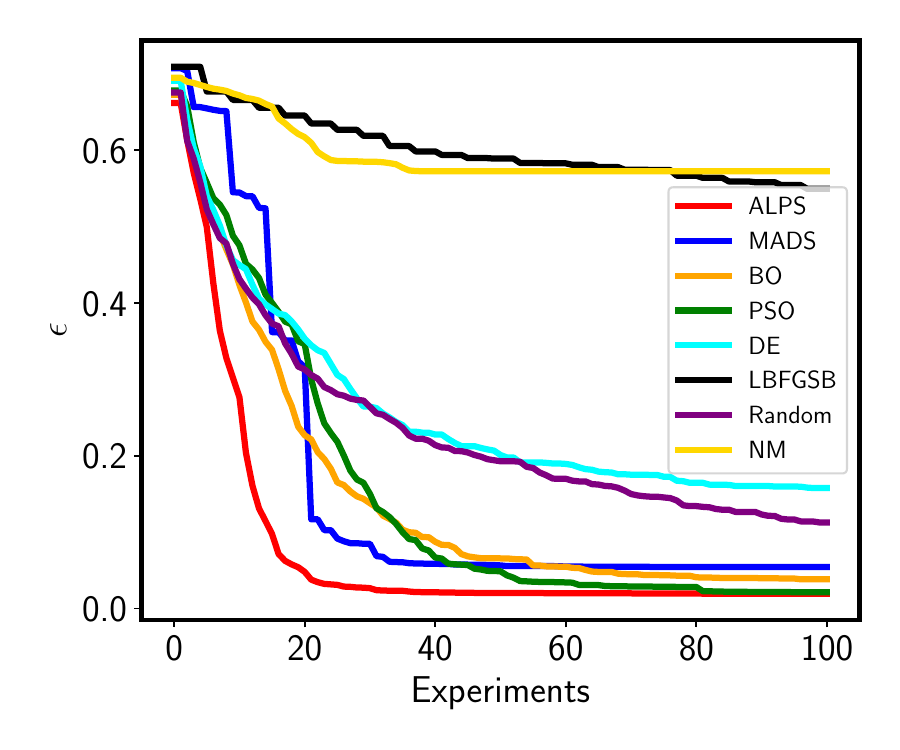}
   \caption{}
   \label{fig:Figure3a}
\end{subfigure}
\begin{subfigure}[b]{0.49\textwidth}
   \includegraphics[width=\linewidth]{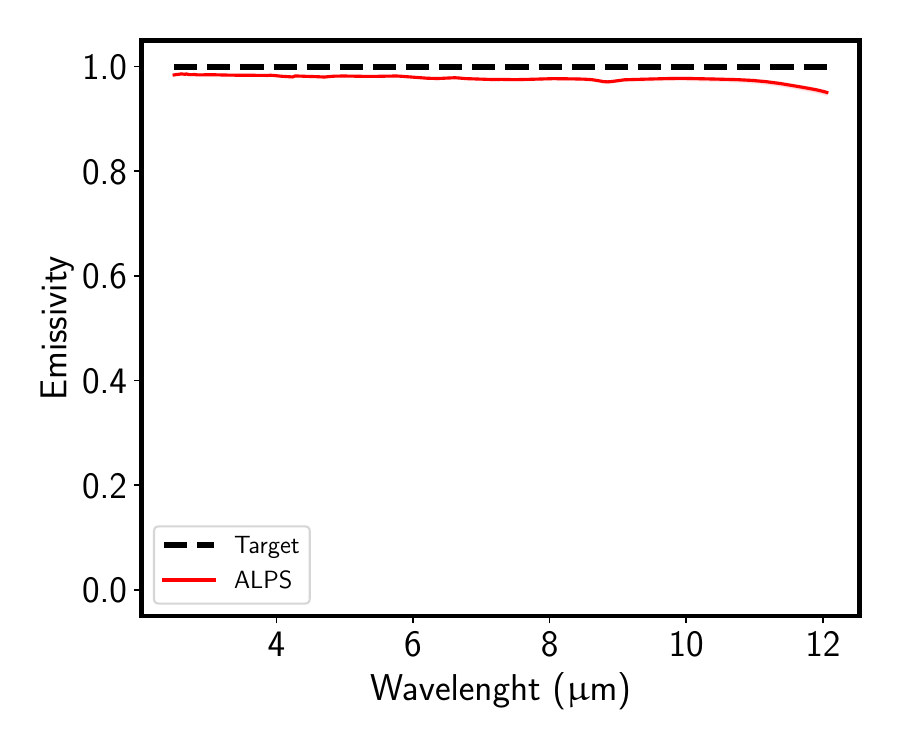}     
   \caption{}
   \label{fig:Figure3e} 
\end{subfigure}
\begin{subfigure}[b]{0.49\textwidth}
   \includegraphics[width=\linewidth]{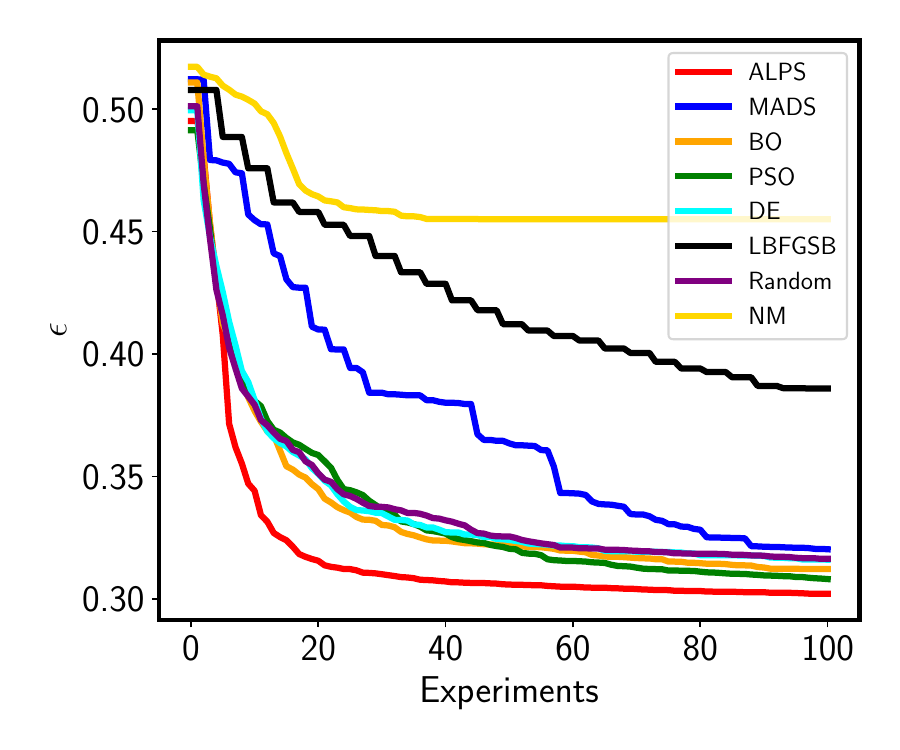}
   \caption{}
   \label{fig:Figure3b}
\end{subfigure}
\begin{subfigure}[b]{0.49\textwidth}
   \includegraphics[width=\linewidth]{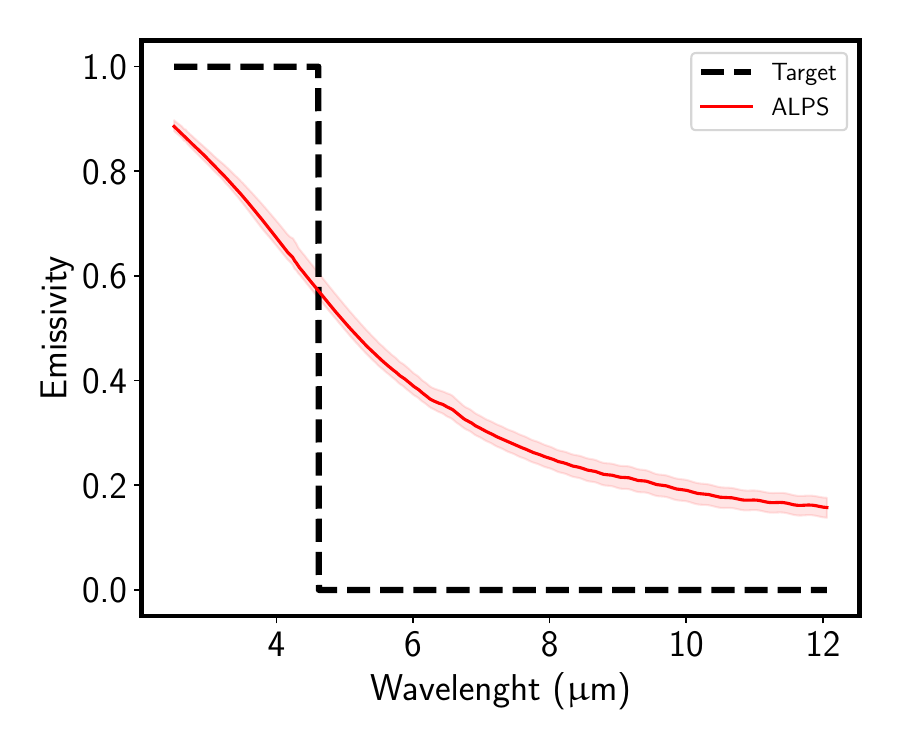}
   \caption{}
   \label{fig:Figure3f}
\end{subfigure}
\caption[]{The inverse design results for the Inconel photonic surface  benchmarks. Convergence graphs for all optimization algorithms are shown in the first column, while in the second, ALPS solution reconstruction graphs are shown: (a) Convergence graphs for the Inconel near-perfect emitter target benchmark. (b) ALPS solution reconstruction graph for the Inconel near-perfect emitter benchmark. (c) Convergence graphs for the Inconel TPV emitter target benchmark. (d) ALPS solution reconstruction graph for the Inconel TPV emitter.}
\label{fig:Figure3_1}
\end{figure}

\begin{figure}[!h]
\centering
\begin{subfigure}[b]{0.49\textwidth}
   \includegraphics[width=\linewidth]{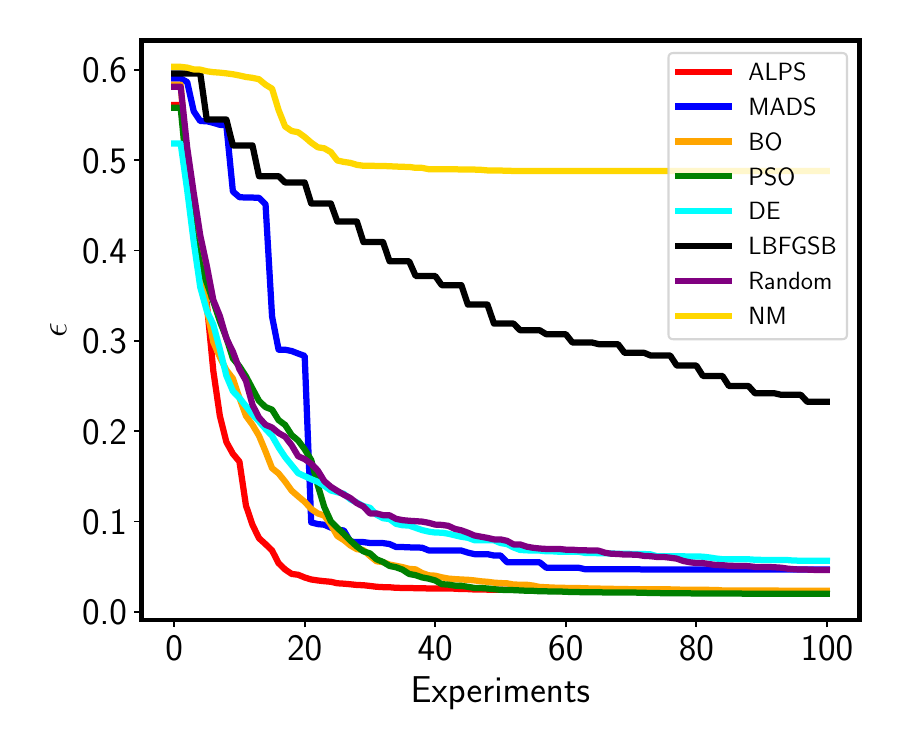}
   \caption{}
   \label{fig:Figure3c}
\end{subfigure}
\begin{subfigure}[b]{0.49\textwidth}
   \includegraphics[width=\linewidth]{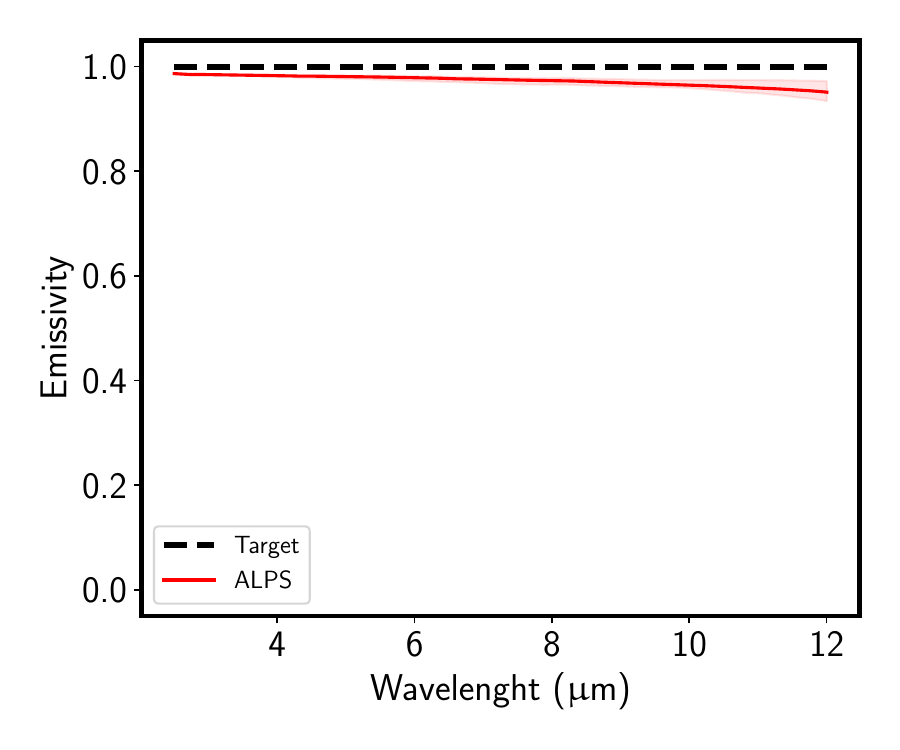}     
   \caption{}
   \label{fig:Figure3g} 
\end{subfigure}
\begin{subfigure}[b]{0.49\textwidth}
   \includegraphics[width=\linewidth]{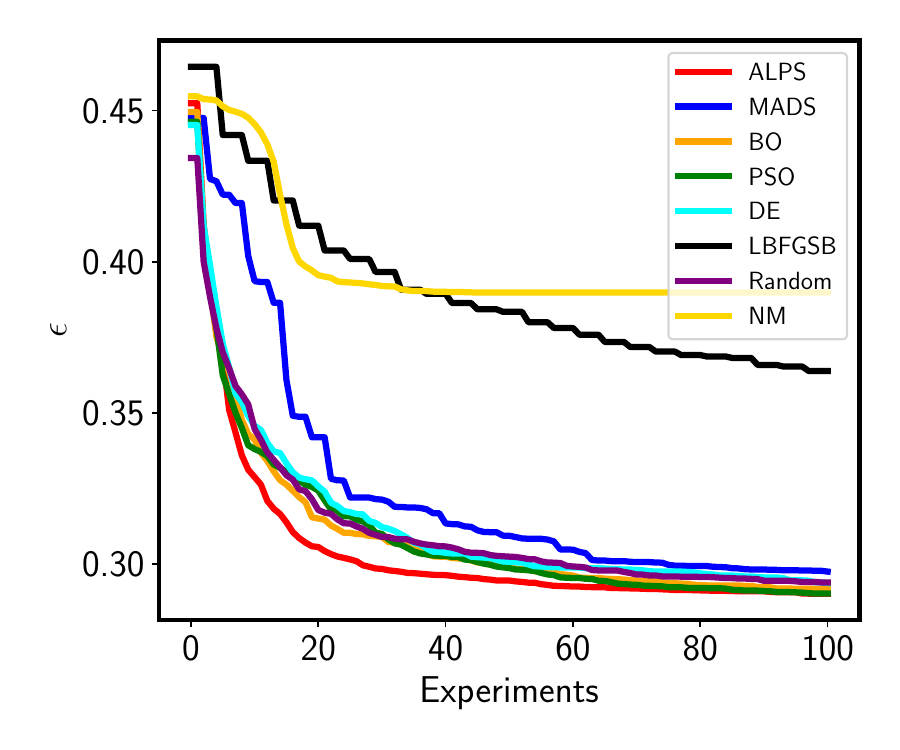}
   \caption{}
   \label{fig:Figure3d}
\end{subfigure}
\begin{subfigure}[b]{0.49\textwidth}
   \includegraphics[width=\linewidth]{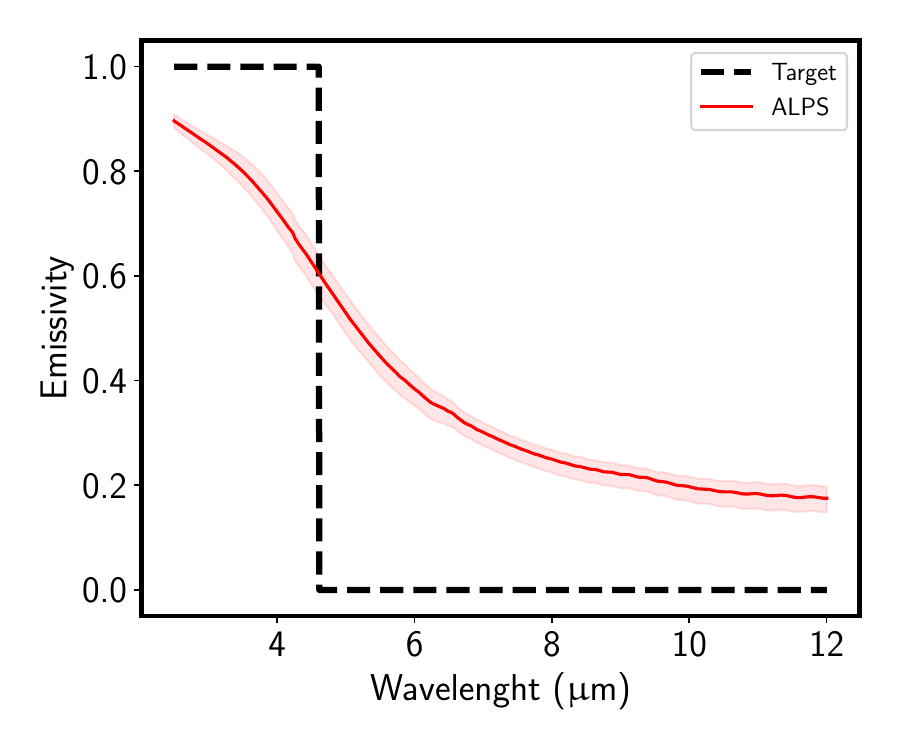}
   \caption{}
   \label{fig:Figure3h}
\end{subfigure}
\caption[]{The inverse design results for the Stainless steel photonic surface benchmarks. Convergence graphs for all optimization algorithms are shown in the first column, while in the second, ALPS solution reconstruction graphs are shown: (a) Convergence graphs for the Stainless steel near-perfect emitter target benchmark. (b) ALPS solution reconstruction graph for the Stainless steel near-perfect emitter benchmark. (c) Convergence graphs for the Stainless steel TPV emitter target benchmark. (d) ALPS solution reconstruction graph for the Stainless steel TPV emitter.}
\label{fig:Figure3_2}
\end{figure}

Fig. \ref{fig:Figure3_1} and Fig. \ref{fig:Figure3_2} show the results of the photonic surface inverse design benchmarks. The first column (Fig. \ref{fig:Figure3a} and \ref{fig:Figure3b}, and Fig. \ref{fig:Figure3c} and \ref{fig:Figure3d}) shows the convergence graphs of the Inconel near-perfect emitter, Inconel TPV emitter, Stainless steel near-perfect emitter, and Stainless steel TPV emitter, respectively. The second column of both figures (Fig. \ref{fig:Figure3e} and \ref{fig:Figure3f}, and Fig. \ref{fig:Figure3g} and \ref{fig:Figure3h}) shows, in the same order, the ALPS solution reconstruction graphs as it is the best performing algorithm for all of the target values. Moreover, in contrast to the synthetic benchmarks, the performance of NM on the photonic surface inverse design is poor. Convergence graphs with uncertainty can be found in the App. \ref{subapp:photonicresults},  Fig. \ref{fig:FigureC4}. The near-perfect emitter reconstructed solutions are extremely close to the target values (for both Inconel and Stainless steel benchmarks), however, the TPV emitter solution reconstruction shows that it is not possible to fully approximate the target. This approximation could be potentially improved by including additional design space parameters, i.e., additional laser fabrication parameters.

Compared to other optimization methods, ALPS exhibits quicker convergence. For a real experimental setting, this is extremely beneficial as the main goal is to achieve inverse design with minimal usage of experimental resources. Most optimization algorithms are able to converge to satisfying solutions (shown in Fig. \ref{fig:FigureC5}-\ref{fig:FigureC8})  (with the exception of NM, and LBFGSB, that are highly sensitive to the initial condition). However, similarly as for the synthetic benchmarks, ALPS exhibits reliability and consistency as it has an overall lowest uncertainty across all benchmarks (see detailed convergence statistics for all runs  presented in Tab. \ref{tab:detailed_results_inconel_near_perfect}-\ref{tab:detailed_results_ss_TPV}). Generally, this property of the optimization algorithm is ideal when only a low number of evaluations are allowed and multiple repeated trials are not possible. 

For the Inconel near-perfect emitter target, the ALPS mean $\epsilon$ is 0.02, which is the same as PSO, however, the standard deviation of $\epsilon$ for ALPS is 0.003, while it is 0.01 for PSO. For the Inconel TPV emitter  target, the ALPS mean and standard deviation of $\epsilon$ are 0.30 and 0.003, respectively,  while the second and third best are the PSO and BO algorithms, with 0.31 and 0.01, and 0.31 and 0.02 mean and standard deviation for $\epsilon$, respectively. For the Stainless steel near-perfect emitter, ALPS mean and standard deviation of $\epsilon$  are 0.02 and 0.004 which is the same as PSO. BO is the second best performing algorithm for this case having achieved the same mean, but a larger standard deviation of 0.02. 

The final target of the photonic surface benchmark, the Stainless steel TPV emitter, is not a hard challenge for most optimization algorithms, as ALPS, BO, PSO, DE and even random sampling have the same $\epsilon$ value, and are the same in all other metrics, except for the maximum $\epsilon$ value obtained in the 100 repeated runs, where the lowest is found by ALPS with 0.30. Finally, the mean $\epsilon$ values for the top three optimization algorithms (ALPS, BO, PSO) are comparable after 100 experimental model evaluations for all benchmarks. However, a significant advantage of ALPS is its ability to accurately infer laser parameters with as few as 20 to 50 experimental model evaluations, as demonstrated in the first column of Fig. \ref{fig:Figure3_1} and Fig. \ref{fig:Figure3_2}, i.e., ALPS finds better solutions faster than all other algorithms. For all these analyses we used a set of unoptimized ALPS and out-of-the-box RF algorithm hyperparameters, however, a detailed analysis on how these parameters influence the performance of ALPS is given in App. \ref{subapp:alps_hyperparameters}.

\newpage
\subsection{Photonic Surfaces Design with Warm Starting}\label{subsec:phodesignresults}

We employ two strategies to accelerate the convergence of ALPS through warm starting. The first strategy, termed cross-target warm starting, uses a model derived from a previous inverse design process with a different target, more specifically, the same material is used but a different target shape. The second, cross-material warm starting, utilizes a model from a prior process where both the target shape and the material differed. As outlined in Alg. \ref{alg:alps}, these pre-trained models are used for generating initial samples instead of the LHS method. These samples are then used to train an initial RF surrogate which is further used in the main ALPS loop.

To demonstrate the effectiveness and versatility of ALPS with both cross-target and cross-material warm starting, we conduct a comparative analysis against ALPS without warm starting using our photonic surface inverse design benchmark. We perform 100 runs for each target to generate the convergence graphs. Each run limits the experimental model evaluations to 100. Additionally, prior to running of ALPS with warm starting, we run the ALPS without warm starting on the different target to generate the model that is used for initializing samples. This generated model is subsequently used for warm starting each ALPS run. This procedure is included in every of the 100 repeated runs for each target to mitigate any variability in the performance of the prior warm starting models.

\begin{figure}[!h]
\centering
\begin{subfigure}[b]{0.49\textwidth}
   \includegraphics[width=\linewidth]{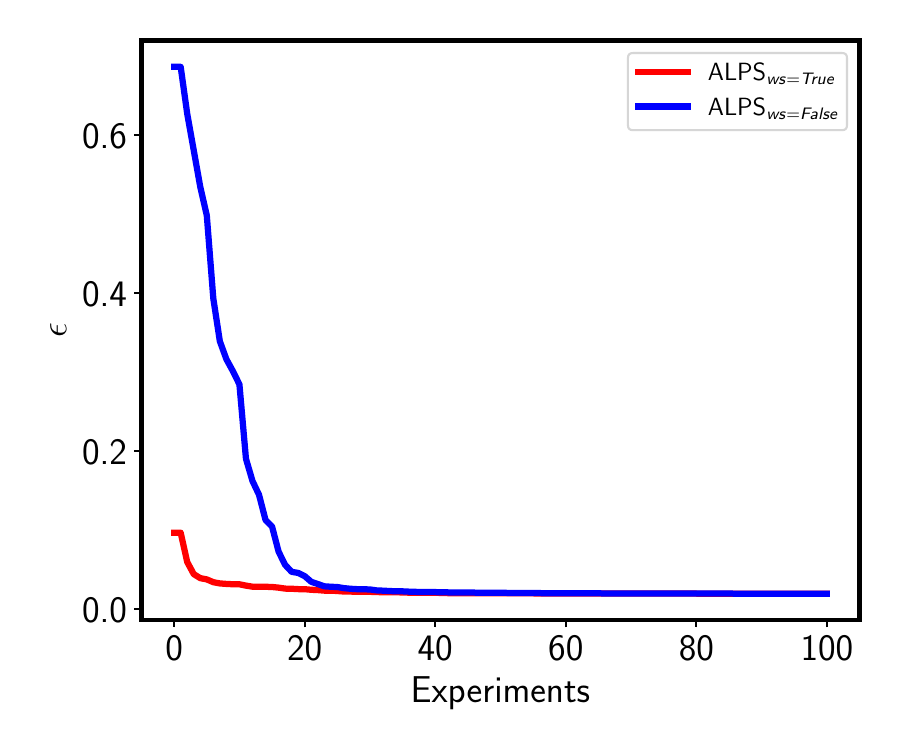}
   \caption{}
   \label{fig:Figure4a}
\end{subfigure}
\begin{subfigure}[b]{0.49\textwidth}
   \includegraphics[width=\linewidth]{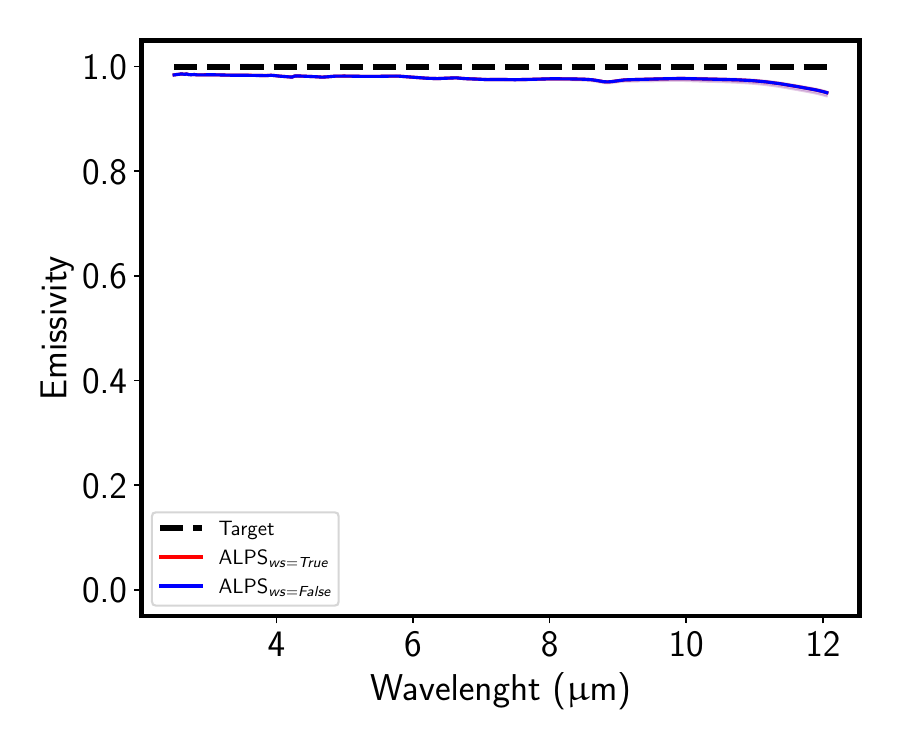}     
   \caption{}
   \label{fig:Figure4e} 
\end{subfigure}
\begin{subfigure}[b]{0.49\textwidth}
   \includegraphics[width=\linewidth]{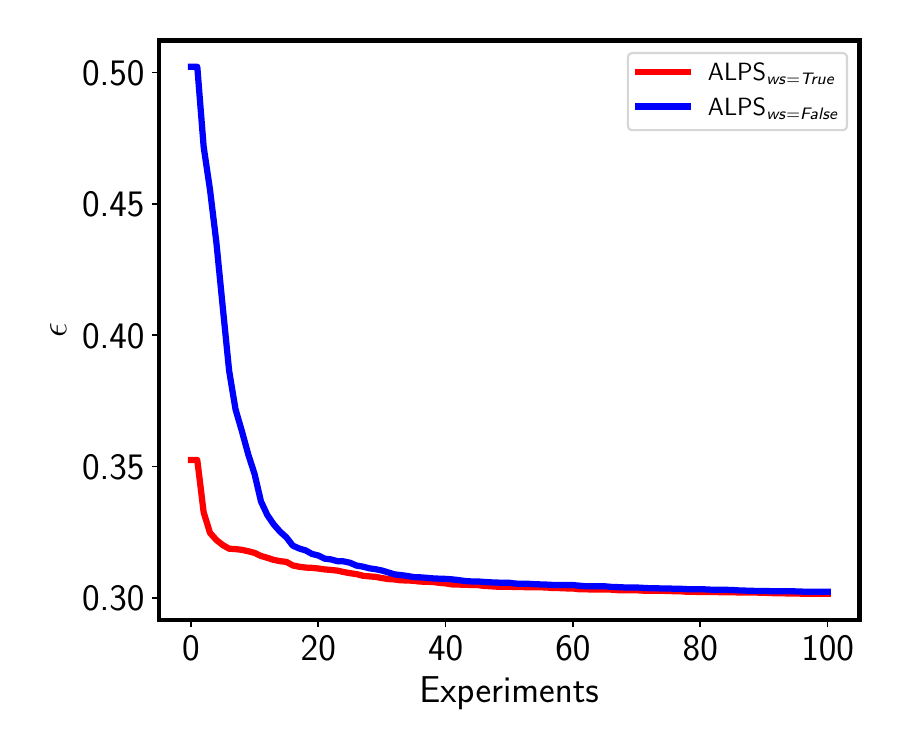}
   \caption{}
   \label{fig:Figure4b}
\end{subfigure}
\begin{subfigure}[b]{0.49\textwidth}
   \includegraphics[width=\linewidth]{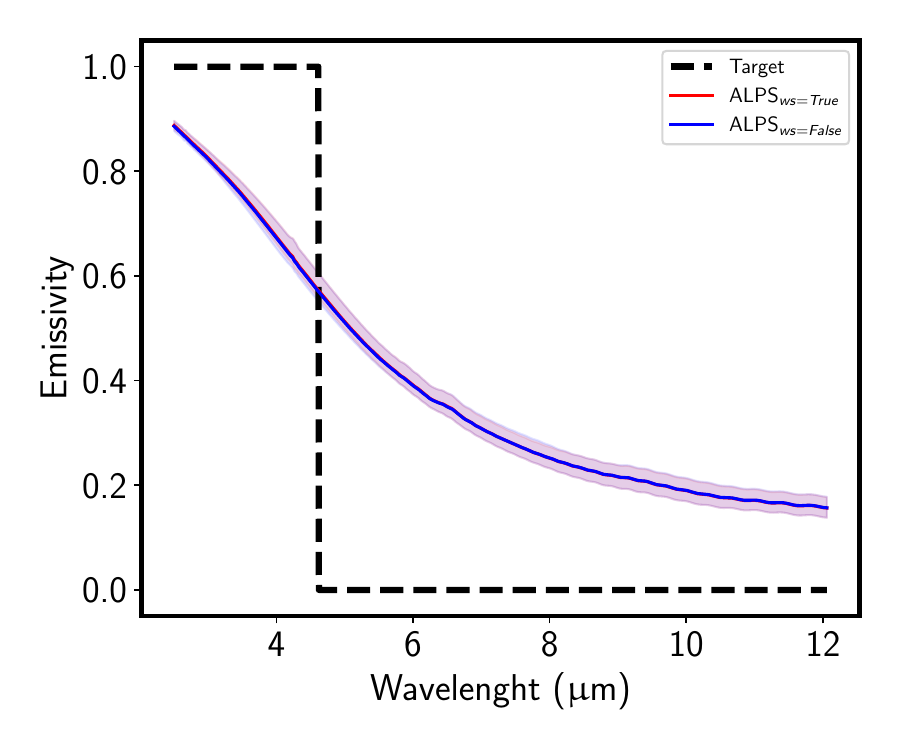}
   \caption{}
   \label{fig:Figure4f}
\end{subfigure}
\caption[]{The inverse design results of ALPS with and without cross-target warm starting for the photonic surface benchmarks (ALPS$_{ws=True}$ means ALPS with warm starting and is shown as the red line). Convergence graphs are shown in the first column, while in the second, the solution reconstruction graphs are shown: (a) Convergence graphs for the Inconel near-perfect emitter target benchmark. (b) Solution reconstruction graph for the Inconel near-perfect emitter benchmark. (c) Convergence graphs for the Inconel TPV emitter target benchmark. (d) Solution reconstruction graph for the Inconel TPV emitter. }
\label{fig:Figure4_1}
\end{figure}

\begin{figure}[!h]
\centering
\begin{subfigure}[b]{0.49\textwidth}
   \includegraphics[width=\linewidth]{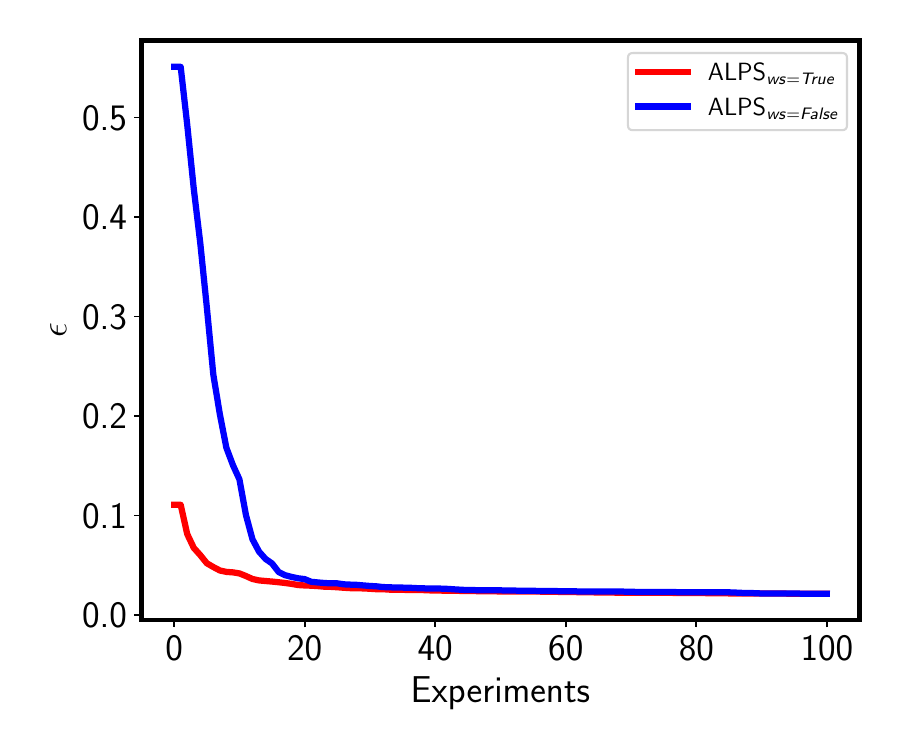}
   \caption{}
   \label{fig:Figure4c}
\end{subfigure}
\begin{subfigure}[b]{0.49\textwidth}
   \includegraphics[width=\linewidth]{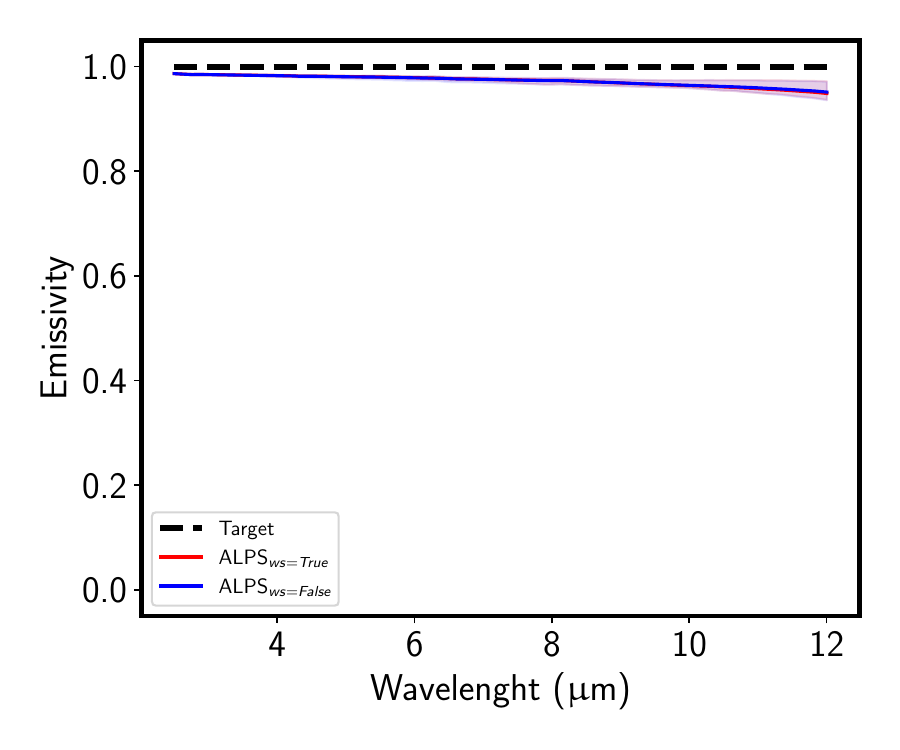}     
   \caption{}
   \label{fig:Figure4g} 
\end{subfigure}
\begin{subfigure}[b]{0.49\textwidth}
   \includegraphics[width=\linewidth]{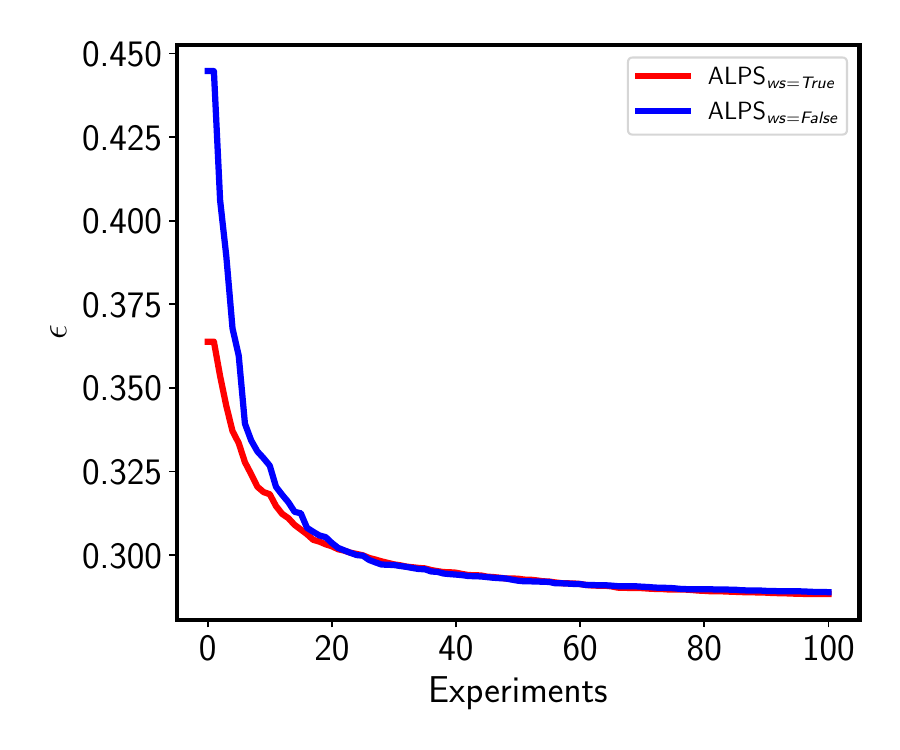}
   \caption{}
   \label{fig:Figure4d}
\end{subfigure}
\begin{subfigure}[b]{0.49\textwidth}
   \includegraphics[width=\linewidth]{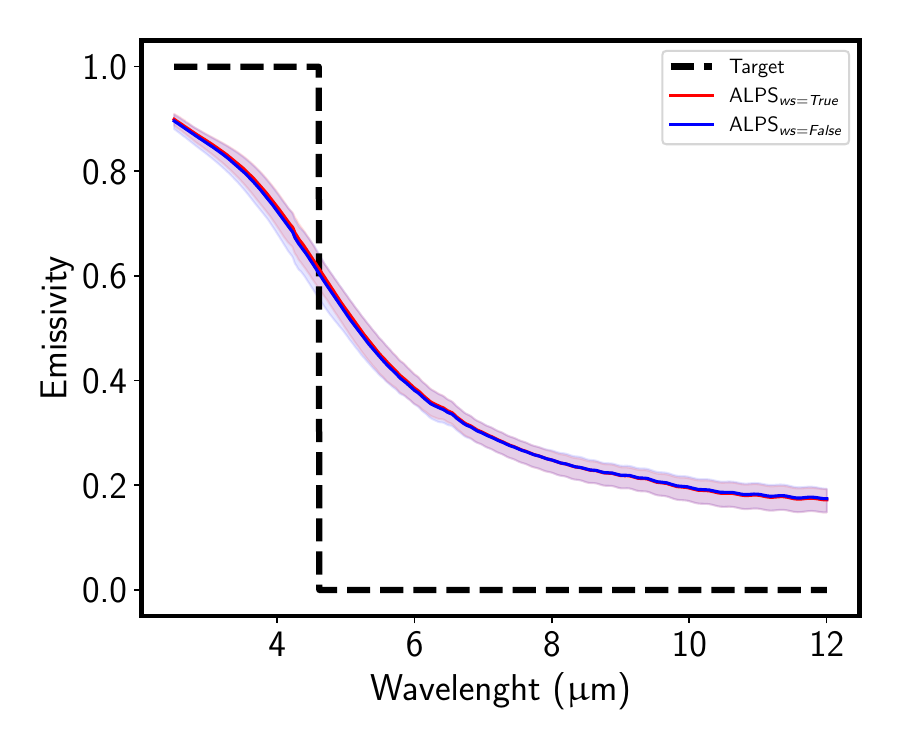}
   \caption{}
   \label{fig:Figure4h}
\end{subfigure}
\caption[]{The inverse design results of ALPS with and without cross-target warm starting for the photonic surface benchmarks (ALPS$_{ws=True}$ means ALPS with warm starting and is shown as the red line). Convergence graphs are shown in the first column, while in the second, the solution reconstruction graphs are shown: (a) Convergence graphs for the Stainless steel near-perfect emitter target benchmark. (b) Solution reconstruction graph for the Stainless steel near-perfect emitter. (c) Convergence graphs for the Stainless steel TPV emitter target benchmark. (d) Solution reconstruction graph for the Stainless steel TPV emitter.}
\label{fig:Figure4_2}
\end{figure}

Fig. \ref{fig:Figure4_1} and Fig. \ref{fig:Figure4_2} show the convergence graphs and the solution reconstruction graphs of both ALPS with and without cross-target warm starting for Inconel and Stainless steel, respecitvely. As an example, when the Inconel TPV emitter is used as the target, ALPS is warm started by a model that is generated after an inverse design of the Inconel near-perfect emitter, and vice-versa. From the convergence graphs (Fig. \ref{fig:Figure4a} and \ref{fig:Figure4c}, and Fig. \ref{fig:Figure4b} and \ref{fig:Figure4d}) it can be seen that for all cases the warm started variant of ALPS (denoted as ALPS$_{ws=True}$) performs better than the ALPS without warm starting. When the targets are Inconel near-perfect emitter, Inconel TPV emitter and Stainless steel near-perfect emitter (top rows of Fig. \ref{fig:Figure4_1} and \ref{fig:Figure4_2}), ALPS with warm starting can approximate the target in 10 to 20 experimental model evaluations. For the Stainless steel TPV emitter target, the ALPS with warm starting is slightly better than the ALPS without warm starting. However, after 100 experimental model evaluations, all variants converge to the approximately the same design. Detailed convergence statistics of all runs are available in App. \ref{subapp:alpswarmstarting}, more specifically, Tab. \ref{tab:detailed_results_cross_target_inconel_100}-\ref{tab:detailed_results_cross_material_ss_step}. 

\begin{figure}[!h]
\centering
\begin{subfigure}[b]{0.49\textwidth}
   \includegraphics[width=\linewidth]{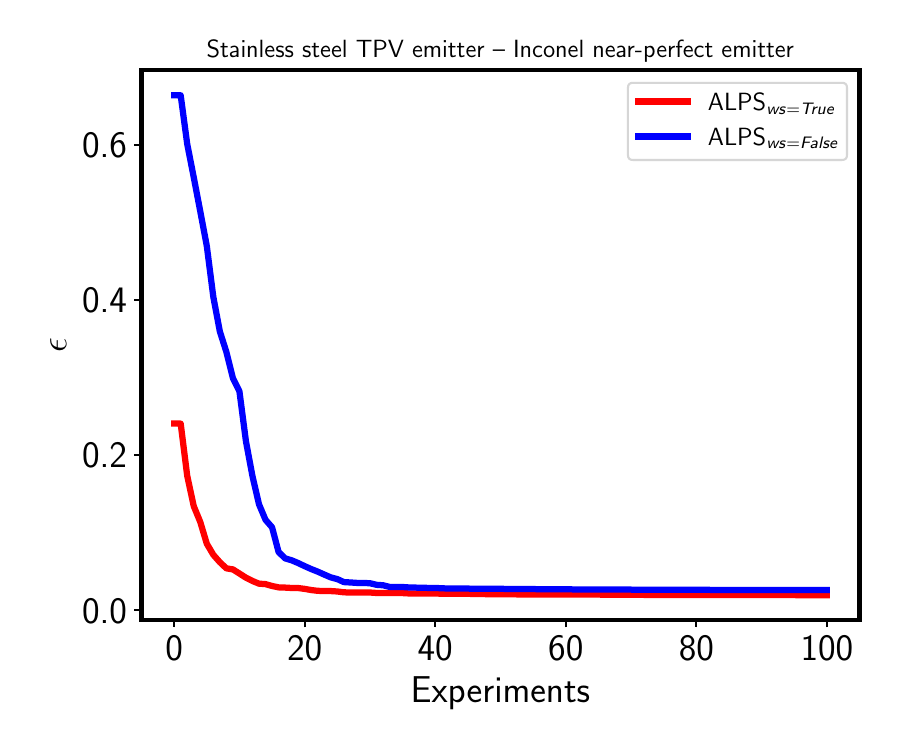}
   \caption{}
   \label{fig:Figure5a}
\end{subfigure}
\begin{subfigure}[b]{0.49\textwidth}
   \includegraphics[width=\linewidth]{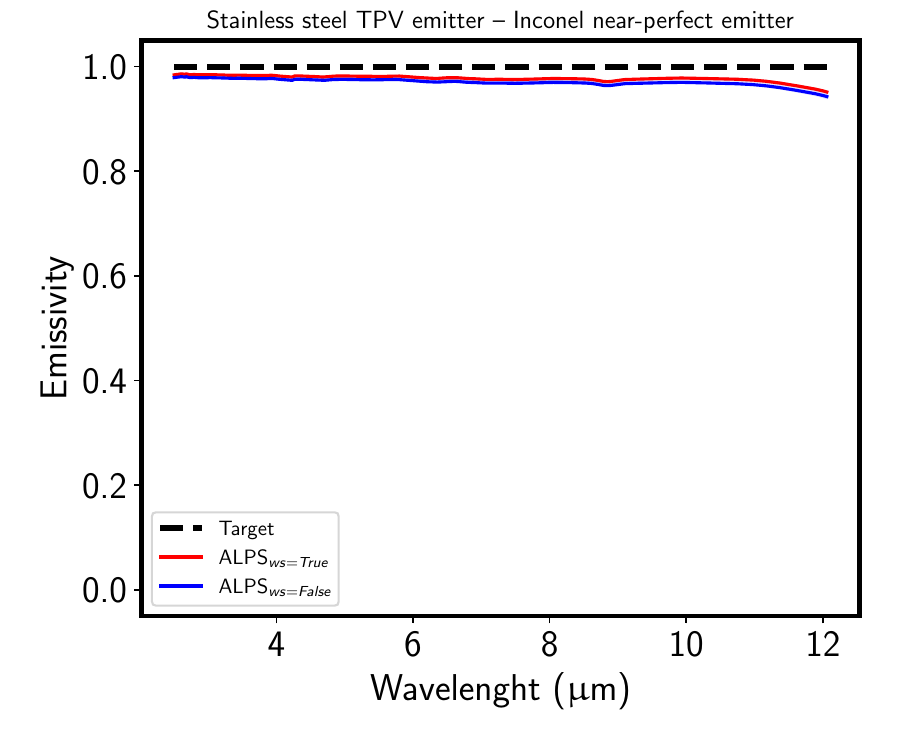}     
   \caption{}
   \label{fig:Figure5e} 
\end{subfigure}
\begin{subfigure}[b]{0.49\textwidth}
   \includegraphics[width=\linewidth]{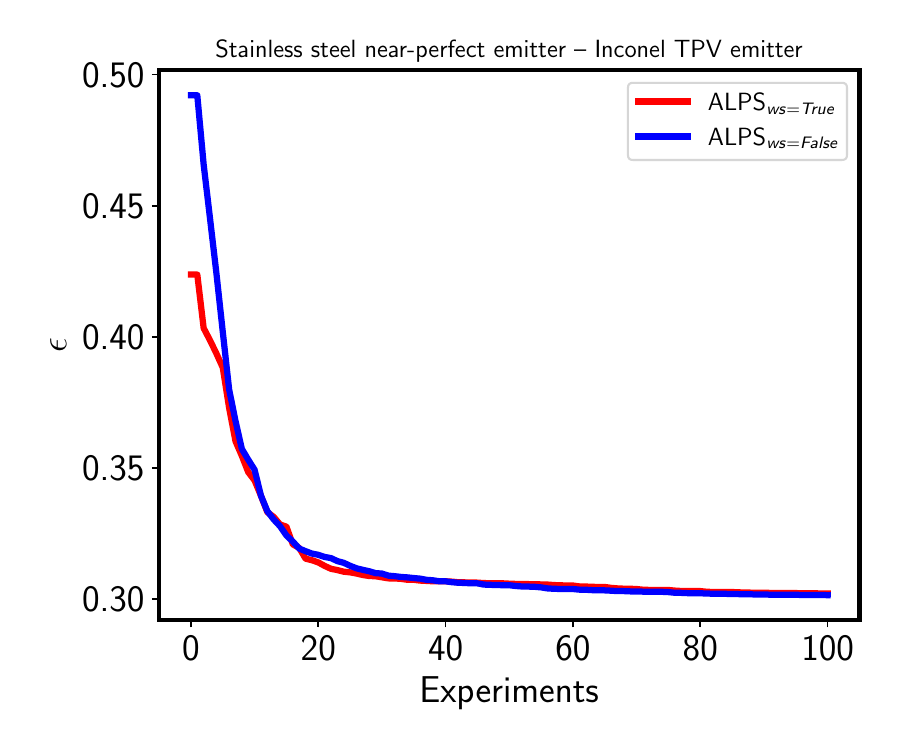}
   \caption{}
   \label{fig:Figure5b}
\end{subfigure}
\begin{subfigure}[b]{0.49\textwidth}
   \includegraphics[width=\linewidth]{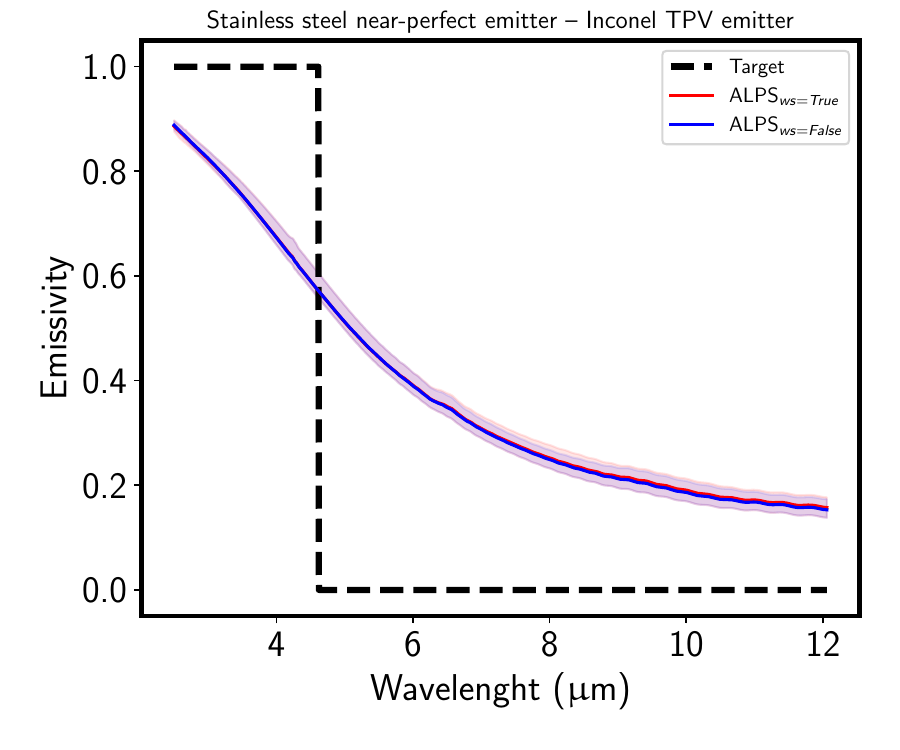}
   \caption{}
   \label{fig:Figure5f}
\end{subfigure}
\caption[]{The inverse design results of ALPS with and without cross-material and cross-target warm starting for the photonic surface benchmarks (ALPS$_{ws=True}$ means ALPS with warm starting and is shown as the red line). The title of each graph presents the warm starting model source (stated first) and inverse design target (stated last). Convergence graphs are shown in the first column, while in the second, the solution reconstruction graphs are shown: (a) Convergence graphs for the Inconel near-perfect emitter target benchmark. (b) Solution reconstruction graph for the Inconel near-perfect emitter benchmark. (c) Convergence graphs for the Inconel TPV emitter target benchmark. (d) Solution reconstruction graph for the Inconel TPV emitter.}
\label{fig:Figure5_1}
\end{figure}

\begin{figure}[!h]
\centering
\begin{subfigure}[b]{0.49\textwidth}
   \includegraphics[width=\linewidth]{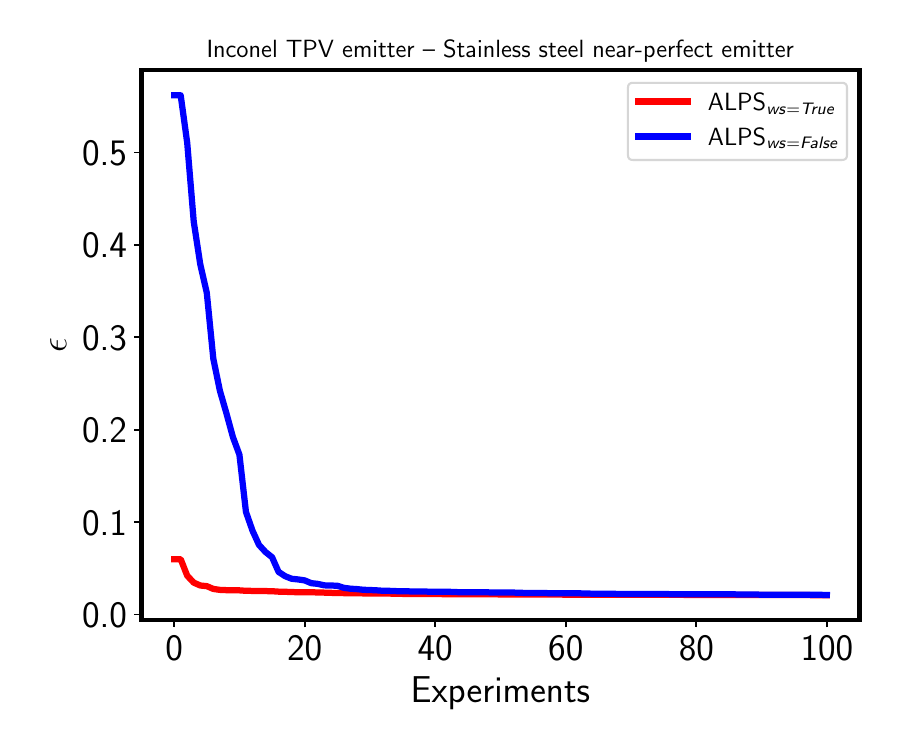}
   \caption{}
   \label{fig:Figure5c}
\end{subfigure}
\begin{subfigure}[b]{0.49\textwidth}
   \includegraphics[width=\linewidth]{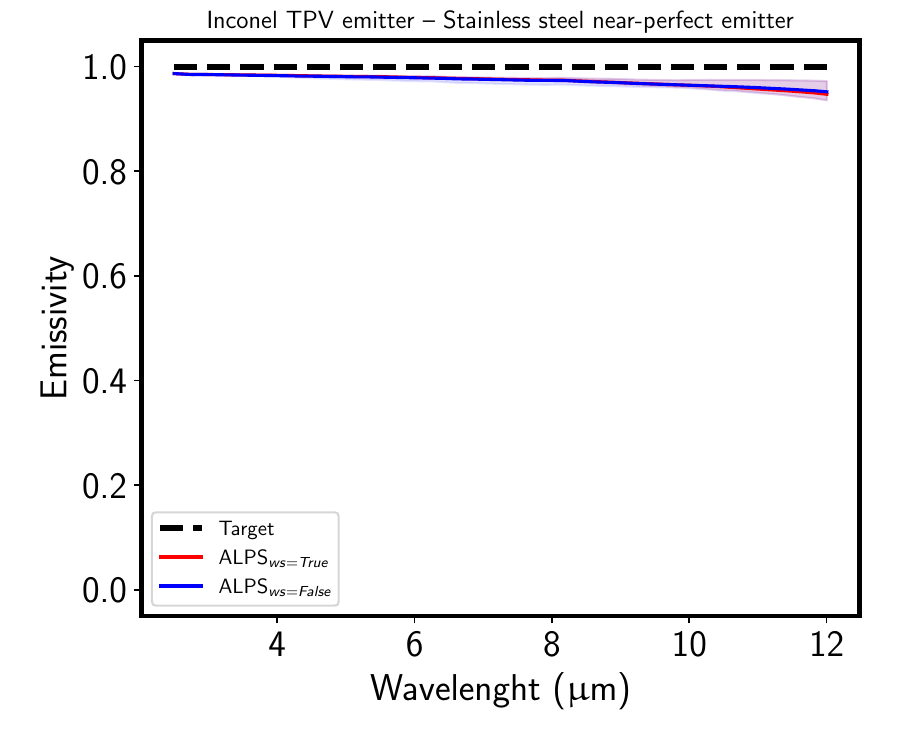}     
   \caption{}
   \label{fig:Figure5g} 
\end{subfigure}
\begin{subfigure}[b]{0.49\textwidth}
   \includegraphics[width=\linewidth]{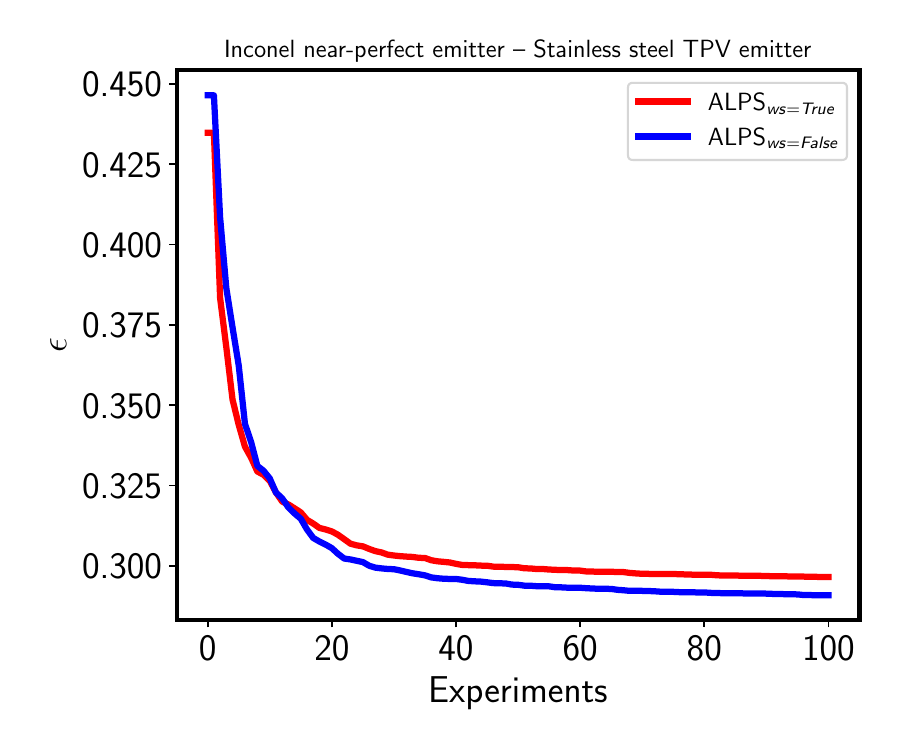}
   \caption{}
   \label{fig:Figure5d}
\end{subfigure}
\begin{subfigure}[b]{0.49\textwidth}
   \includegraphics[width=\linewidth]{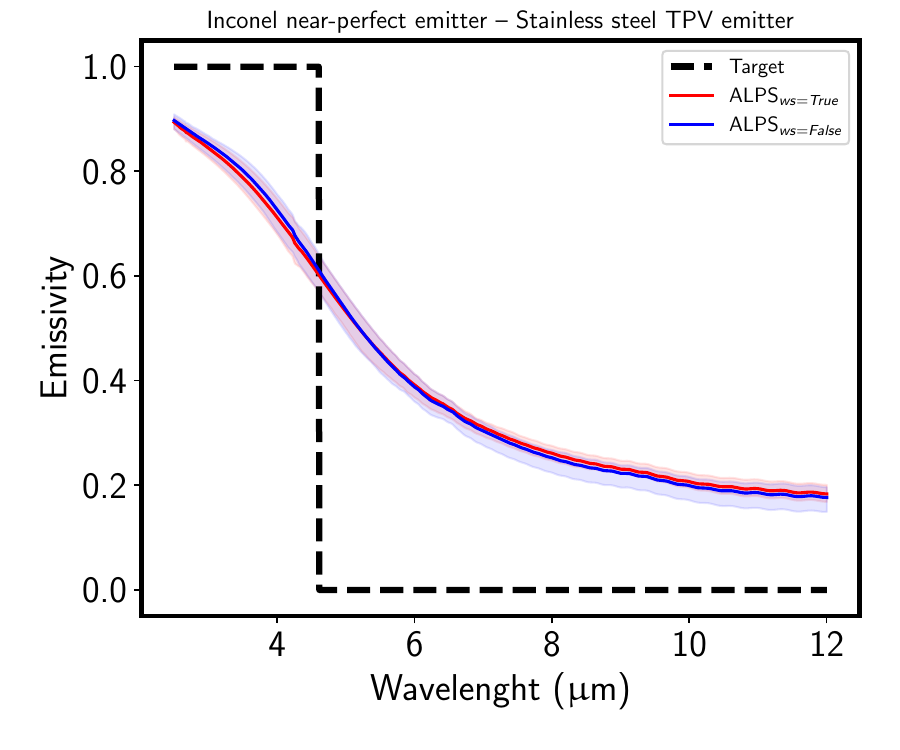}
   \caption{}
   \label{fig:Figure5h}
\end{subfigure}
\caption[]{The inverse design results of ALPS with and without cross-material and cross-target warm starting for the photonic surfaces benchmarks (ALPS$_{ws=True}$ means ALPS with warm starting and is shown as the red line). The title of each graph presents the warm starting model source (stated first) and inverse design target (stated last). Convergence graphs are shown in the first column, while in the second, the solutions reconstruction graphs are shown: (a) Convergence graphs for the Stainless steel near-perfect emitter target benchmark. (b) Solution reconstruction graph for the Stainless steel near-perfect emitter. (c) Convergence graphs for the Stainless steel TPV emitter target benchmark. (d) Solution reconstruction graph for the Stainless steel TPV emitter.}
\label{fig:Figure5_2}
\end{figure}

Fig. \ref{fig:Figure5_1} and Fig. \ref{fig:Figure5_2} show the convergence and solution reconstruction graphs for ALPS with and without cross-material warm starting. In this case, additionally to switching to a different inverse design target, we also include switching to a different material, i.e., when the Stainless steel near-perfect emitter is the target, ALPS is warm started by a model that is generated after an inverse design of the Inconel TPV emitter, and vice-versa. The convergence graphs in the first column of Fig. \ref{fig:Figure5_1} and \ref{fig:Figure5_2} show that cross-material warm starting can be very beneficial (Fig. \ref{fig:Figure5a} and Fig. \ref{fig:Figure5c}), or it can either bring no significant boost in performance (Fig. \ref{fig:Figure5b}), or even slightly worsen the performance (Fig. \ref{fig:Figure5d}). The worsened performance of ALPS with warm starting for the Stainless steel TPV emitter target can be explained with the difference between lower and upper bound difference of the spacing parameter of the experimental models between Inconel and Stainless steel. The best found laser parameters for the Stainless steel TPV target by the ALPS model without warm starting are 0.9 W for laser power, 338 mm/s for scanning speed, and 3.37 $\upmu$m for the spacing. The spacing domain for the Inconel experimental model ranges from 15 $\upmu$m to 28 $\upmu$m, as shown in Fig. \ref{fig:FigureA2}. Therefore, the initial solutions generated by the warm starting model are limited to laser parameters within this range. In this context, employing LHS is advantageous as it provides an unbiased approach, enabling the uniform generation of initial samples across the entire spacing domain.

\newpage

\newpage

\newpage

\newpage

\newpage
\section{Conclusion}\label{sec:conclusion}

We introduce and investigate the ALPS algorithm, developed for the inverse design of photonic surfaces. ALPS utilizes a greedy surrogate sampling method, employing the RF algorithm as a surrogate, and has demonstrated superior performance over established optimization algorithms in all benchmark tests. Notably, ALPS offers reliable convergence with limited experimental model evaluations, a significant advantage in photonic surface inverse design practical applications. Additionally, its capability for cross-target and cross-material warm starting enhances its utility in repetitive design processes across similar or identical design spaces. Future research could explore the application of ALPS in real-world autonomous experimental settings for further advancing photonic surface design, as well as generalizing the proposed approach to other problem settings.

\section*{Acknowledgments}
This work was supported by the Laboratory Directed Research and Development Program of Lawrence Berkeley National Laboratory under U.S. Department of Energy Contract No. DE-AC02-05CH11231.  M\"uller's time was supported under U.S. Department of Energy Contract  No. DE-AC36-08GO28308, U.S. Department of Energy Office of Science, Office of Advanced Scientific Computing Research, Scientific Discovery
through Advanced Computing (SciDAC) program through the FASTMath Institute to the National Renewable Energy Laboratory.

\section*{Author Contributions}
L.G. wrote the manuscript, developed the methods, developed the code, designed the numerical experiments, and analyzed the performance of the algorithms, M.P. provided the experimental data and edited the manuscript, J.M. supervised the research and edited the manuscript, V.Z. and W.A.J. supervised the research, provided funding and edited the manuscript.

\section*{Declaration of Competing Interest}
The authors declare that they have no known competing financial interests or personal relationships that could have appeared to influence the work reported in this paper.

\section*{Data Availability}
The machine learning models, data and benchmarks needed to reproduce the study can be found on the following repository: \url{https://github.com/lukagrbcic/ALPS-Data}
\newpage
\begin{appendices}

\section{Experimental Datasets and Model}\label{app:experimental_model}

In this section we describe the experimental dataset used to train the experimental models used to assess the performance of our methods. We also include the details of the RF-PCA algorithm, the model validation procedures and assess the accuracy of the models.

\subsection{RF-PCA Algorithm}\label{subapp:dataset}

The experimental model is based on a combination of the RF algorithm and the PCA algorithm used to compress the large dimensionality of the spectral emissivity space. The RF algorithm (introduced by \citet{breiman2001random}) is an ensemble learning method for ML classification and regression. It constructs randomly defined decision trees during training, outputting the mean prediction of individual trees. This randomness, introduced by selecting a subset of features for each tree split, helps de-correlate trees, reducing overfitting and allowing for interpretability. It is particularly suited for modeling problems that have well structured features like laser fabrication parameters. 

\begin{figure}[!h]
  \includegraphics[width=\linewidth]{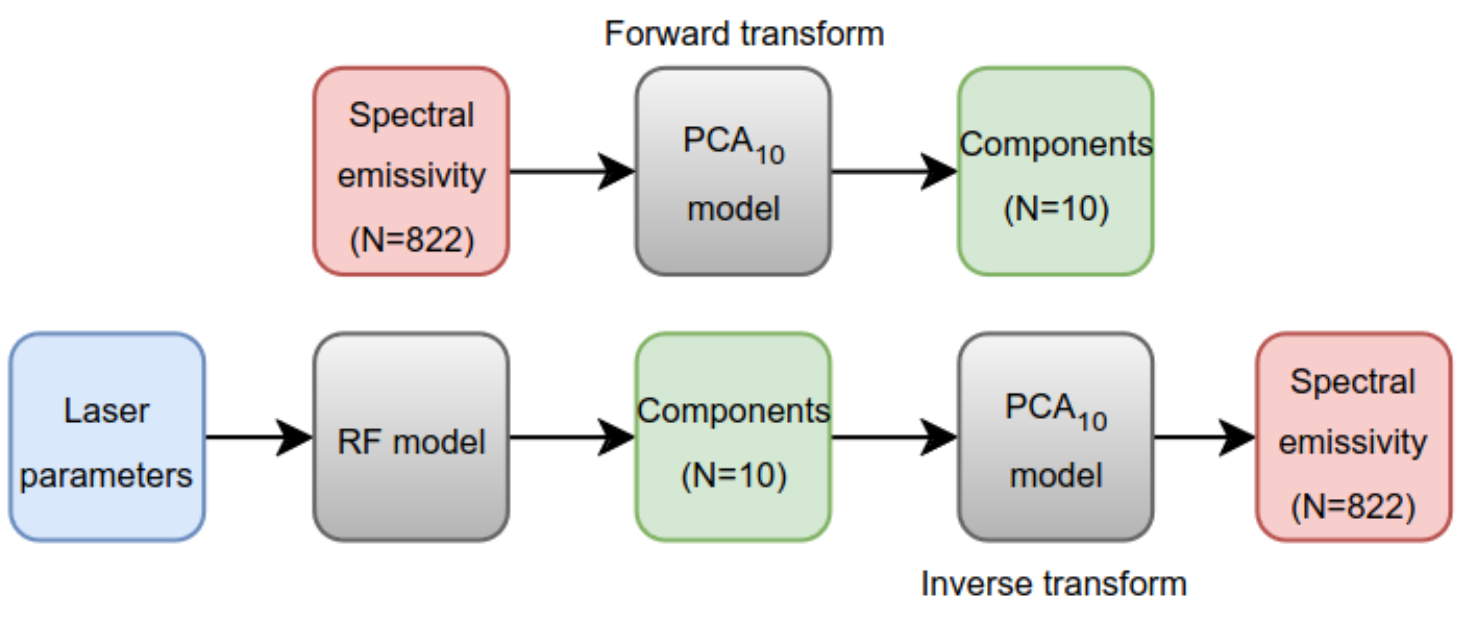}
  \caption{The RF-PCA algorithm is employed to train the experimental models. The top section of the figure illustrates the PCA model's pipeline (PCA$_10$ implies that the model transforms inputs into 10 components), which processes spectral emissivity values and converts them into a ten-dimensional principal component space. The bottom section depicts the complete RF-PCA model pipeline, where laser parameters serve as inputs to the RF model to predict PCA components. These components are then inversely transformed back into the spectral emissivity space.}
  \label{fig:FigureA1}
\end{figure}

The RF-PCA algorithm comprises two main parts, each with its own pipeline, as shown in Fig. \ref{fig:FigureA1}. Initially, the PCA model is trained to convert spectral emissivity curves into ten-dimensional PCA components. Subsequently, the RF model is trained using laser parameters to predict these PCA components. The PCA model is then employed to inversely transform these predictions back into the spectral emissivity space, thus modeling the relationship between the laser parameters and spectral emissivity. This approach builds on the results discussed in our earlier publication (\citet{grbcic2024ensemble}). It is important to clarify that incorporating the PCA model is not essential for accurately modeling the relationship between laser parameters and spectral emissivity. Rather, the primary function of PCA compression is to enhance the computational efficiency of the entire experimental model and to significantly decrease the storage size of the experimental model used for generating function evaluations in our ALPS numerical experiments.

The Python module scikit-learn 1.2.2. implementation of both the RF and PCA algorithms were used to train the experimental models (\citet{pedregosa2011scikit}). Most hyperparameters of the RF algorithm are set to default values, however, the number of estimators (trees) was set to 450, and the max depth value was set to 10, as they are determined to be optimal for this task through numerical experimentation.

\subsection{Experimental Datasets}\label{subapp:dataset}

\begin{figure}[!h]
\centering
\begin{subfigure}[b]{0.49\textwidth}
   \includegraphics[width=\linewidth]{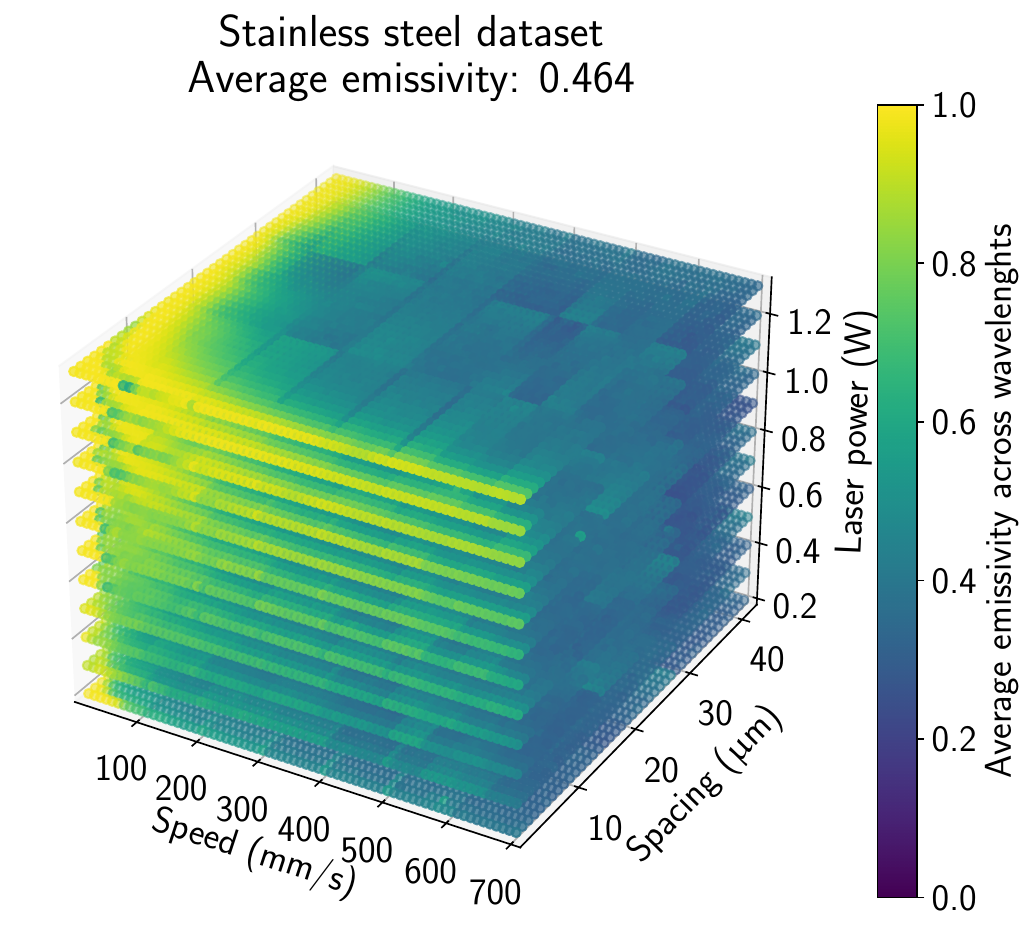}
   \caption{}
   \label{fig:FigureA2a}
\end{subfigure}
\begin{subfigure}[b]{0.49\textwidth}
   \includegraphics[width=\linewidth]{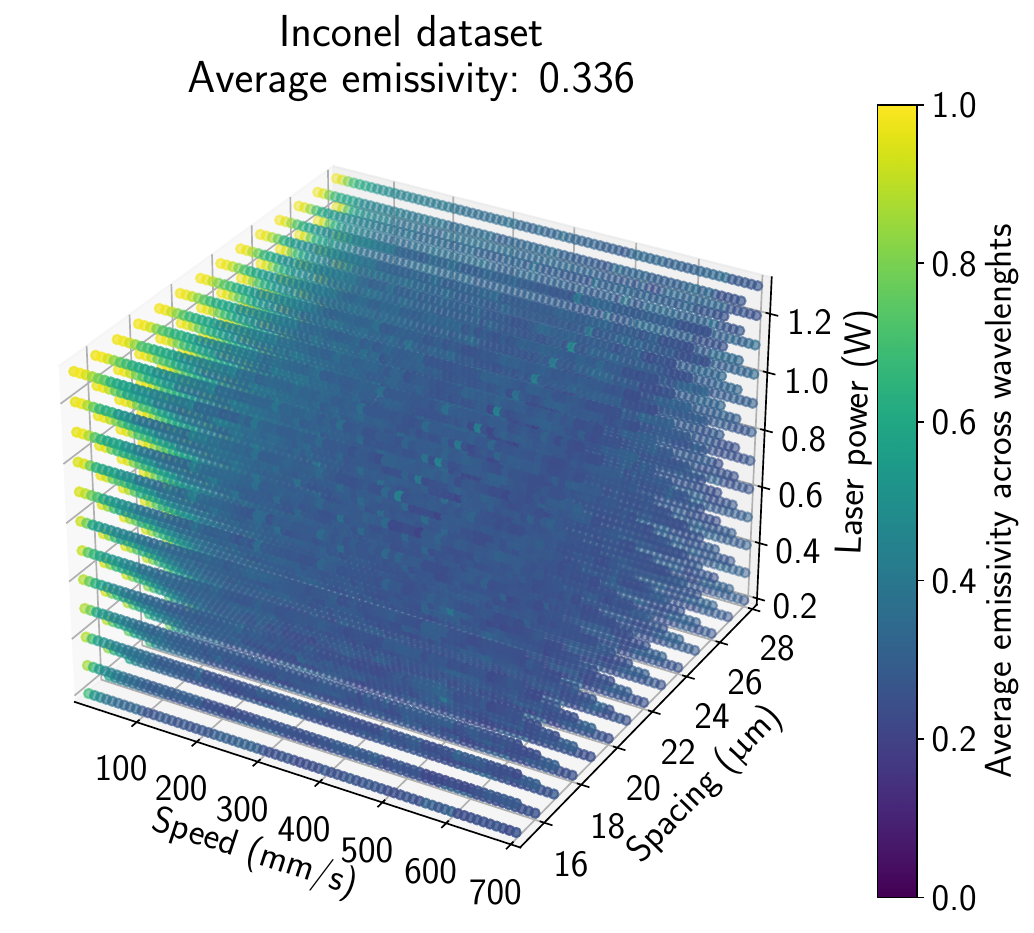}     
   \caption{}
   \label{fig:FigureA2b} 
\end{subfigure}
\begin{subfigure}[b]{0.49\textwidth}
   \includegraphics[width=\linewidth]{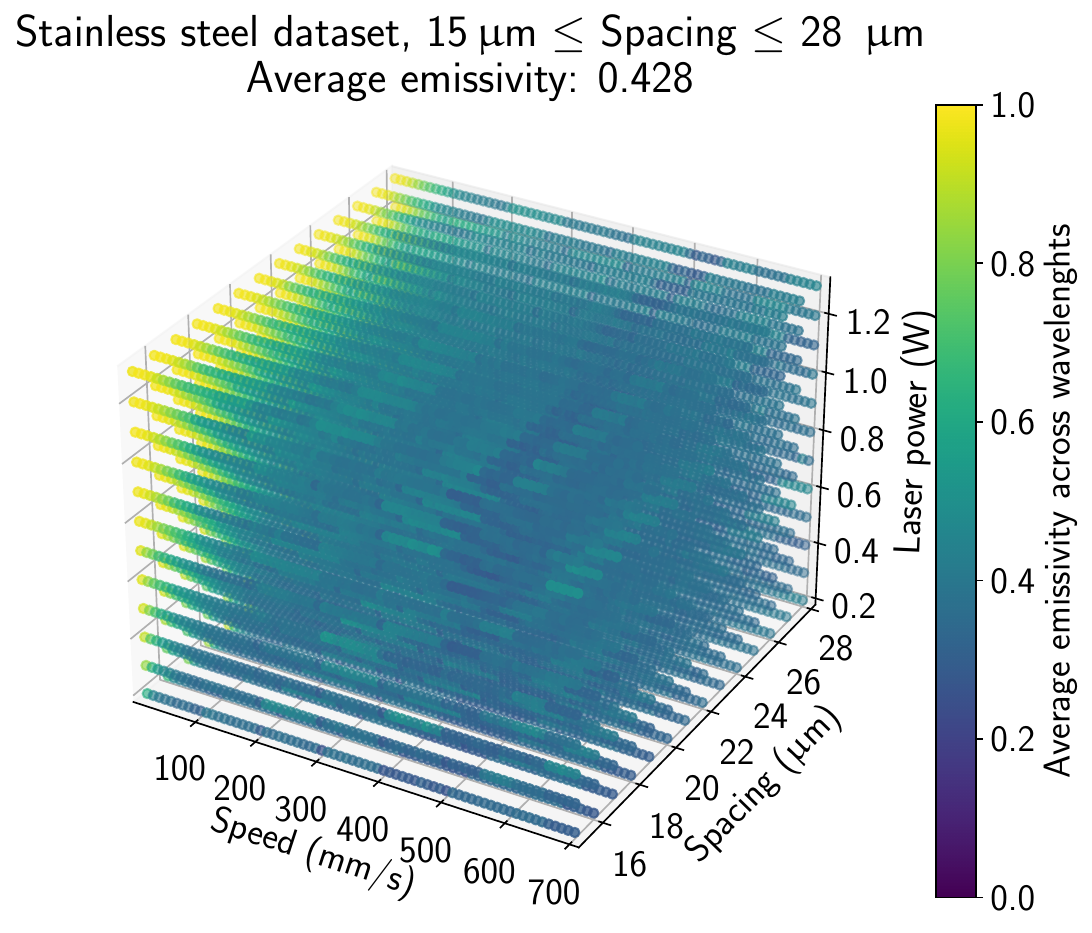}     
   \caption{}
   \label{fig:FigureA2c}
\end{subfigure}
\caption{Experimental dataset distributions in laser fabrication parameter space, colored by the average emissivity value across all wavelengths for all data: (a) Stainless steel dataset where the total number of data is 35,326. (b) Inconel dataset where the total number of data is 11,759. (c) Subset of the Stainless steel dataset using only the same range for the spacing parameter as the Inconel dataset to showcase similarity.}
\label{fig:FigureA2}
\end{figure}

Fig. \ref{fig:FigureA2} illustrates the experimental dataset distributions for both materials, namely Inconel and Stainless steel. The laser fabrication parameter space distribution of dataset instances for the Stainless steel material can be observed in Fig. \ref{fig:FigureA2a}, while Fig. \ref{fig:FigureA2b} displays the Inconel dataset distribution. All data instances are color-coded based on their corresponding average emissivity values across the wavelength domain (2.5 $\upmu$m to 12 $\upmu$m). The total number of data points in the Stainless steel dataset is 35,326, whereas in the Inconel dataset, it is 11,759. The  range of the spacing parameter in the Inconel dataset is narrower than that of the Stainless steel dataset, and the Stainless steel dataset exhibits a higher overall average emissivity value (computed as the average of all instances' individual average emissivities) compared to the Inconel dataset. The experimental sampling for both datasets is uniform and thoroughly covers the entire laser parameter space.

Fig. \ref{fig:FigureA2c} presents a subset of the Stainless steel dataset with the range of the spacing parameter constrained to match that of the Inconel dataset. This comparison aims to highlight the resemblance between the two datasets. It is observed that there are similarities in the gradient of the average emissivity across the laser fabrication space. However, it is noteworthy that the average emissivity of the Stainless steel dataset remains higher (0.428) compared to that of the Inconel dataset (0.336). The similarity of the gradient potentially implies that using a pretrained model for Stainless steel could be used for warm starting the photonic surface inverse design process that is based on the Inconel material.

For full details of the experimental procedure that is used to generate these datasets, readers are encouraged to refer to our prior publications (\citet{park2024tnn} and \citet{grbcic2024ensemble}).

\newpage
\subsection{Experimental Models Validation}\label{subapp:ml_exp_results}

To train and validate the experimental models, we randomly selected 20,000 instances from the entire Stainless steel dataset and 11,000 instances from the Inconel dataset. Utilizing the whole dataset to train the models is not considered since an increase in data yields an increase in the storage space of the models. We divided each dataset into 75\% for training and 25\% for testing. Next, we conducted a learning curve analysis on the 75\% of the training dataset using a K-fold ($K=3$) cross-validation method to assess the accuracy and uncertainty of the predictions of both models. The primary metric used to assess the models is the RMSE (defined in Eq. (\ref{eqn:rmse})).

\begin{figure}[!h]
\centering
\begin{subfigure}[b]{0.49\textwidth}
   \includegraphics[width=\linewidth]{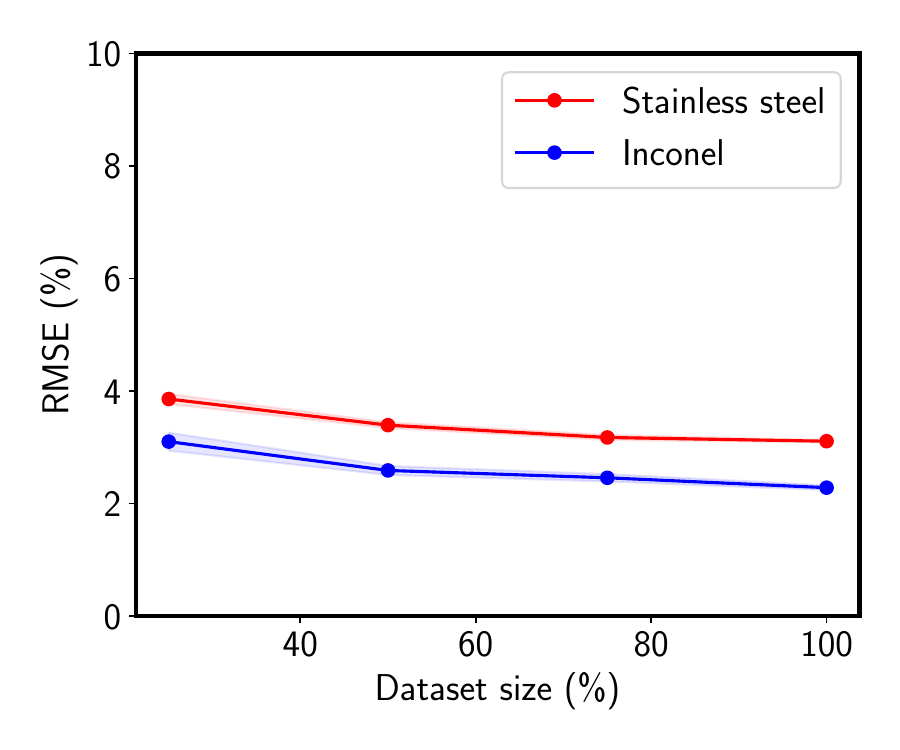}
   \caption{}
   \label{fig:FigureA3a}
\end{subfigure}
\begin{subfigure}[b]{0.49\textwidth}
   \includegraphics[width=\linewidth]{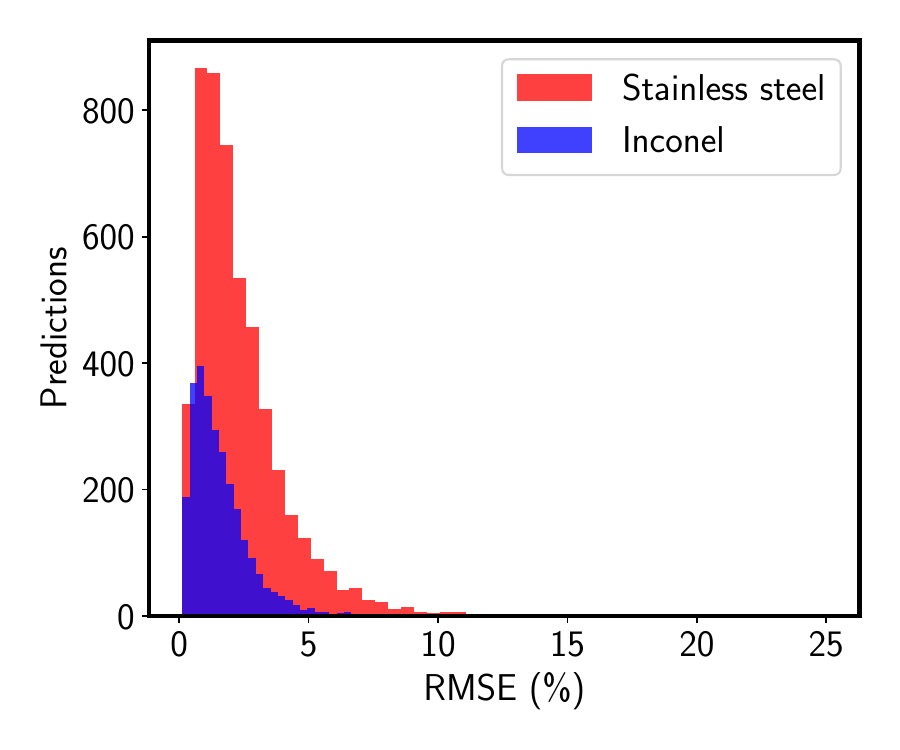}     
   \caption{}
   \label{fig:FigureA3b} 
\end{subfigure}
\caption{Results of the RF-PCA experimental model validation: (a) Learning curve of both Stainless steel and Inconel models. (b) Histogram visualization of RMSE values for each predicted instance when compared to the testing set instances. }
\label{fig:FigureA3}
\end{figure}

Fig. \ref{fig:FigureA3} presents the outcomes of the analysis for both models. Fig. \ref{fig:FigureA3a} displays the learning curves for both models, illustrating the impact of varying dataset sizes—25\%, 50\%, 75\%, and 100\% of the training set (75\% of the whole dataset)—on model accuracy through a cross-validation analysis. The results indicate that increasing the dataset size does not significantly enhance the accuracy of either model. The Stainless Steel model exhibits a marginally higher RMSE, which may be attributed to a broader range in laser parameter space (i.e., a larger range of the  spacing parameter). Fig. \ref{fig:FigureA3b} depicts the results obtained when the entirety of the training data (75\% of the whole dataset) is used to train the models, and the testing data (25\% of the whole dataset) is utilized for model assessment. This is illustrated through a histogram of RMSE values derived from comparing the model predictions against the test set instances. For both models, the majority of RMSE values are concentrated below 5\%. The Stainless Steel model has average, maximum, and standard deviation RMSE values of 2.3\%, 24.0\%, and 1.8\%, respectively. For the Inconel model, these values are 1.6\%, 14.4\%, and 1.4\%. Although the maximum RMSE values indicate poor performance on some outliers, the models generally demonstrate excellent performance based on their average RMSE values. The models that are used for inference during the ALPS process are trained only with the training set (75\% of each selected experimental dataset) to further reduce the storage size of the models.

Finally, we evaluated the accuracy of the PCA model for reconstructing spectral emissivity. The spectral emissivity curves from the training set are utilized to train the PCA model with 10 components. The trained PCA model then transformed the test set spectral emissivity curves into a ten-dimensional principal component space. These components are subsequently inversely transformed back to the spectral emissivity space and compared with the original curves from the test set using RMSE. The average RMSE for the Stainless steel PCA model (calculated by averaging the RMSE between each test set instance and its PCA reconstruction) is 0.08\%, and for the Inconel model, it is 0.074\%. These extremely low RMSE values suggest that the PCA model effectively compresses the spectral emissivity curves.

\newpage

\section{Synthetic Benchmarks and Numerical Experiments Setup}\label{app:exp_setup}

In this section, we introduce the supplementary inverse design synthetic benchmarks employed to evaluate the performance of the ALPS algorithm and its comparison against other optimization techniques. Additionally, we detail the algorithms selected for comparison, including their implementation specifics and hyperparameters.

\subsection{Synthetic Benchmarks}\label{subapp:synbench}

\begin{figure}[!h]
\centering
\begin{subfigure}[b]{0.49\textwidth}
   \includegraphics[width=\linewidth]{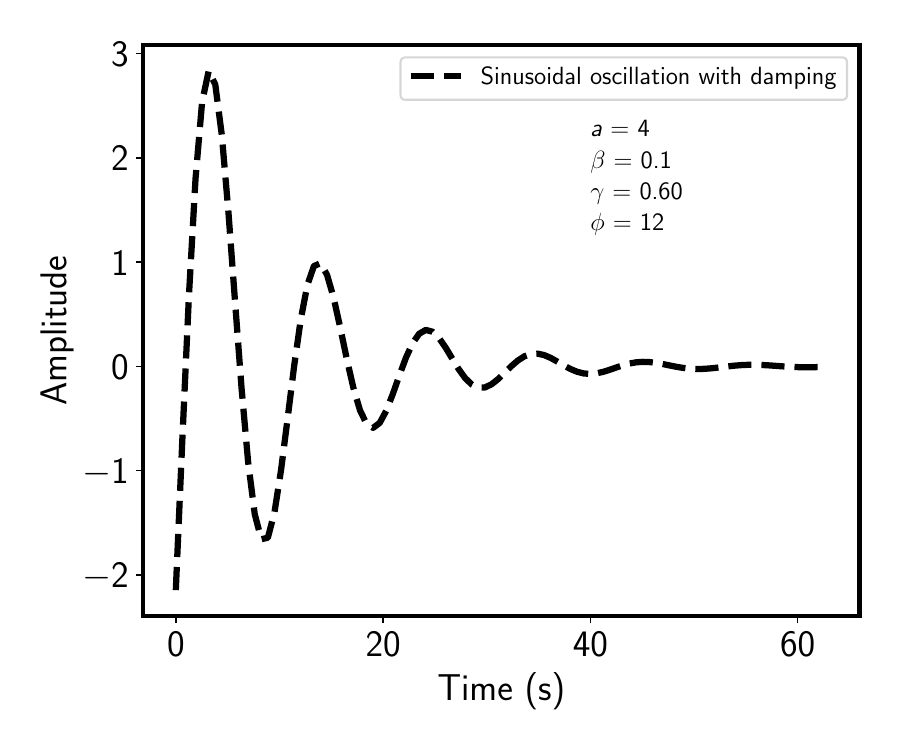}
   \caption{}
   \label{fig:FigureB1a}
\end{subfigure}
\begin{subfigure}[b]{0.49\textwidth}
   \includegraphics[width=\linewidth]{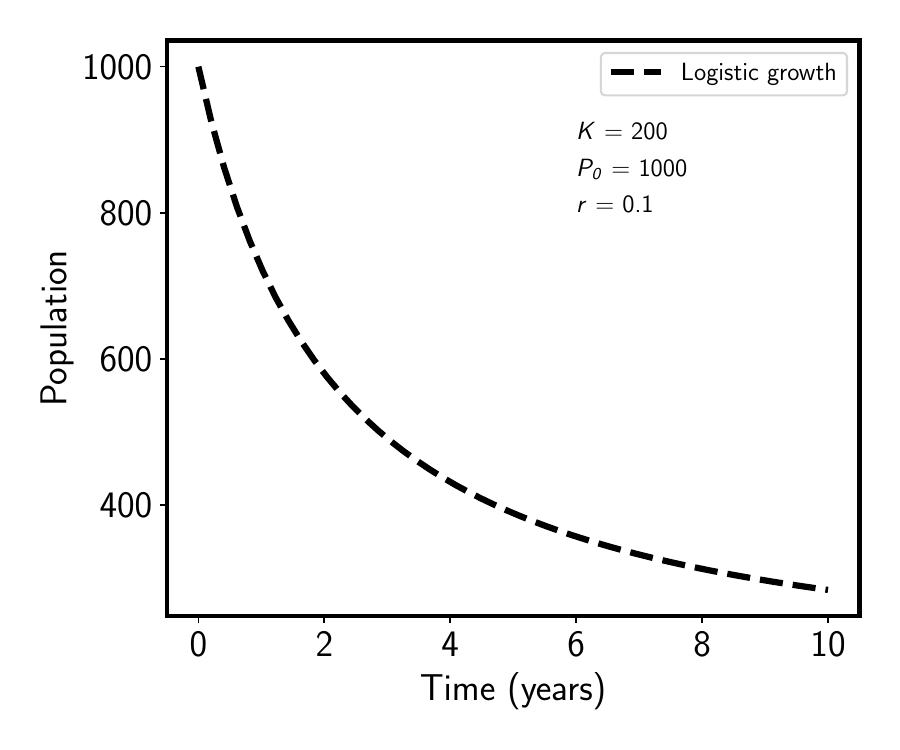}     
   \caption{}
   \label{fig:FigureB1b} 
\end{subfigure}
\caption{The synthetic benchmark targets with annotated model coefficients needed to reconstruct them: (a) The sinusoidal oscillation with damping model. (b) The logistic growth model.}
\label{fig:FigureB1}
\end{figure}

The reason for inclusion of inverse design synthetic benchmarks is to showcase that the ALPS approach generalizes beyond the photonic surface inverse design problem.  Fig. \ref{fig:FigureB1} shows the two additional targets, namely the sinusoidal oscillation with damping model and  the logistic growth model used for population modeling. These synthetic benchmarks are chosen since the target value is a curve that is parameterized with a low dimensional design vector, similarly to the photonic surface inverse design benchmarks. Both synthetic benchmarks are defined through the inverse design mathematical framework presented in Eq. (\ref{eqn:inverse_design_definition}), (\ref{eqn:inverse_optimization}) and (\ref{eqn:rmse}). 

The mathematical expression for the sinusoidal oscillator with damping is defined in Eq. (\ref{eqn:oscillator}) as:

\begin{equation}
\begin{aligned}
\mathbf{y} = A(t) = a \cdot e^{-\beta \cdot t} \cdot \sin{(\gamma \cdot t + \phi)}
\label{eqn:oscillator}
\end{aligned}
\end{equation}

where $\mathit{a}, \beta, \gamma,$ and $\phi$ are the model coefficients that need to be determined in order to reconstruct the target vector $\mathbf{y}$ (displayed in Fig. \ref{fig:FigureB1a} with the coefficient values). For this case, the design vector is $\mathbf{x} = [a, \beta, \gamma, \phi]^T$. The amplitude $A$ is discretized with 100 time steps $t$ in the domain from 0 to 20$\pi$.

The mathematical expression for the logistic growth model is defined in Eq. (\ref{eqn:logistic}) as:

\begin{equation}
\begin{aligned}
\mathbf{y} = P(t) = \frac{K}{1 + (\frac{K - P_0}{P_0}) \cdot e^{-r \cdot t}}
\label{eqn:logistic}
\end{aligned}
\end{equation}

where $\mathit{K}, P_0,$ and $\mathit{r}$ are model coefficients (values annotated in Fig. \ref{fig:FigureB1b}) that define the design vector $\mathbf{x} = [K, P_0, r]^T$. The population $P$  target value vector is discretized using 50 points of time $t$ with a total range from 0 to 10 years, i.e. the target vector $\mathbf{y}$ has 50 components. The lower and upper boundaries for inverse design optimization (Eq. (\ref{eqn:inverse_optimization})) of the sinusoidal oscillation with damping and the logistic growth model, $\mathbf{x_{lb}}$ and $\mathbf{x_{ub}}$, are defined in Tab. \ref{tab:syn_bounds}.
     
\begin{table}[!h]
\centering
\caption{The lower and upper bounds of the inverse design optimization problem, $\mathbf{x_{lb}}$ and $\mathbf{x_{ub}}$, respectively, for both synthetic benchmarks (defined in the first column).}
\begin{tabular}{ c | c | c}
 \textbf{Synthetic benchmark}&  $\mathbf{x_{lb}} $ &  $\mathbf{x_{ub}}$ \\
\hline
\hline
Sinusoidal oscillation & $\mathit{a}$ = 2,  $\beta$ = 0.05, $\gamma$ = 0, $\phi$ = 3 & $\mathit{a}$ = 5,  $\beta$ = 0.4, $\gamma$ = 2, $\phi$ = 15  \\
Logistic growth & $\mathit{K}$ = 100,  $P_0$ = 100, $\mathit{r}$ = 0.01 & $\mathit{K}$ = 1200,  $P_0$ = 1400, $\mathit{r}$ = 0.4 \\
\end{tabular}

\label{tab:syn_bounds}
\end{table}

\newpage

\subsection{Algorithms used for Comparison}\label{subapp:algorithms}

To highlight the performance  of  ALPS, we conduct a comparative analysis with other established algorithms from various optimization algorithm categories. From the population-based algorithms category, we select PSO and DE (\citet{kennedy1995particle}, \citet{storn1997differential}). Mesh Adaptive MADS and NM are chosen from the direct search algorithms category (\citet{audet2006mesh}, \citet{lagarias1998convergence}). L-BFGS-B is included as a representative of the gradient-based algorithms, while BO is a   the model-based optimizer (\citet{shahriari2015taking}).

\begin{table}[!h]
\centering
\caption{The optimization algorithms and their corresponding categories, that are used for comparison with ALPS on the photonic surface inverse design problem, as well as the two synthetic benchmark problems. }
\begin{tabular}{ c | c | c}
 \textbf{Algorithm} &  \textbf{Category} & \textbf{Python module} \\
\hline
\hline
PSO & Population-based & Indago v0.5.0 (\citet{Indago}) \\
DE & Population-based & Indago v0.5.0 (\citet{Indago})\\
MADS & Direct search & OMADS v2.1.0 (\citet{OMADS_AB}) \\
NM & Direct search & scipy v1.11.4 (\citet{virtanen2020scipy}) \\
LBFGSB & Gradient-based  & scipy v1.11.4 (\citet{virtanen2020scipy}) \\
BO & Model-based  & Our implementation \\
\end{tabular}

\label{tab:optimization_algorithms}
\end{table}

The hyperparameters of all used algorithms (except BO and MADS) are set to their default recommended values by each Python module used, and can be found in  each module's documentation. All MADS hyperparameters are set as default (\citet{OMADS_AB}) except for the scaling factor which is set to 10 as it is found to be a good value through numerical experimentation for all benchmarks. MADS, NM and LBFGSB are initialized using uniform random sampling.

Finally, the BO Python implementation is custom and it is based on the GP model with the Expected Improvement (EI) acquisition function (\citet{zhan2020expected}). The exploration hyperparameter of the EI is set to 0.01. The scikit-learn 1.2.2. (\citet{pedregosa2011scikit}) implementation of the GP model is utilized with the Matern kernel and the length scale parameter set to 1. The smoothness parameter $\nu$ is set to 2.5. This configuration was selected as it is one of the most popular default configurations used in Python BO modules such as scikit-optimize (\citet{head_2020_4014775}). Furthermore, in the GP model, the number of optimizer restarts is set to 2, and the target normalization feature is set to True. The L-BFGS-B optimizer with the Python module default parameters is used to find the maximum of EI in order to determine the best sample during each iteration (\citet{virtanen2020scipy}). The initial samples (initial sample size is 5 samples) of each BO run is generated using LHS to be as similar as possible to the ALPS framework, and the termination criterion is also set as the maximum number of allowed experimental model evaluations.

For a fair comparison with ALPS, we employ the custom BO implementation, avoiding the advanced model hyperparameter techniques often used in standard BO modules, usually with increased computational cost. These techniques could potentially enhance ALPS as well; however, our study aims to highlight the simplicity and effectiveness of our approach for inverse design. This approach performs well even when using out-of-the-box ML models with default hyperparameters.

\newpage
\section{Inverse Design Convergence Results and Hyperparameter Analysis}\label{app:fullresults}

In this section we present the complete results for all benchmarks and algorithms used for our numerical experiments, presented through both convergence graphs with uncertainty, and through solution reconstruction graphs. We also include a detailed analysis of the hyperparameters that are used in ALPS and how they affect its performance.

\subsection{Synthetic Benchmarks Convergence}\label{subapp:synbencresults}

\begin{figure}[!h]
\centering
\begin{subfigure}[b]{0.49\textwidth}
   \includegraphics[width=\linewidth]{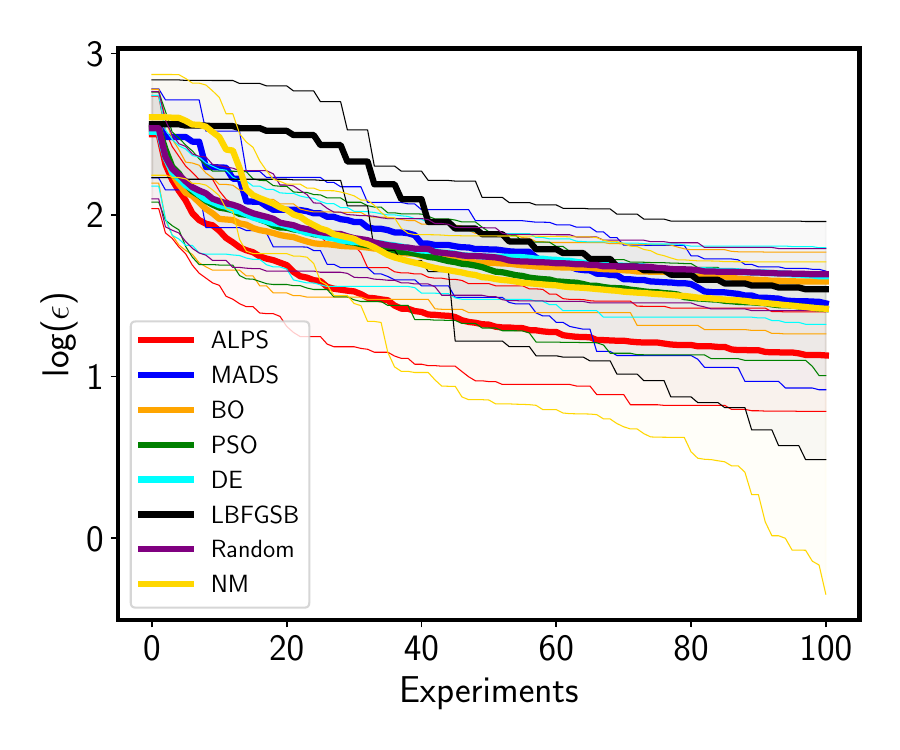}
   \caption{}
   \label{fig:FigureC1a}
\end{subfigure}
\begin{subfigure}[b]{0.49\textwidth}
   \includegraphics[width=\linewidth]{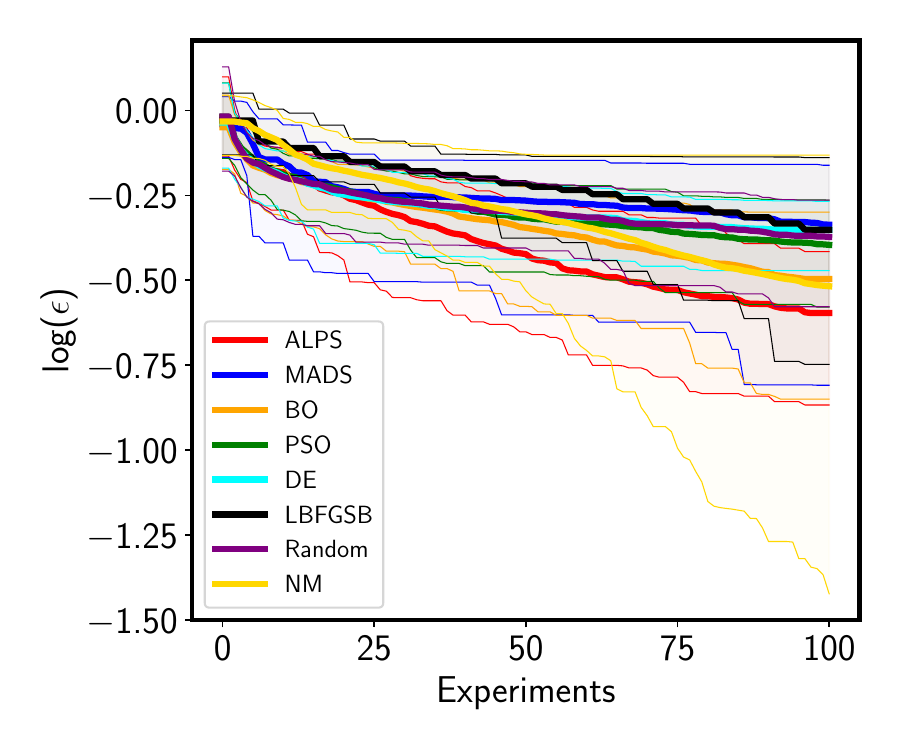}     
   \caption{}
   \label{fig:FigureC1b} 
\end{subfigure}
\caption{Convergence graphs display the base-10 logarithm of the mean error value (Eq. (\ref{eqn:inverse_design_definition})) for the synthetic benchmarks. The plots feature mean error values of 100 runs, represented by thick lines. Uncertainty bounds (shaded regions) are indicated by the 10$^{th}$ and 90$^{th}$ percentiles of the error across all algorithms: (a) Logistic growth target. (b) Sinusoidal oscillation with damping target.}
\label{fig:FigureC1}
\end{figure}

\begin{figure}[!ht]
\centering
\adjustbox{center}{
\begin{minipage}{1.5\textwidth}
\begin{subfigure}[b]{0.24\textwidth}
   \includegraphics[width=\linewidth]{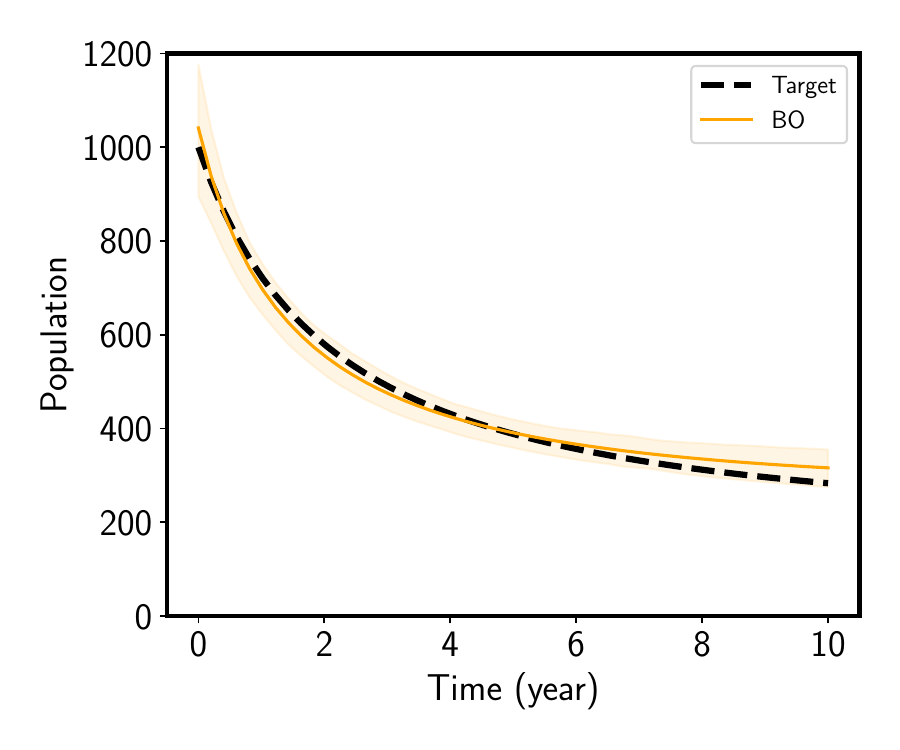}
   \caption{}
   \label{fig:FigureC2a}
\end{subfigure}
\hfill 
\begin{subfigure}[b]{0.24\textwidth}
   \includegraphics[width=\linewidth]{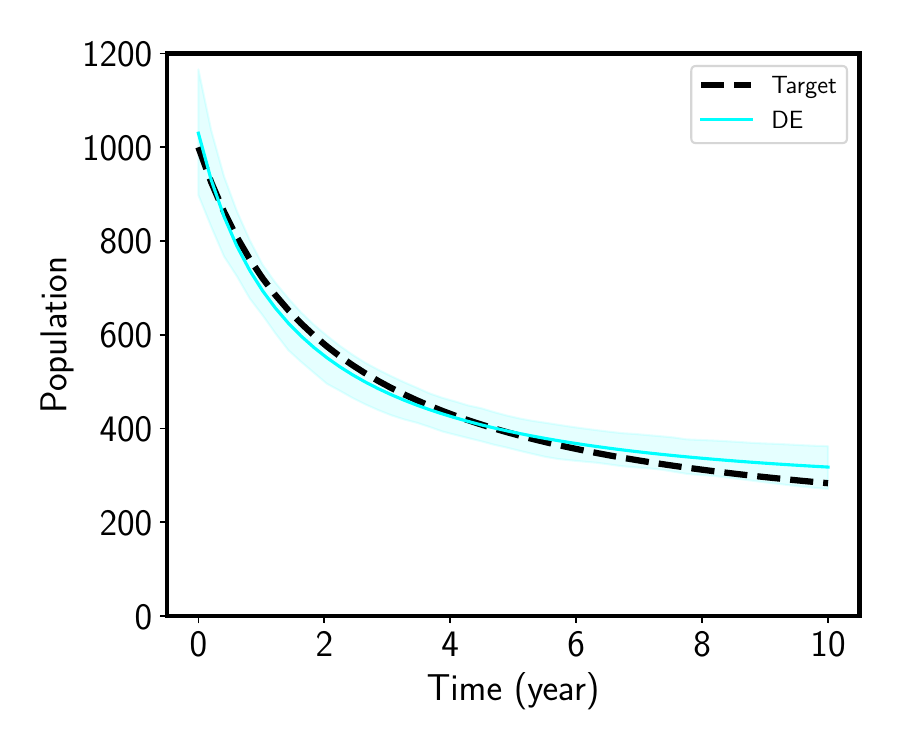}     
   \caption{}
   \label{fig:FigureC2b} 
\end{subfigure}
\hfill
\begin{subfigure}[b]{0.24\textwidth}
   \includegraphics[width=\linewidth]{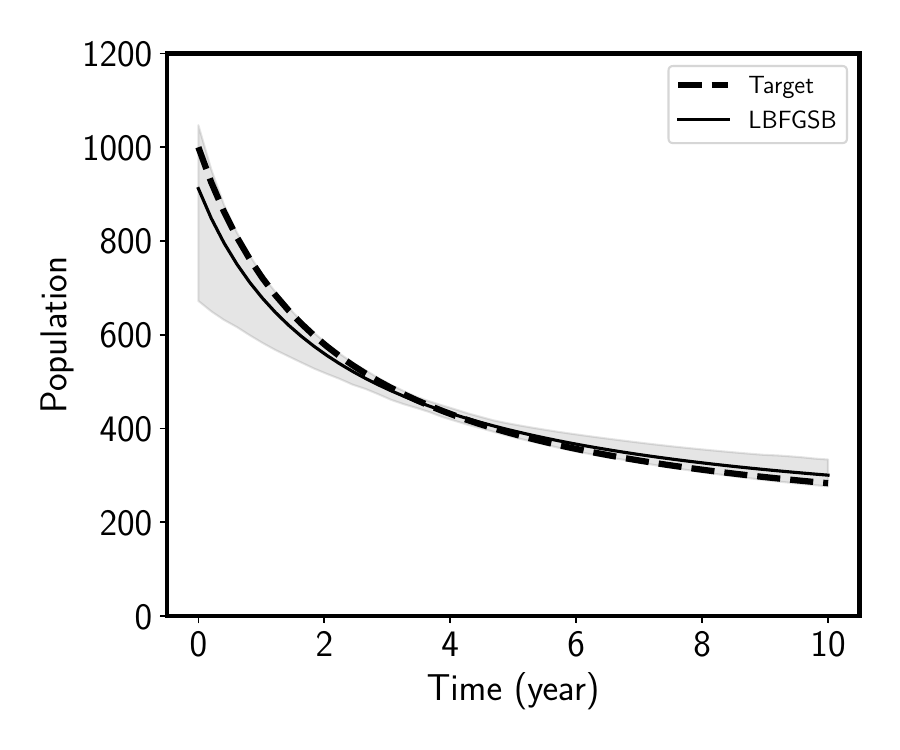}     
   \caption{}
   \label{fig:FigureC2c} 
\end{subfigure}
\hfill
\begin{subfigure}[b]{0.24\textwidth}
   \includegraphics[width=\linewidth]{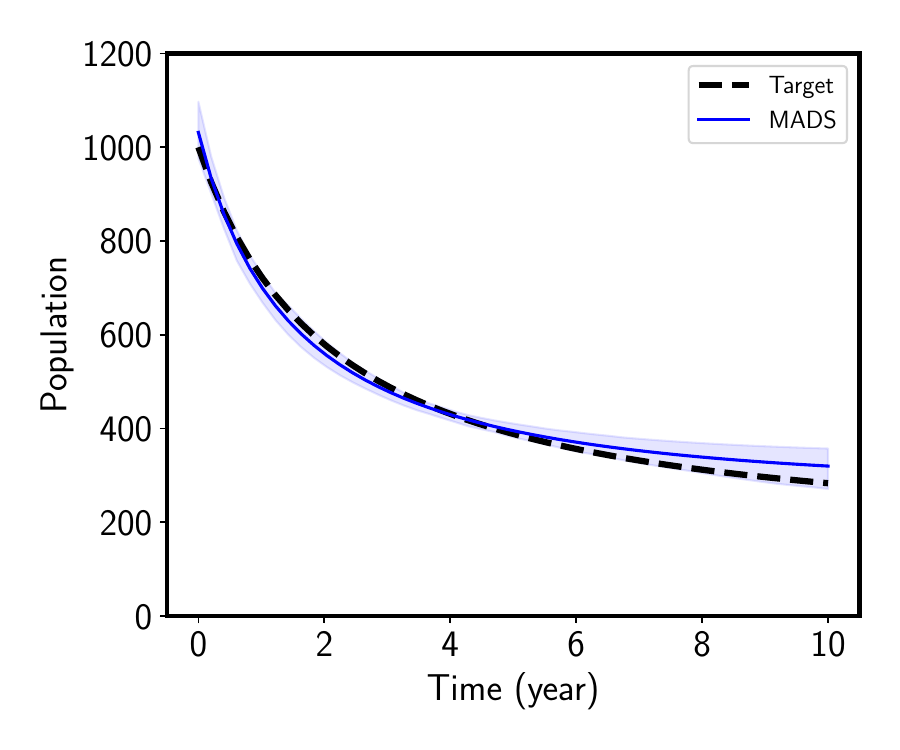}     
   \caption{}
   \label{fig:FigureC2d} 
\end{subfigure}
\end{minipage}
}

\adjustbox{center}{
\begin{minipage}{1.5\textwidth}
\begin{subfigure}[b]{0.24\textwidth}
   \includegraphics[width=\linewidth]{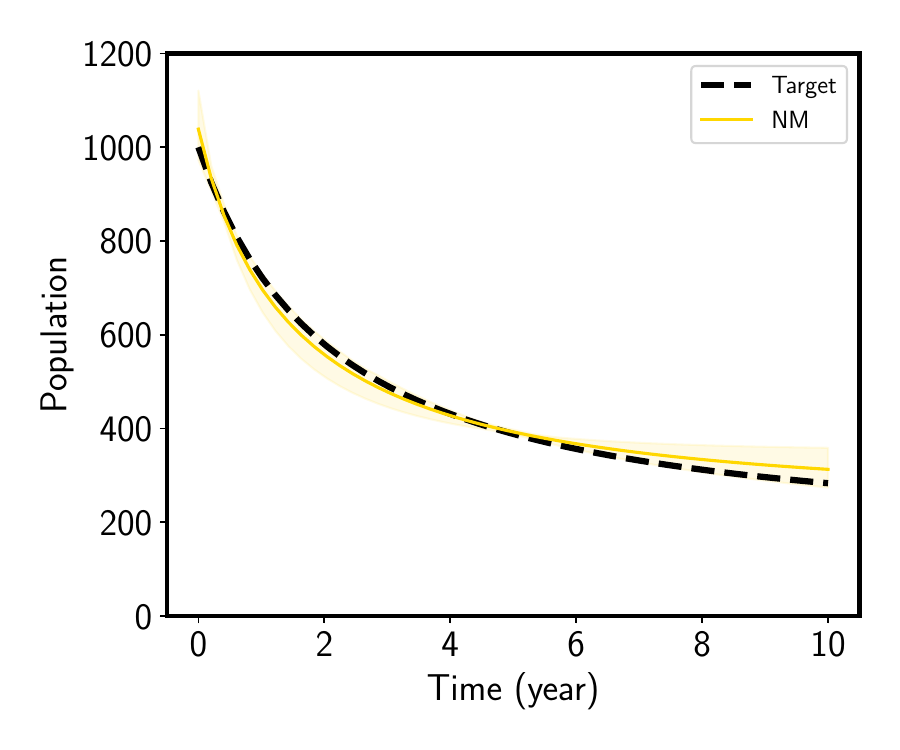}     
   \caption{}
   \label{fig:FigureC2e} 
\end{subfigure}
\hfill
\begin{subfigure}[b]{0.24\textwidth}
   \includegraphics[width=\linewidth]{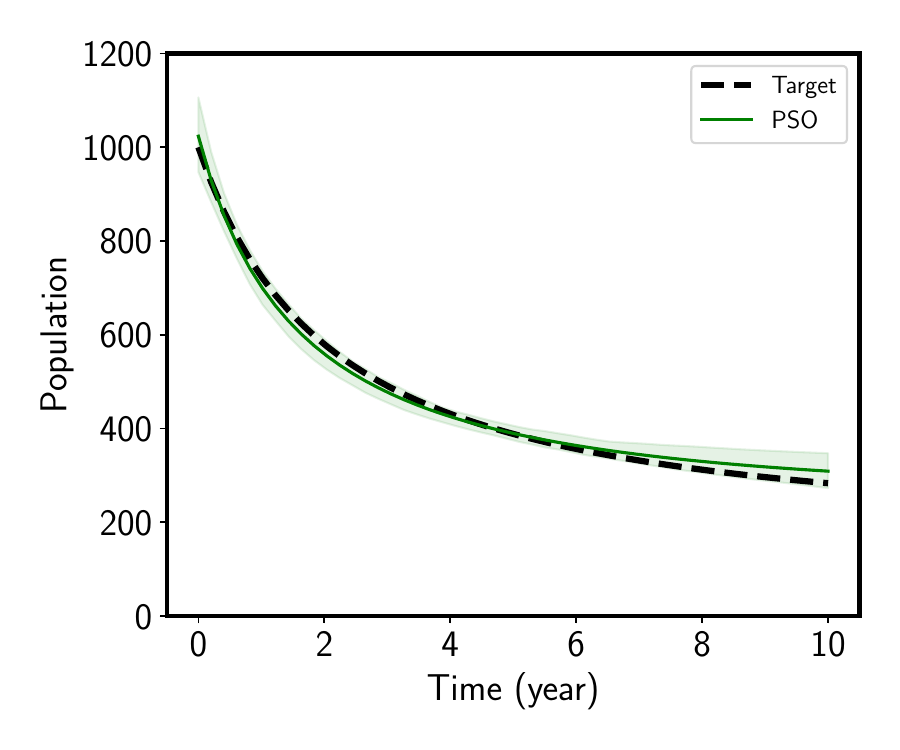}     
   \caption{}
   \label{fig:FigureC2f} 
\end{subfigure}
\hfill
\begin{subfigure}[b]{0.24\textwidth}
   \includegraphics[width=\linewidth]{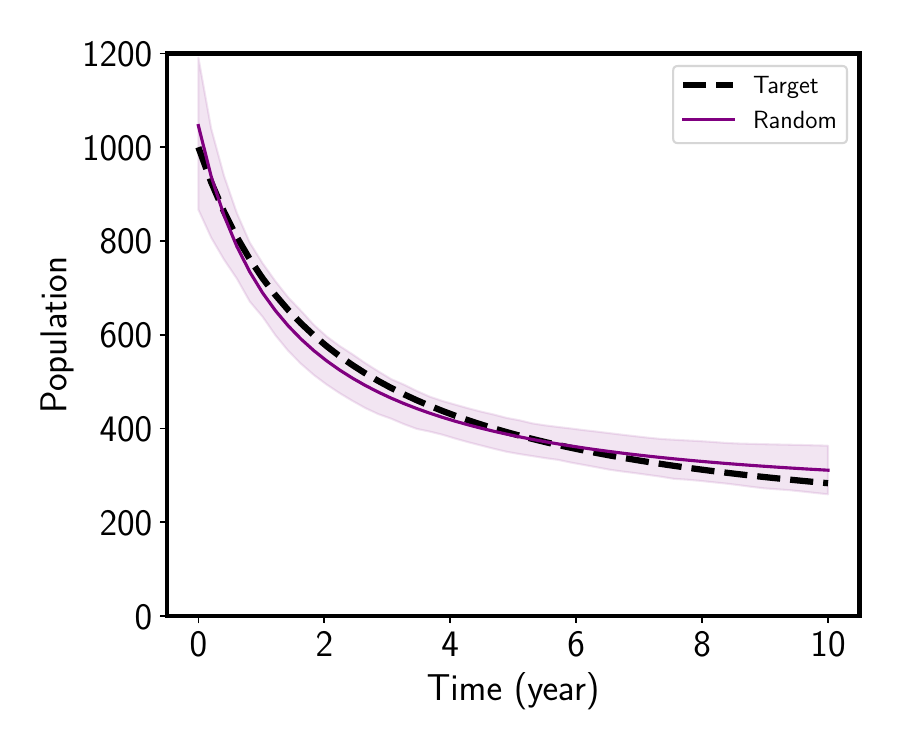}     
   \caption{}
   \label{fig:FigureC2g} 
\end{subfigure}
\hfill
\begin{subfigure}[b]{0.24\textwidth}
   \includegraphics[width=\linewidth]{Figures/Figure2b.pdf}     
   \caption{}
   \label{fig:FigureC2h} 
\end{subfigure}
\end{minipage}
}
\caption{The logistic growth solution reconstruction graphs of all optimization algorithms are presented: (a) BO algorithm. (b) DE algorithm. (c) LBFGSB algorithm. (d) MADS algorithm. (e) NM algorithm. (f) PSO algorithm. (g) Random sampling. (h) ALPS algorithm.}
\label{fig:FigureC2}
\end{figure}

\begin{figure}[!ht]
\centering
\adjustbox{center}{
\begin{minipage}{1.5\textwidth}
\begin{subfigure}[b]{0.24\textwidth}
   \includegraphics[width=\linewidth]{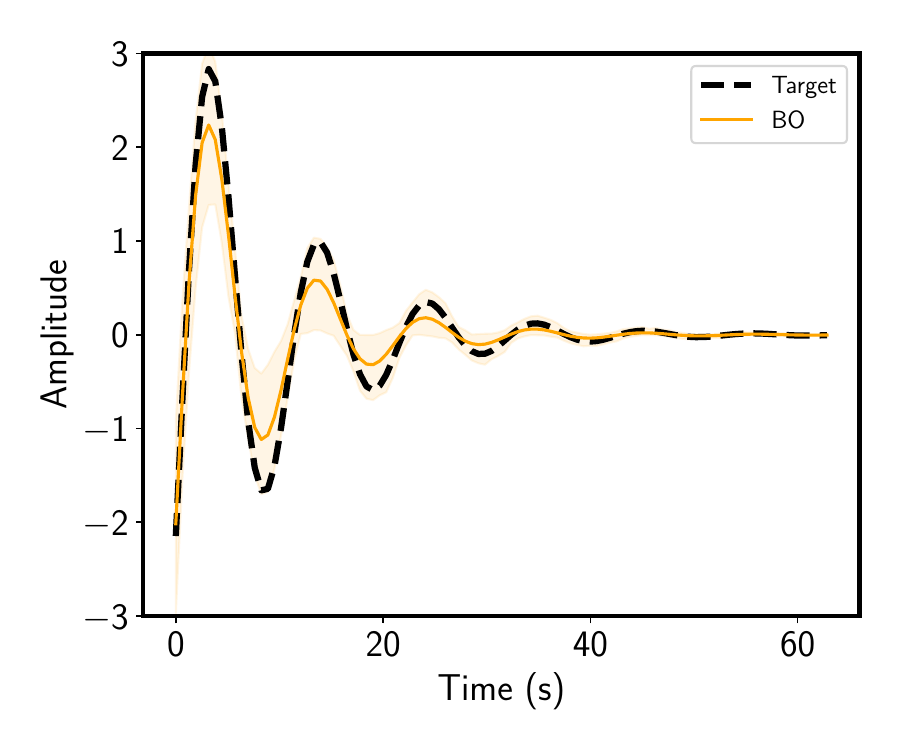}
   \caption{}
   \label{fig:FigureC3a}
\end{subfigure}
\hfill 
\begin{subfigure}[b]{0.24\textwidth}
   \includegraphics[width=\linewidth]{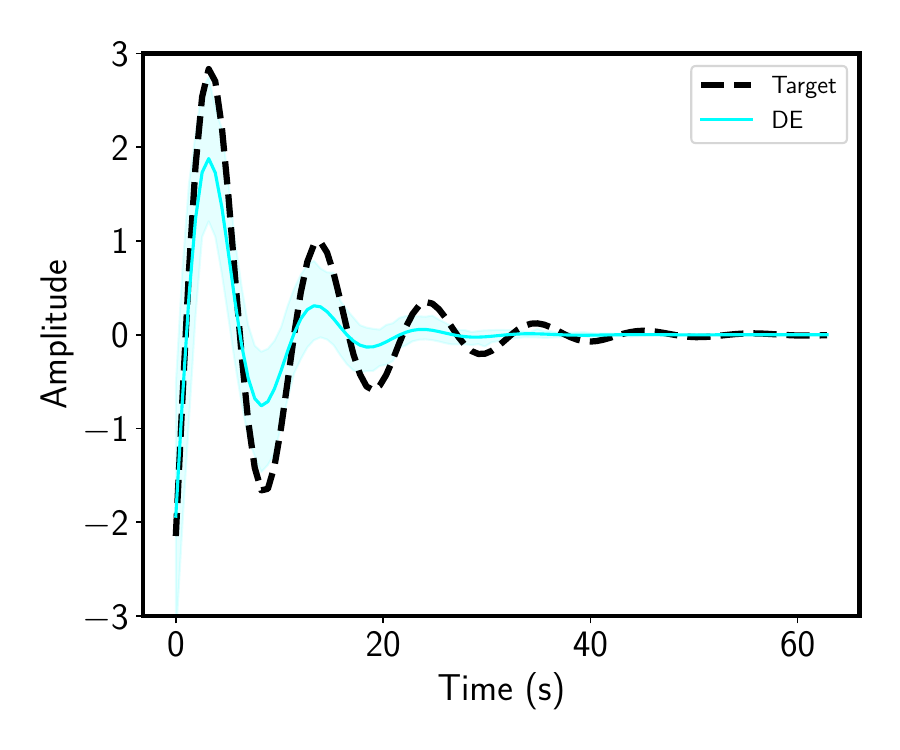}     
   \caption{}
   \label{fig:FigureC3b} 
\end{subfigure}
\hfill
\begin{subfigure}[b]{0.24\textwidth}
   \includegraphics[width=\linewidth]{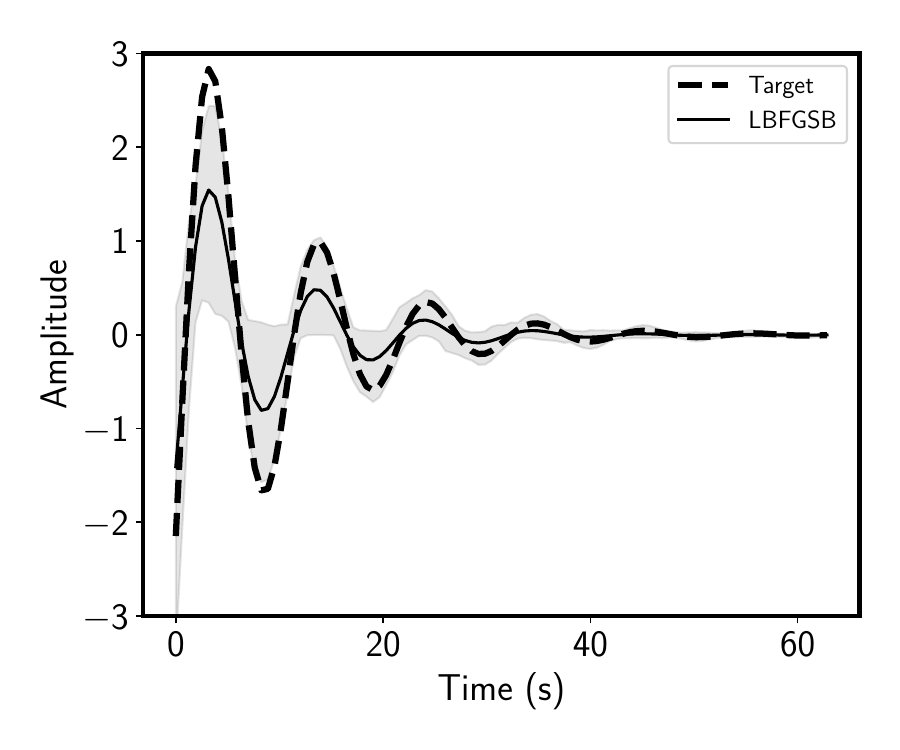}     
   \caption{}
   \label{fig:FigureC3c} 
\end{subfigure}
\hfill
\begin{subfigure}[b]{0.24\textwidth}
   \includegraphics[width=\linewidth]{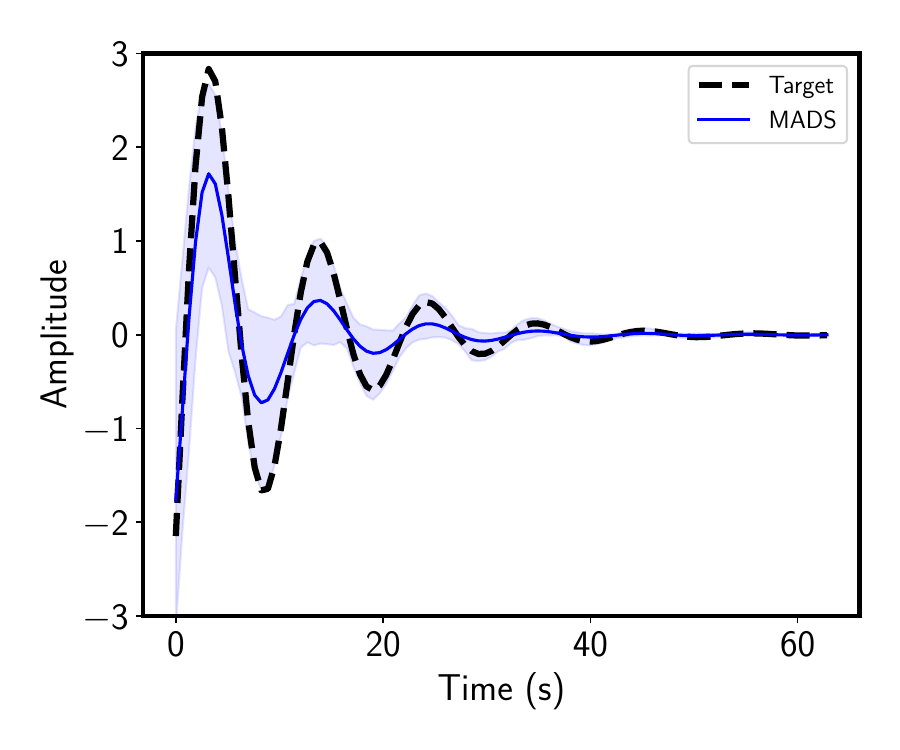}     
   \caption{}
   \label{fig:FigureC3d} 
\end{subfigure}
\end{minipage}
}

\adjustbox{center}{
\begin{minipage}{1.5\textwidth}
\begin{subfigure}[b]{0.24\textwidth}
   \includegraphics[width=\linewidth]{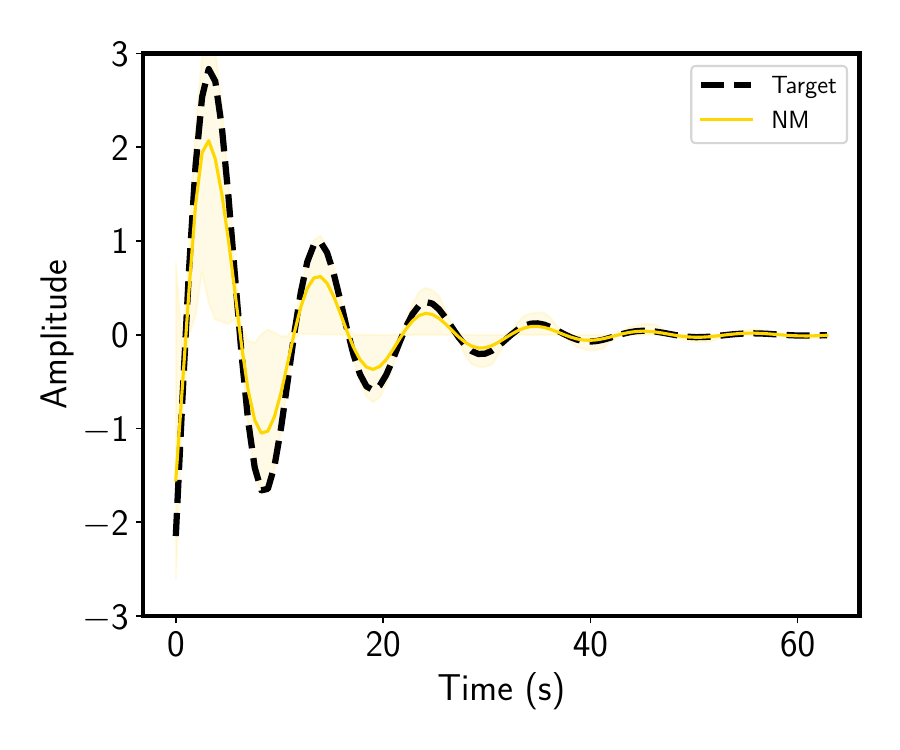}     
   \caption{}
   \label{fig:FigureC3e} 
\end{subfigure}
\hfill
\begin{subfigure}[b]{0.24\textwidth}
   \includegraphics[width=\linewidth]{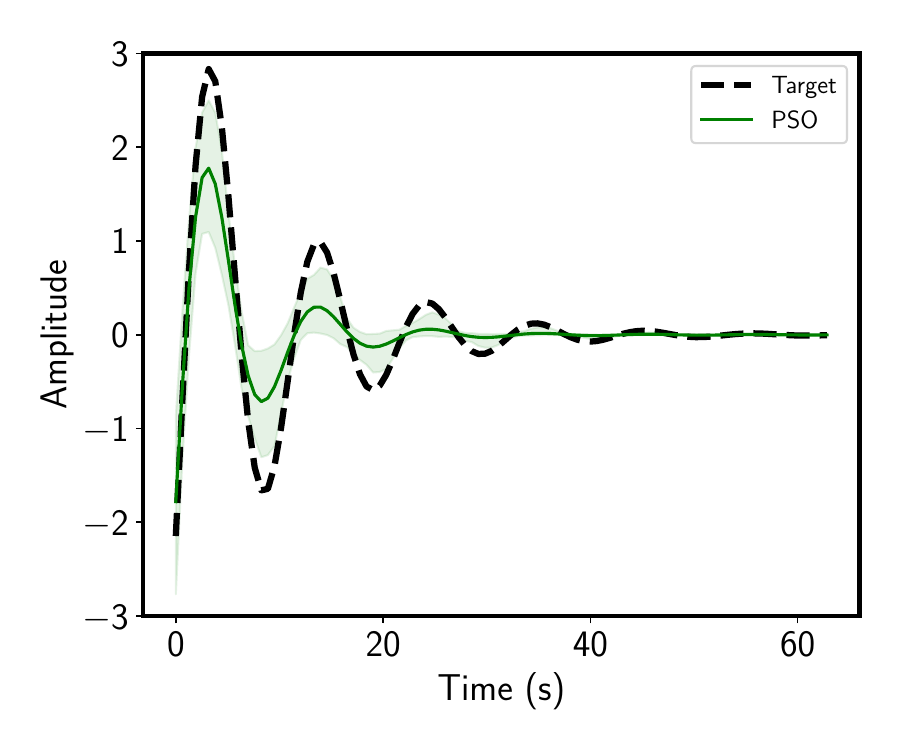}     
   \caption{}
   \label{fig:FigureC3f} 
\end{subfigure}
\hfill
\begin{subfigure}[b]{0.24\textwidth}
   \includegraphics[width=\linewidth]{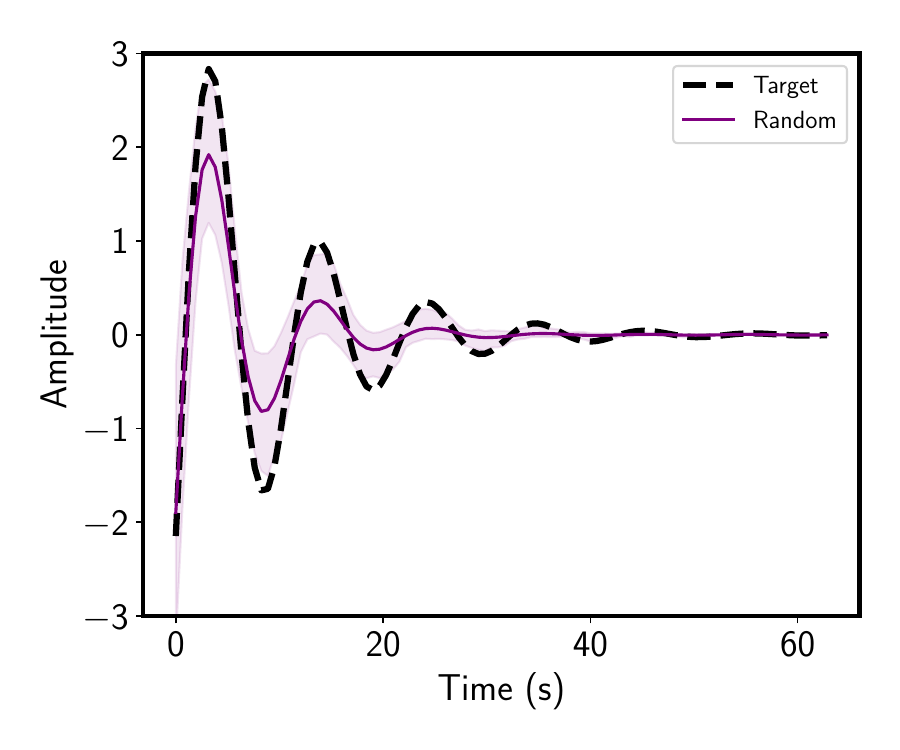}     
   \caption{}
   \label{fig:FigureC3g} 
\end{subfigure}
\hfill
\begin{subfigure}[b]{0.24\textwidth}
   \includegraphics[width=\linewidth]{Figures/Figure2d.pdf}     
   \caption{}
   \label{fig:FigureC3h} 
\end{subfigure}
\end{minipage}
}
\caption{The sinusoidal oscillation with damping solution reconstruction graphs of all optimization algorithms are presented: (a) BO algorithm. (b) DE algorithm. (c) LBFGSB algorithm. (d) MADS algorithm. (e) NM algorithm. (f) PSO algorithm. (g) Random sampling. (h) ALPS algorithm.}
\label{fig:FigureC3}
\end{figure}

This section presents detailed convergence graphs and solution reconstruction graphs for the logistic growth and sinusoidal oscillation with damping benchmarks. Fig. \ref{fig:FigureC1} displays convergence graphs that include uncertainty metrics for both benchmarks across all algorithms. Fig. \ref{fig:FigureC2} and \ref{fig:FigureC3} illustrate the reconstructed solutions obtained by all optimization algorithms, compared with the target design. Tab. \ref{tab:detailed_results_logistic} and \ref{tab:detailed_results_sinusoidal} provide comprehensive convergence statistics from 100 repeated runs for each optimization algorithm.

\newpage

\newpage

\newpage

\begin{table}[!h]
\centering
\caption{Detailed convergence statistics of the $\epsilon$ value for the logistic growth benchmark, for all algorithms. The results are for 100 repeated runs. Smaller values are better.}
\begin{tabular}{ c | c | c | c | c}
 \textbf{Algorithm}  & \textbf{Mean $\epsilon$} &  \textbf{Std $\epsilon$} & \textbf{Maximum $\epsilon$} & \textbf{Minimum $\epsilon$}\\
\hline
\hline
ALPS & $\mathbf{13.52}$ & $\mathbf{9.030}$ & 67.08 & 2.91  \\
MADS & 28.40& 14.00 & $\mathbf{57.93}$ &1.50  \\
BO & 38.79& 15.91  & 97.68 &  7.34 \\
PSO & 26.20&16.14&108.25& 3.23 \\
DE &41.77 & 16.38 & 112.56 &  4.72 \\
LBFGSB& 34.70&34.83 &114.52 &0.22  \\
Random &42.96&15.99 &  100.55 &14.71  \\
NM & 26.30 & 28.95 & 131.88& $\mathbf{0.10}$  \\
\end{tabular}
\label{tab:detailed_results_logistic}
\end{table}

\begin{table}[!h]
\centering
\caption{Detailed convergence statistics of the $\epsilon$ value for the sinusoidal oscillation with damping growth benchmark, for all algorithms. The results are for 100 repeated runs. Smaller values are better.}
\begin{tabular}{ c | c | c | c | c}
 \textbf{Algorithm}  & \textbf{Mean $\epsilon$} &  \textbf{Std $\epsilon$} & \textbf{Maximum $\epsilon$} & \textbf{Minimum $\epsilon$}\\
\hline
\hline
ALPS & $\mathbf{0.25}$ & 0.11 &0.62& 0.08  \\
MADS &0.46 &0.19 &0.72& 0.06 \\
BO& 0.31  &0.13 &0.63& 0.07 \\ 
PSO &0.40 &0.11 &$\mathbf{0.60}$& 0.10  \\ 
DE& 0.44 & $\mathbf{0.08}$ & 0.63 &0.23  \\ 
LBFGSB& 0.44 &0.19 &0.79& 0.05  \\
Random& 0.42 &0.10& 0.62& 0.16  \\
NM& 0.30&0.26 &0.79& $\mathbf{0.01}$  \\
\end{tabular}

\label{tab:detailed_results_sinusoidal}
\end{table}

\newpage
\subsection{Photonic Surface Inverse Design Convergence}\label{subapp:photonicresults}

\begin{figure}[!h]
\centering
\begin{subfigure}[b]{0.49\textwidth}
   \includegraphics[width=\linewidth]{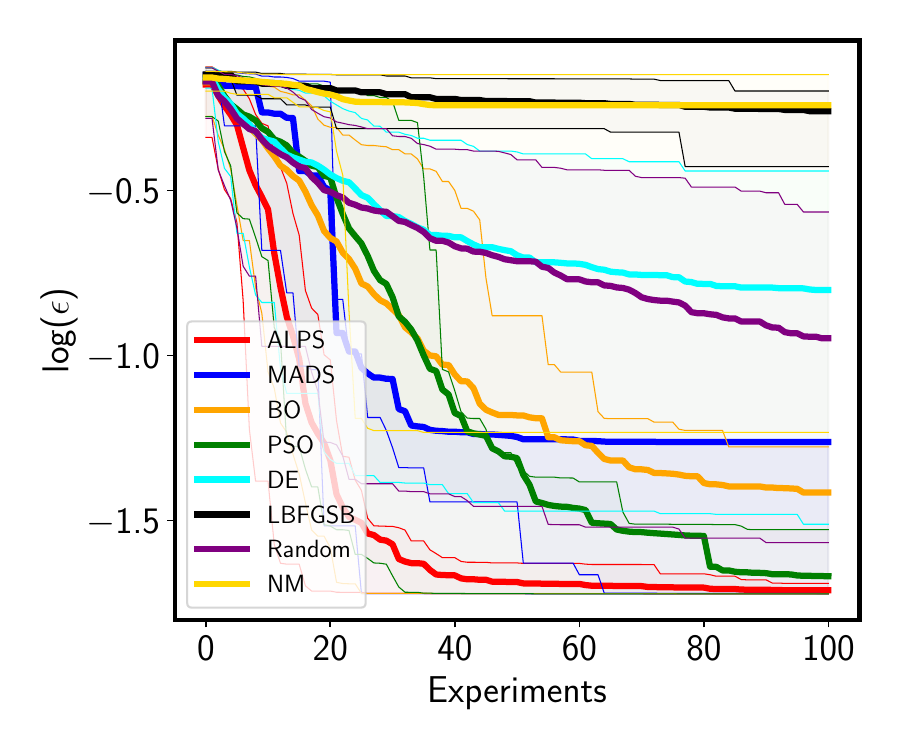}
   \caption{}
   \label{fig:FigureC4a}
\end{subfigure}
\begin{subfigure}[b]{0.49\textwidth}
   \includegraphics[width=\linewidth]{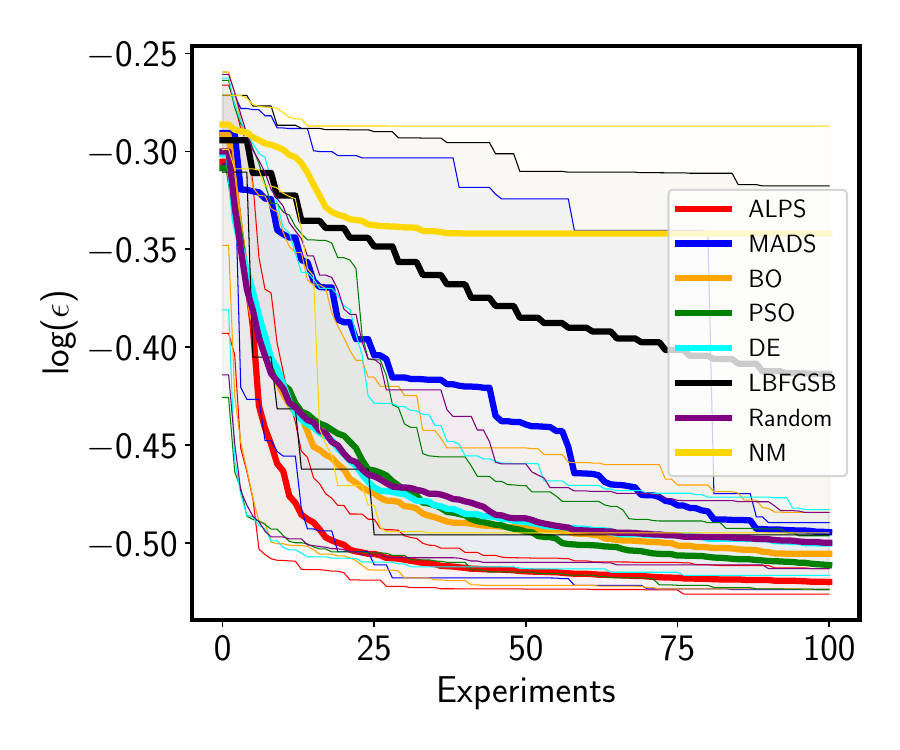}     
   \caption{}
   \label{fig:FigureC4b} 
\end{subfigure}
\begin{subfigure}[b]{0.49\textwidth}
   \includegraphics[width=\linewidth]{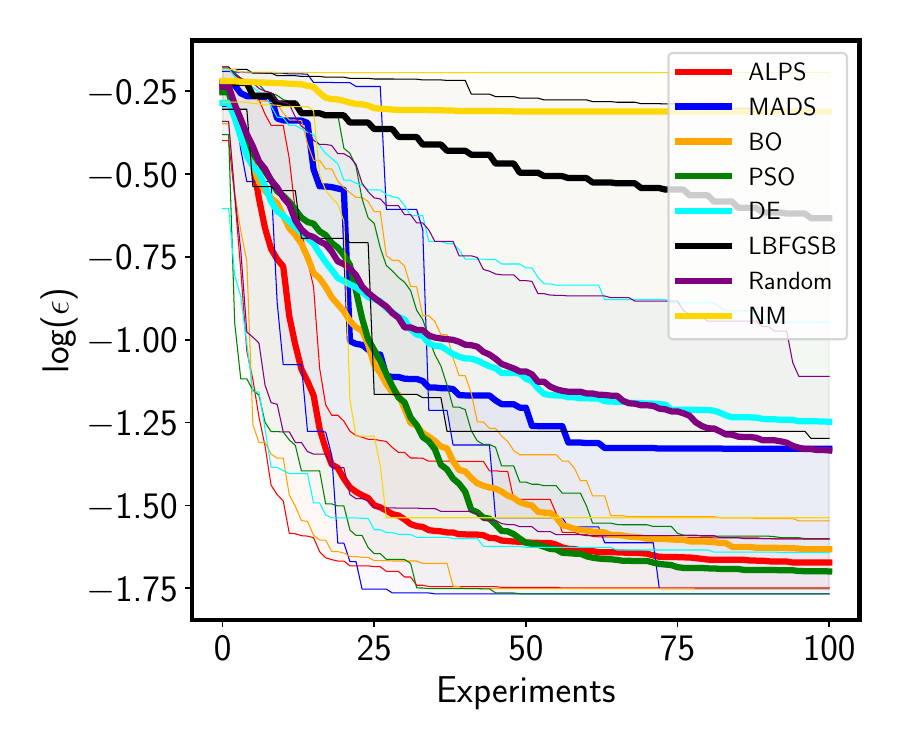}
   \caption{}
   \label{fig:FigureC4c}
\end{subfigure}
\begin{subfigure}[b]{0.49\textwidth}
   \includegraphics[width=\linewidth]{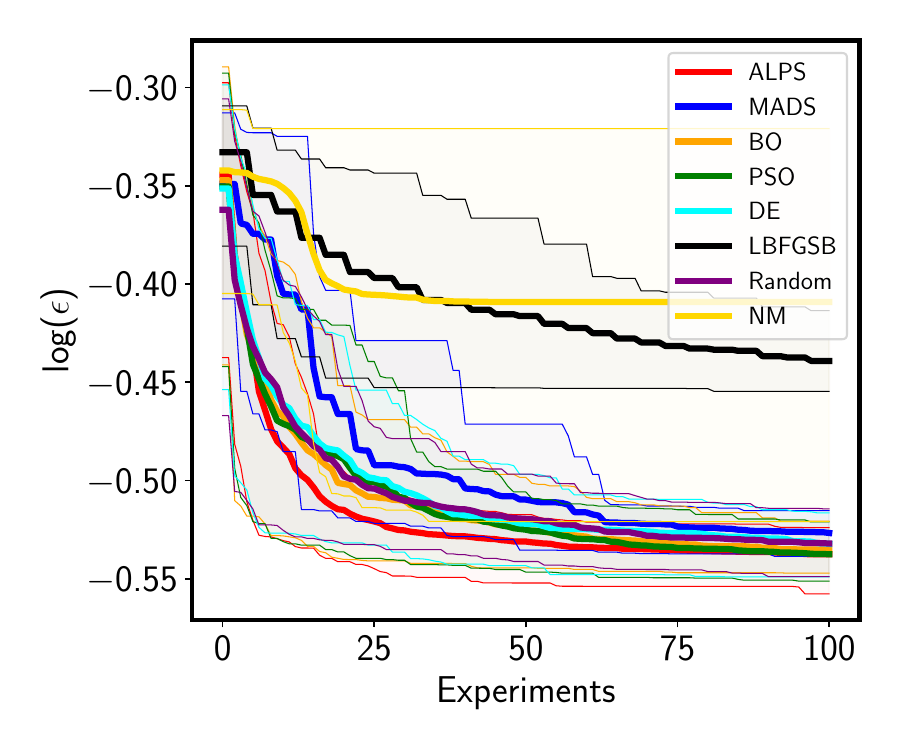}
   \caption{}
   \label{fig:FigureC4d}
\end{subfigure}
\caption[]{Convergence graphs display the base-10 logarithm of the mean error value (Eq. (\ref{eqn:inverse_design_definition})) for the photonic surface inverse design benchmarks, for all algorithms: (a) Inconel near-perfect emitter target. (b) Inconel TPV emitter target. (c) Stainless steel near-perfect emitter target. (d) Stainless steel TPV emitter target.}
\label{fig:FigureC4}
\end{figure}

\begin{figure}[H]
\centering
\adjustbox{center}{
\begin{minipage}{1.5\textwidth}
\begin{subfigure}[b]{0.24\textwidth}
   \includegraphics[width=\linewidth]{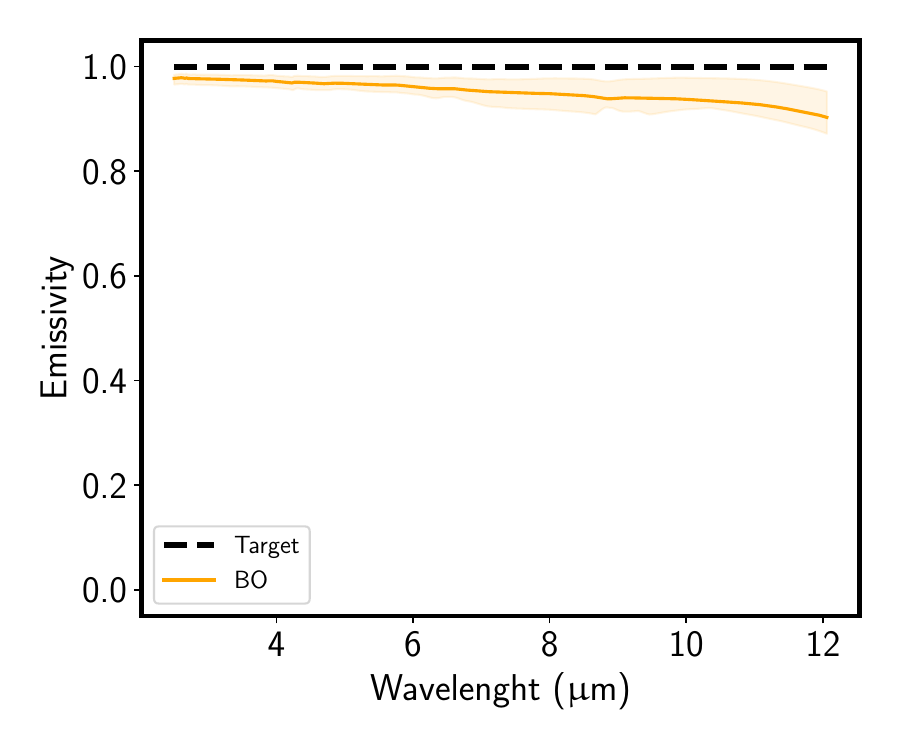}
   \caption{}
   \label{fig:FigureC5a}
\end{subfigure}
\hfill
\begin{subfigure}[b]{0.24\textwidth}
   \includegraphics[width=\linewidth]{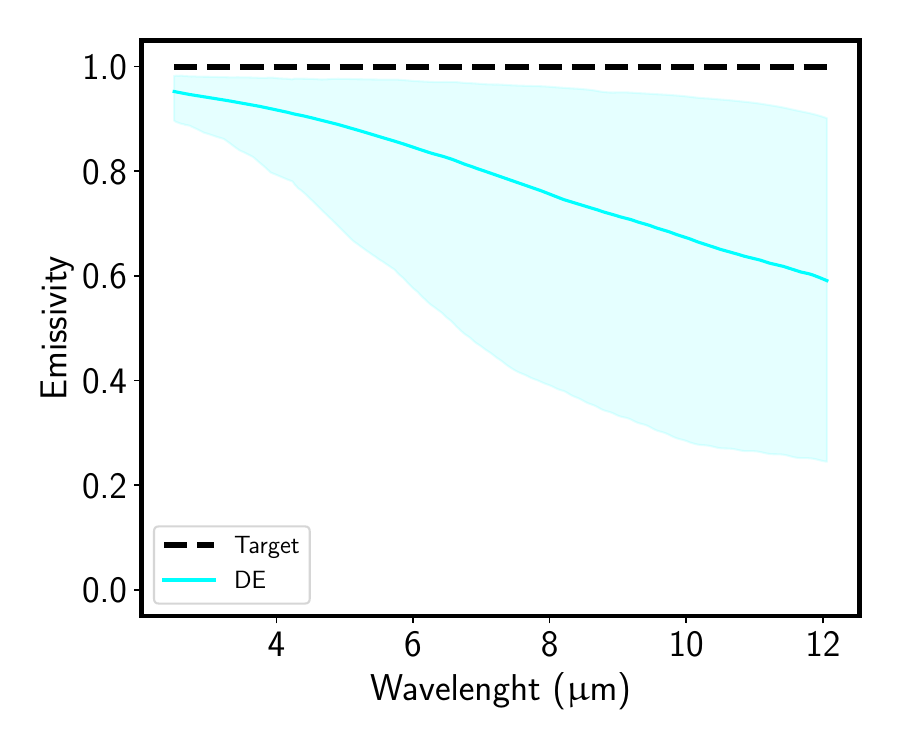}     
   \caption{}
   \label{fig:FigureC5b} 
\end{subfigure}
\hfill
\begin{subfigure}[b]{0.24\textwidth}
   \includegraphics[width=\linewidth]{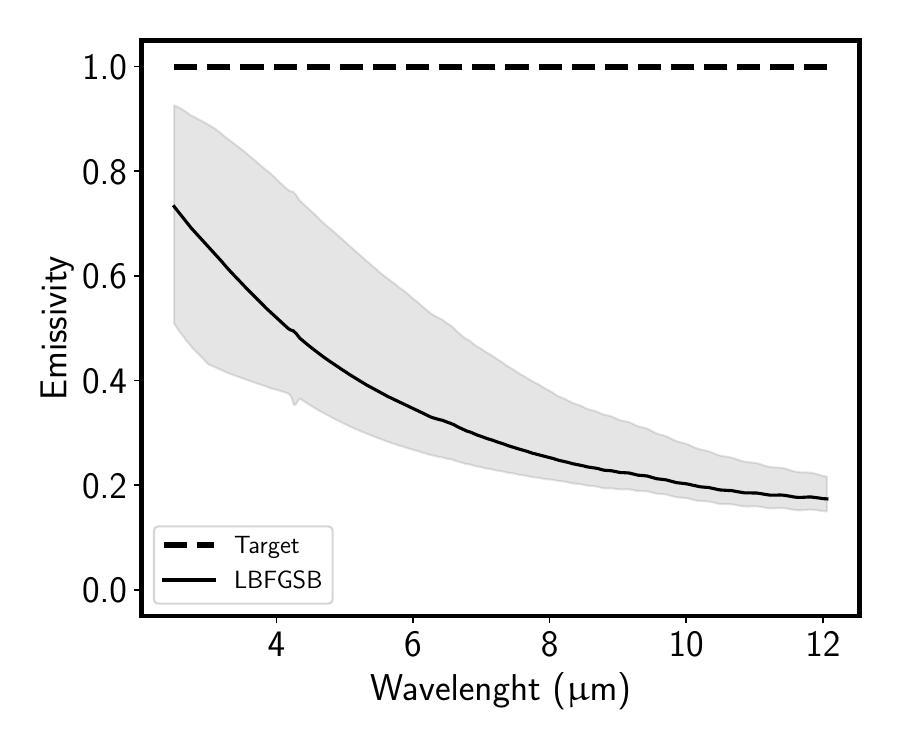}     
   \caption{}
   \label{fig:FigureC5c} 
\end{subfigure}
\hfill
\begin{subfigure}[b]{0.24\textwidth}
   \includegraphics[width=\linewidth]{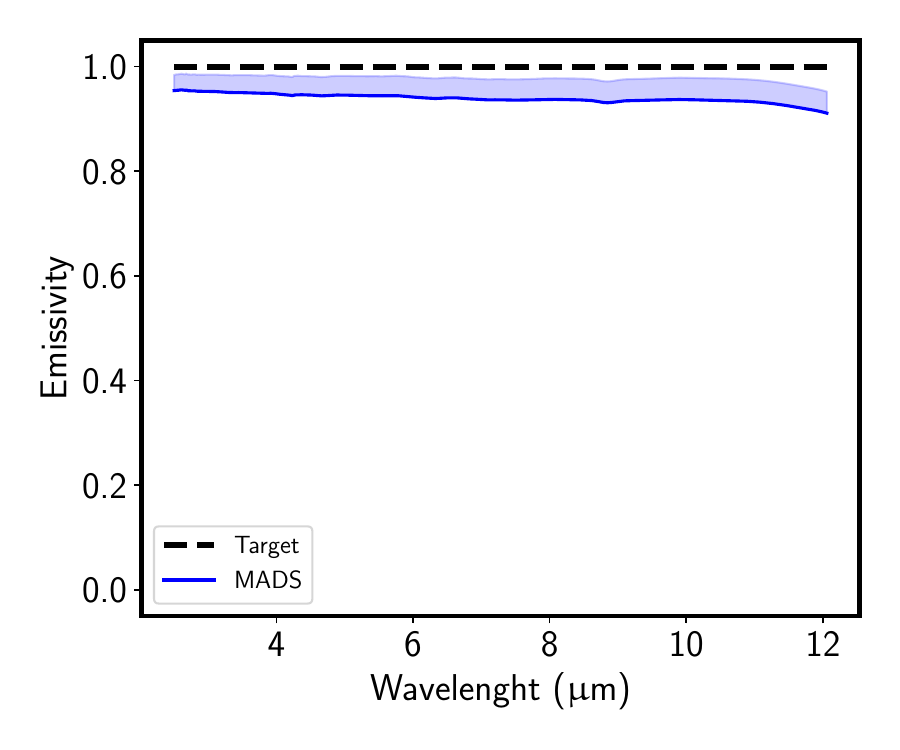}     
   \caption{}
   \label{fig:FigureC5d} 
\end{subfigure}
\end{minipage}
}

\adjustbox{center}{
\begin{minipage}{1.5\textwidth}
\begin{subfigure}[b]{0.24\textwidth}
   \includegraphics[width=\linewidth]{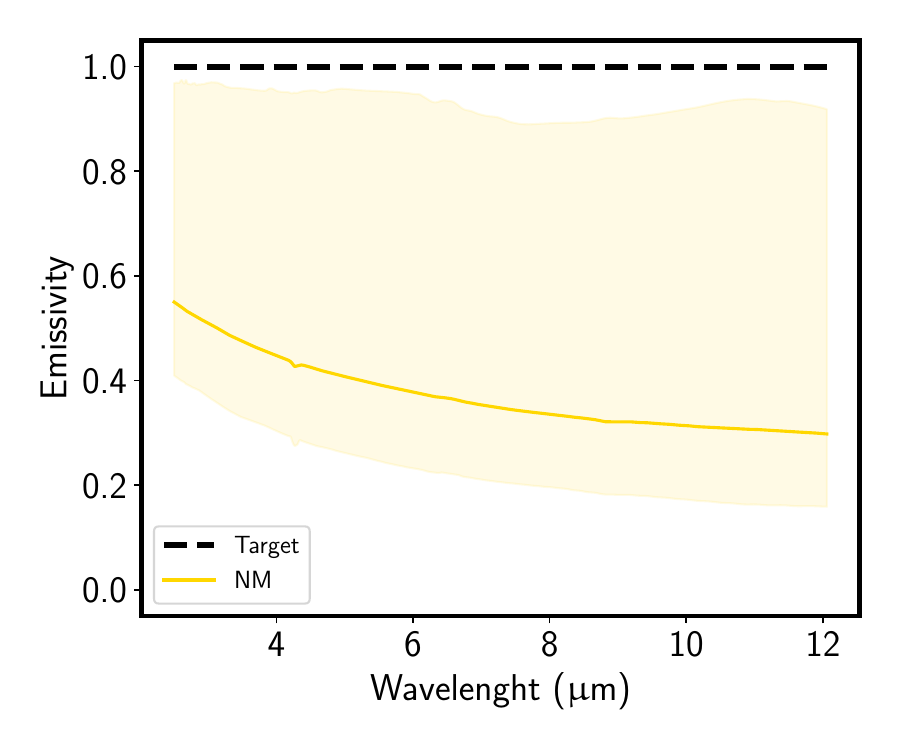}     
   \caption{}
   \label{fig:FigureC5e} 
\end{subfigure}
\hfill
\begin{subfigure}[b]{0.24\textwidth}
   \includegraphics[width=\linewidth]{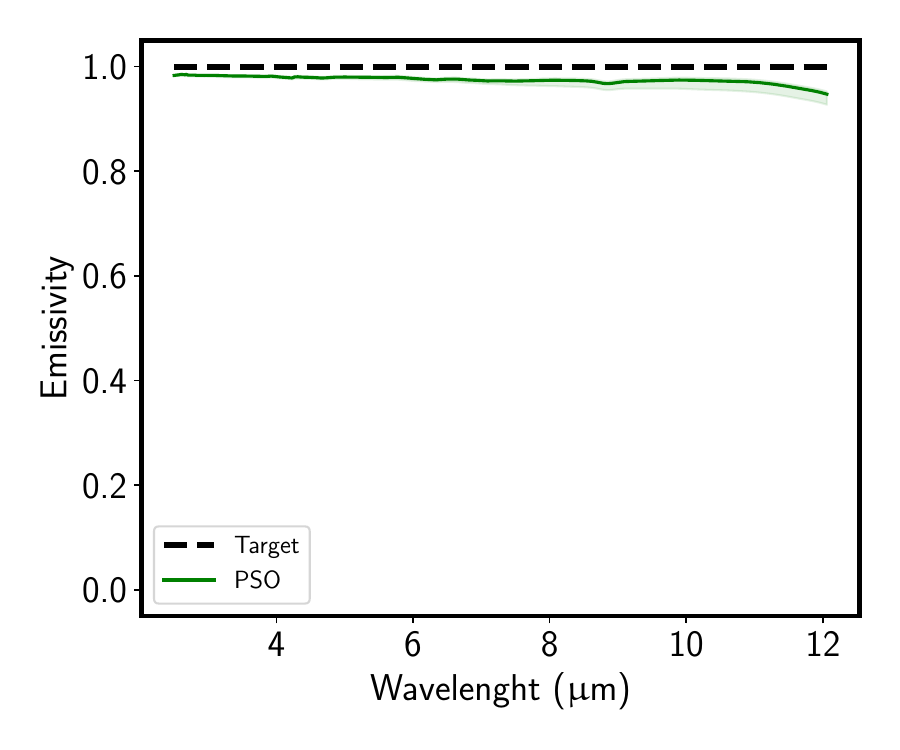}     
   \caption{}
   \label{fig:FigureC5f} 
\end{subfigure}
\hfill
\begin{subfigure}[b]{0.24\textwidth}
   \includegraphics[width=\linewidth]{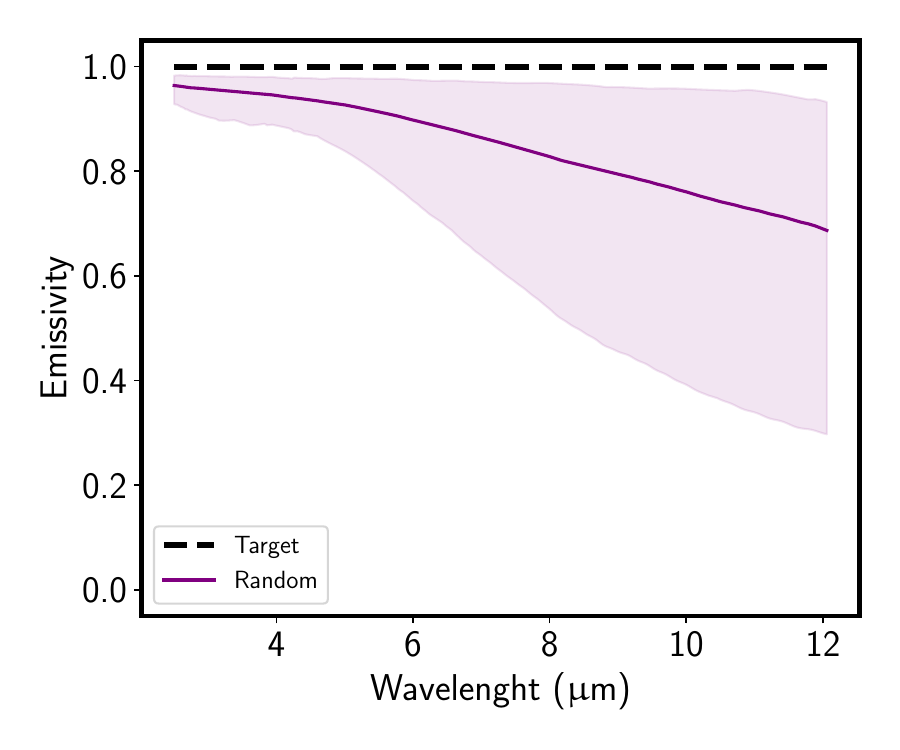}     
   \caption{}
   \label{fig:FigureC5g} 
\end{subfigure}
\hfill
\begin{subfigure}[b]{0.24\textwidth}
   \includegraphics[width=\linewidth]{Figures/Figure3e.pdf}     
   \caption{}
   \label{fig:FigureC5h} 
\end{subfigure}
\end{minipage}
}
\caption{The solution reconstruction graphs for the Inconel near-perfect emitter target for all optimization algorithms: (a) BO algorithm. (b) DE algorithm. (c) LBFGSB algorithm. (d) MADS algorithm. (e) NM algorithm. (f) PSO algorithm. (g) Random sampling. (h) ALPS algorithm.}
\label{fig:FigureC5}
\end{figure}

\begin{figure}[H]
\centering
\adjustbox{center}{
\begin{minipage}{1.5\textwidth}
\begin{subfigure}[b]{0.24\textwidth}
   \includegraphics[width=\linewidth]{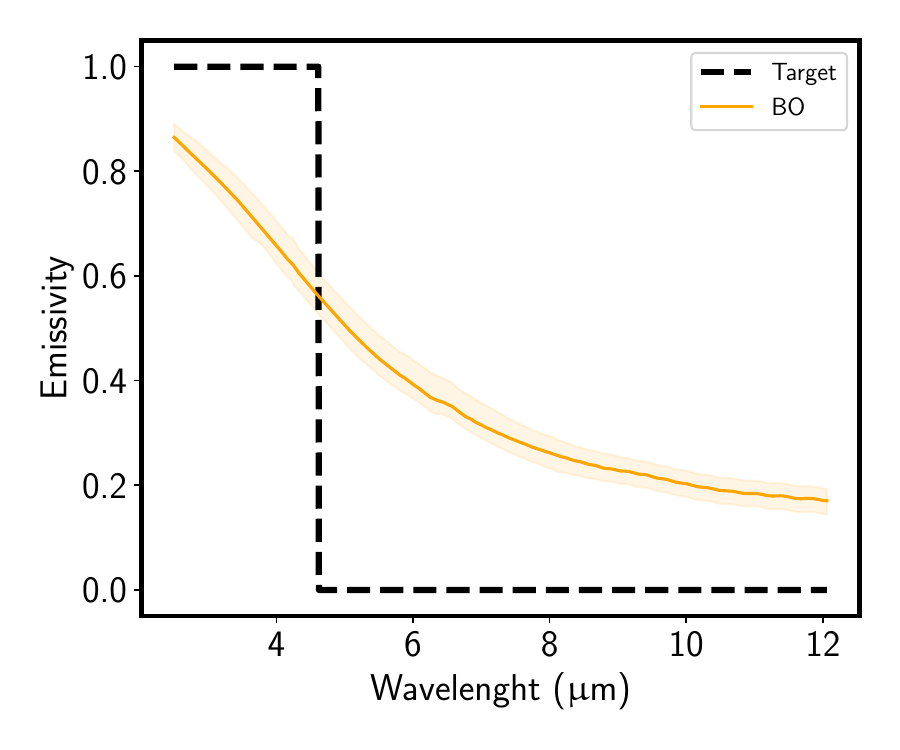}
   \caption{}
   \label{fig:FigureC6a}
\end{subfigure}
\hfill
\begin{subfigure}[b]{0.24\textwidth}
   \includegraphics[width=\linewidth]{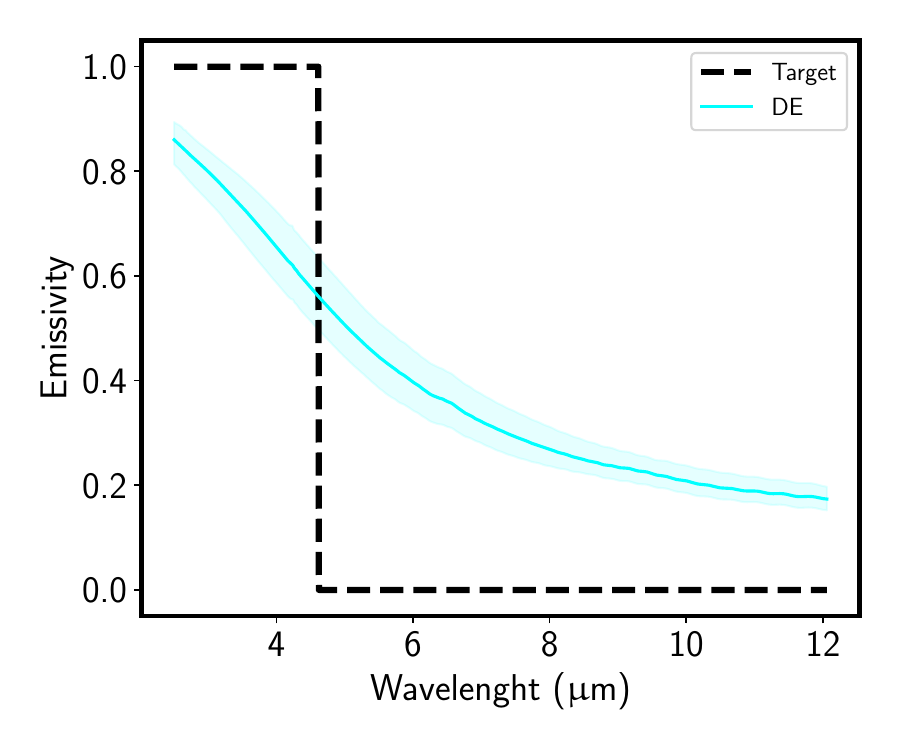}     
   \caption{}
   \label{fig:FigureC6b} 
\end{subfigure}
\hfill
\begin{subfigure}[b]{0.24\textwidth}
   \includegraphics[width=\linewidth]{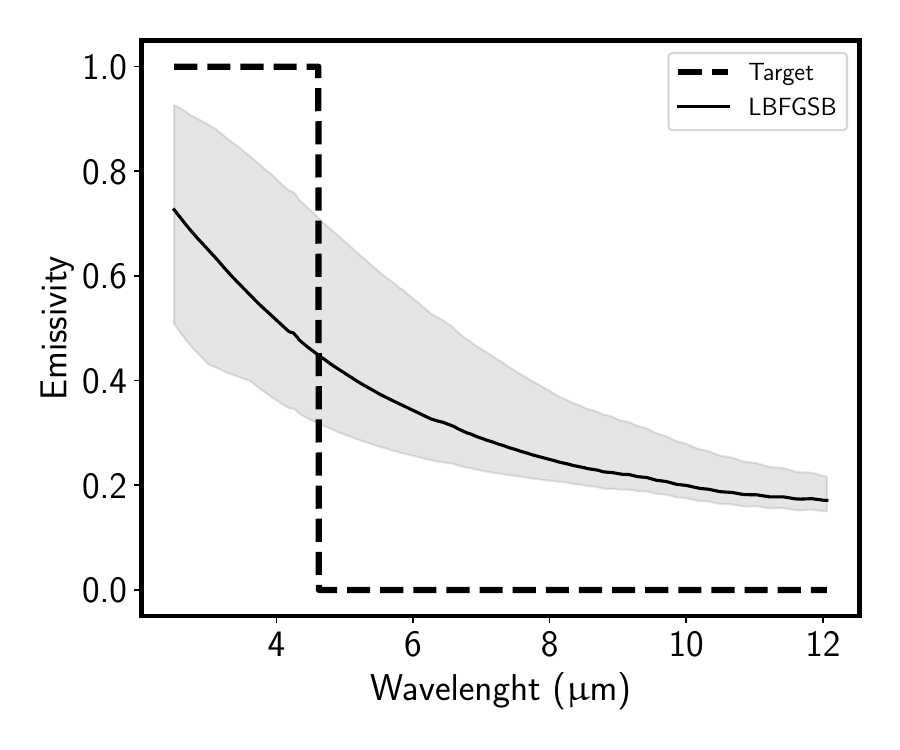}     
   \caption{}
   \label{fig:FigureC6c} 
\end{subfigure}
\hfill
\begin{subfigure}[b]{0.24\textwidth}
   \includegraphics[width=\linewidth]{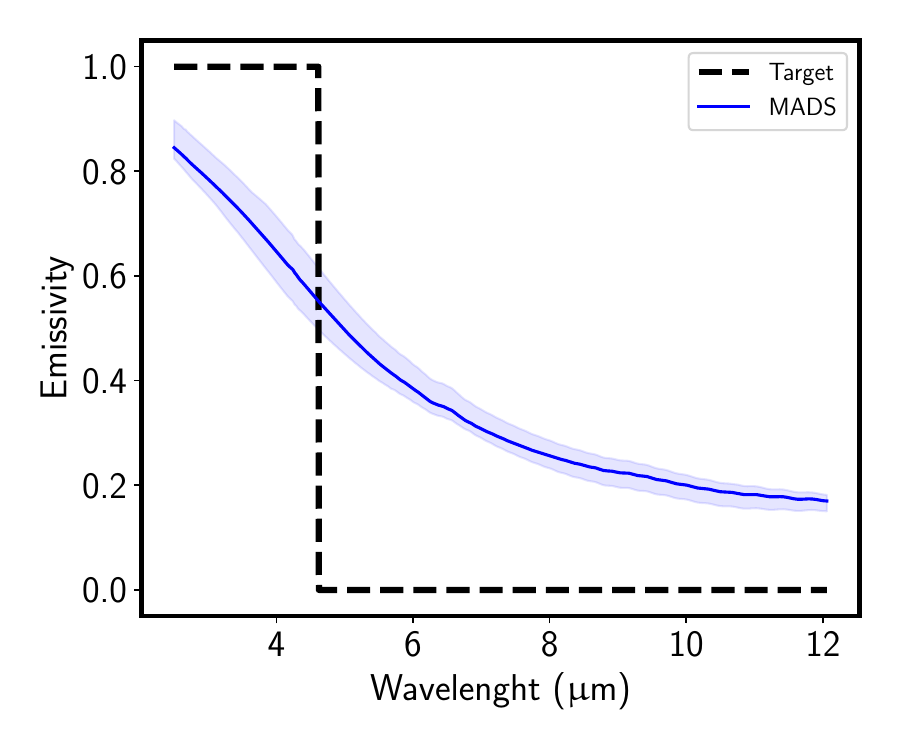}     
   \caption{}
   \label{fig:FigureC6d} 
\end{subfigure}
\end{minipage}
}

\adjustbox{center}{
\begin{minipage}{1.5\textwidth}
\begin{subfigure}[b]{0.24\textwidth}
   \includegraphics[width=\linewidth]{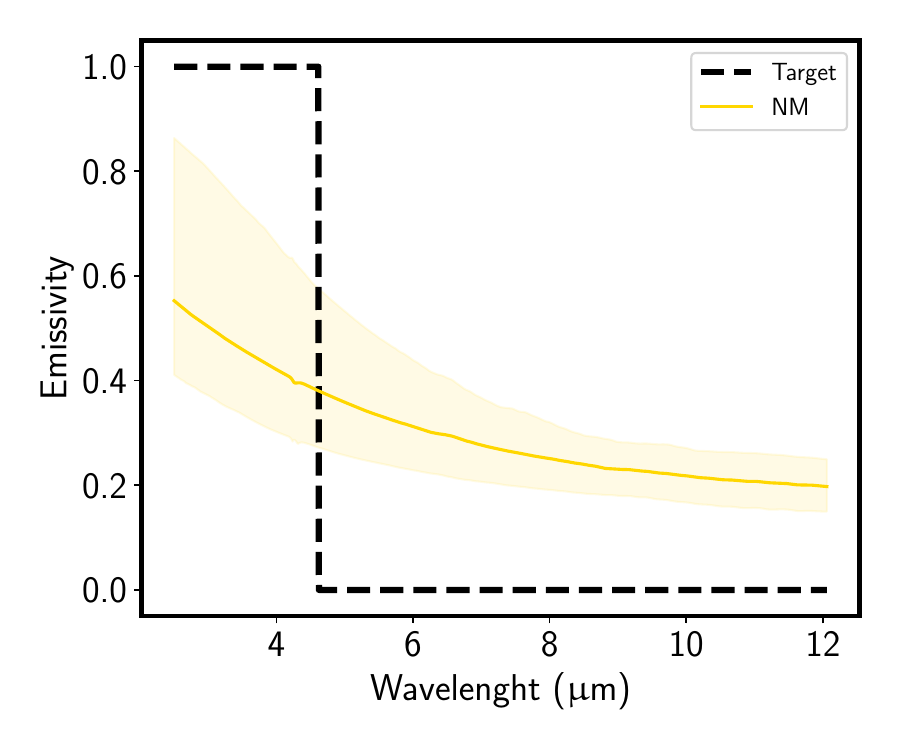}     
   \caption{}
   \label{fig:FigureC6e} 
\end{subfigure}
\hfill
\begin{subfigure}[b]{0.24\textwidth}
   \includegraphics[width=\linewidth]{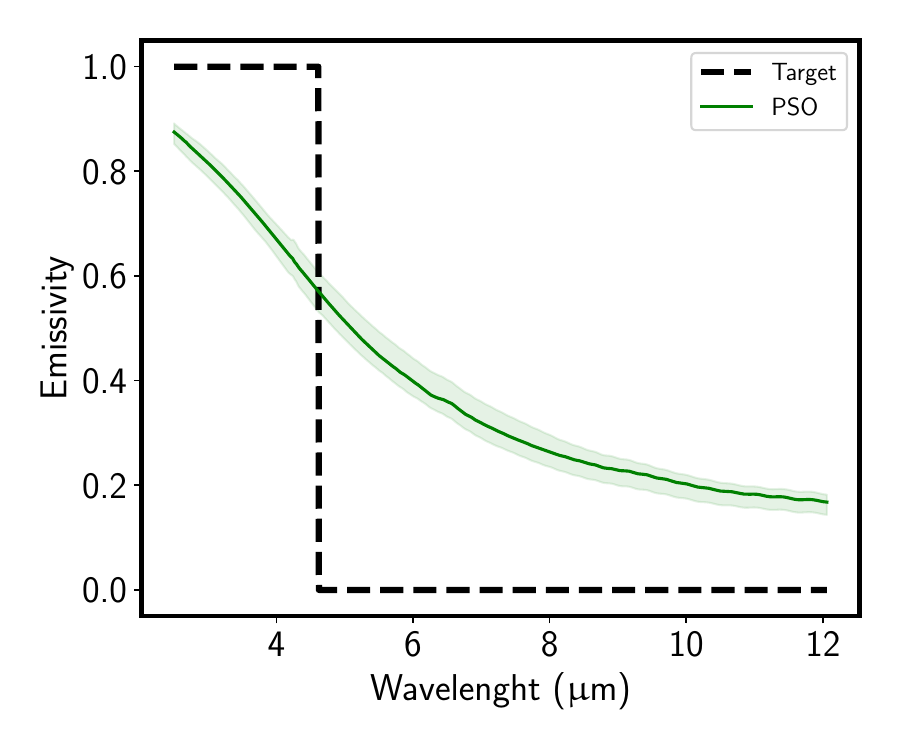}     
   \caption{}
   \label{fig:FigureC6f} 
\end{subfigure}
\hfill
\begin{subfigure}[b]{0.24\textwidth}
   \includegraphics[width=\linewidth]{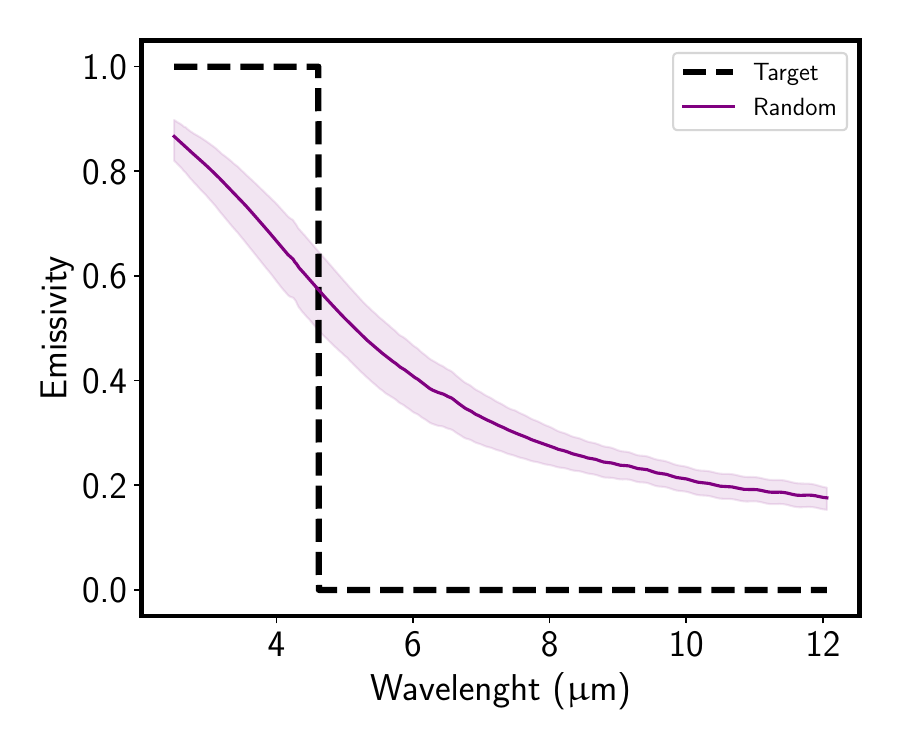}     
   \caption{}
   \label{fig:FigureC6g} 
\end{subfigure}
\hfill
\begin{subfigure}[b]{0.24\textwidth}
   \includegraphics[width=\linewidth]{Figures/Figure3f.pdf}     
   \caption{}
   \label{fig:FigureC6h} 
\end{subfigure}
\end{minipage}
}
\caption{The solution reconstruction graphs of the Inconel TPV emitter target for all optimization algorithms: (a) BO algorithm. (b) DE algorithm. (c) LBFGSB algorithm. (d) MADS algorithm. (e) NM algorithm. (f) PSO algorithm. (g) Random sampling. (h) ALPS algorithm.}
\label{fig:FigureC6}
\end{figure}

\begin{figure}[H]
\centering
\adjustbox{center}{
\begin{minipage}{1.5\textwidth}
\begin{subfigure}[b]{0.24\textwidth}
   \includegraphics[width=\linewidth]{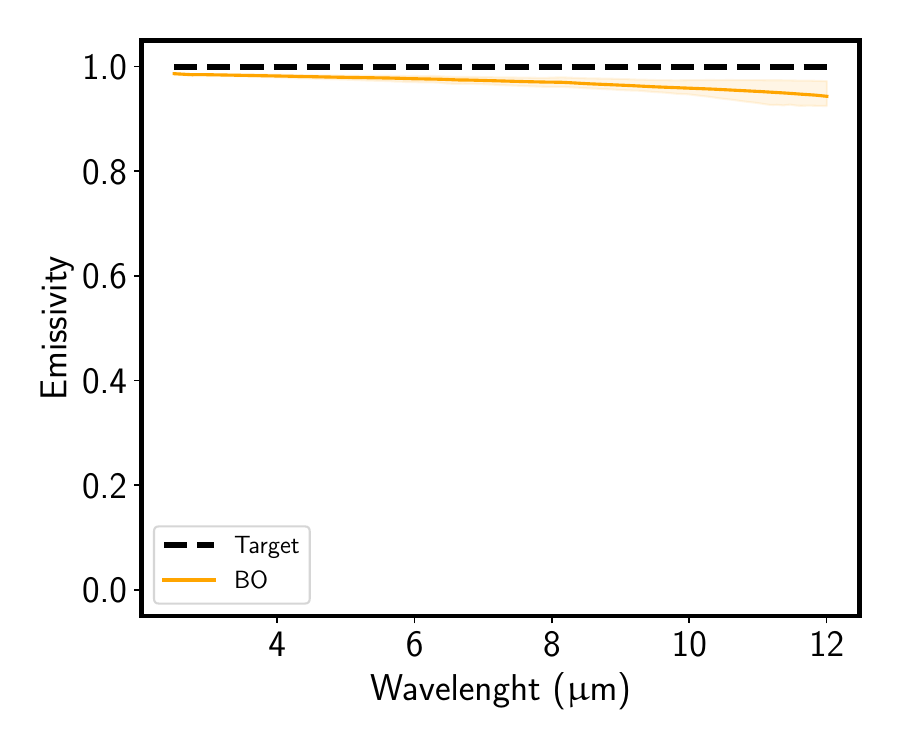}
   \caption{}
   \label{fig:FigureC7a}
\end{subfigure}
\hfill
\begin{subfigure}[b]{0.24\textwidth}
   \includegraphics[width=\linewidth]{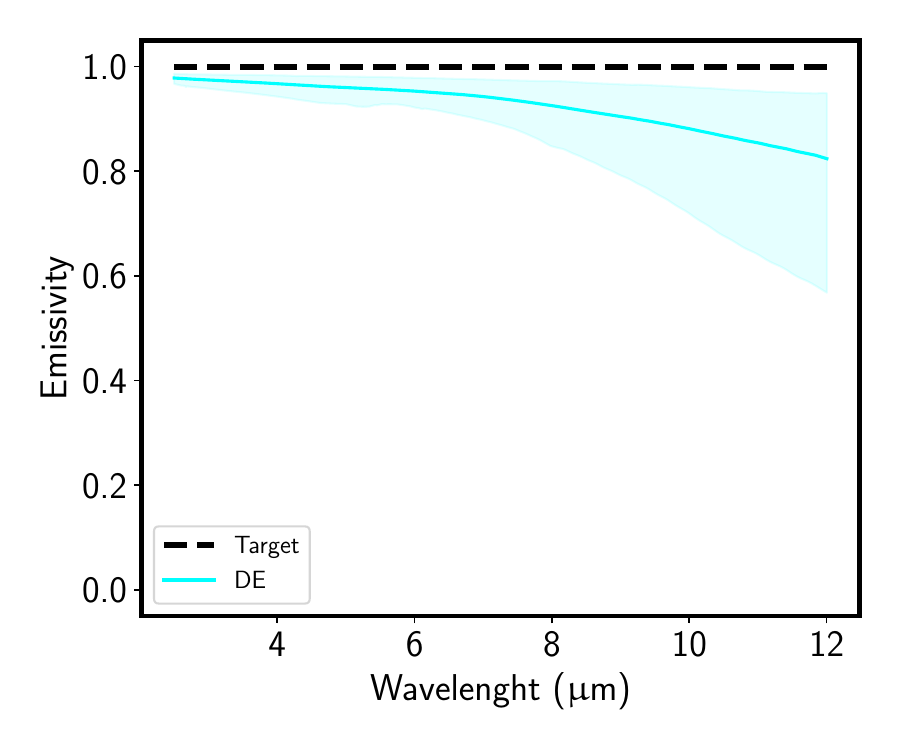}     
   \caption{}
   \label{fig:FigureC7b} 
\end{subfigure}
\hfill
\begin{subfigure}[b]{0.24\textwidth}
   \includegraphics[width=\linewidth]{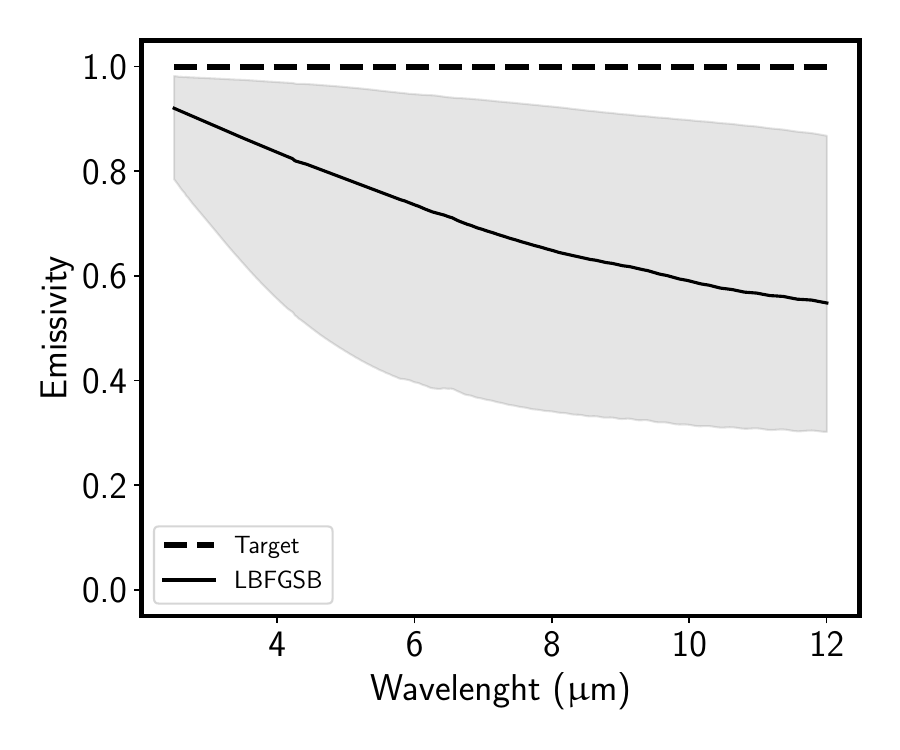}     
   \caption{}
   \label{fig:FigureC7c} 
\end{subfigure}
\hfill
\begin{subfigure}[b]{0.24\textwidth}
   \includegraphics[width=\linewidth]{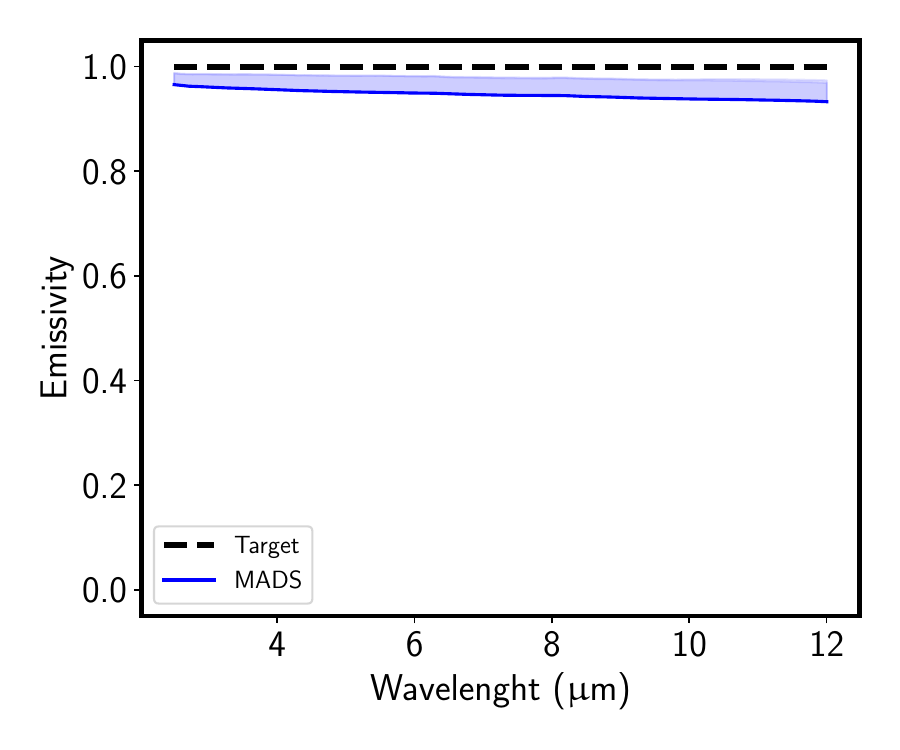}     
   \caption{}
   \label{fig:FigureC7d} 
\end{subfigure}
\end{minipage}
}

\adjustbox{center}{
\begin{minipage}{1.5\textwidth}
\begin{subfigure}[b]{0.24\textwidth}
   \includegraphics[width=\linewidth]{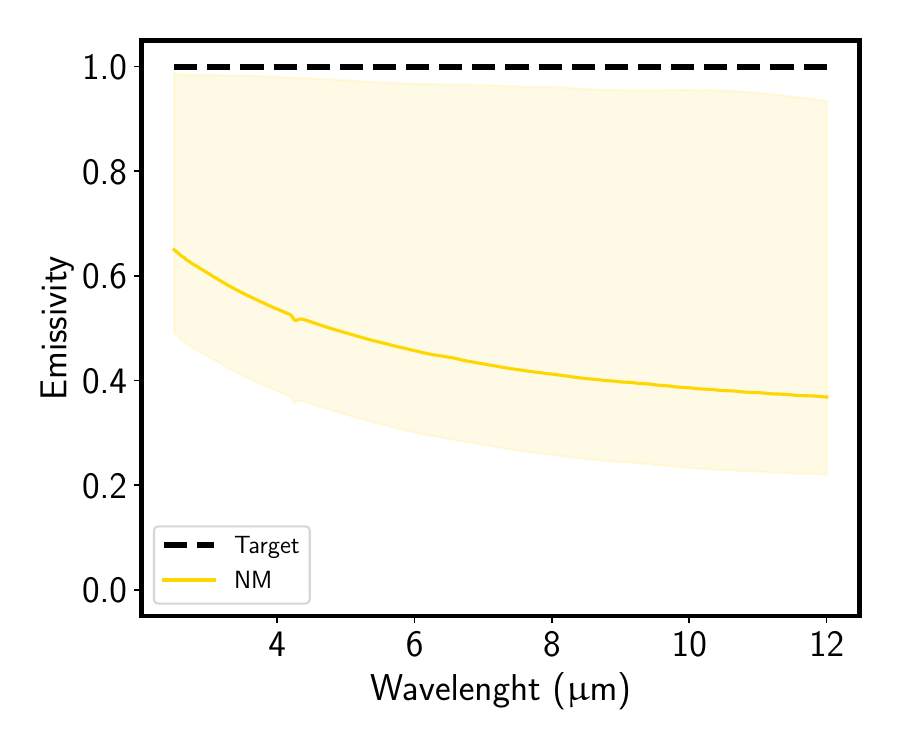}     
   \caption{}
   \label{fig:FigureC7e} 
\end{subfigure}
\hfill
\begin{subfigure}[b]{0.24\textwidth}
   \includegraphics[width=\linewidth]{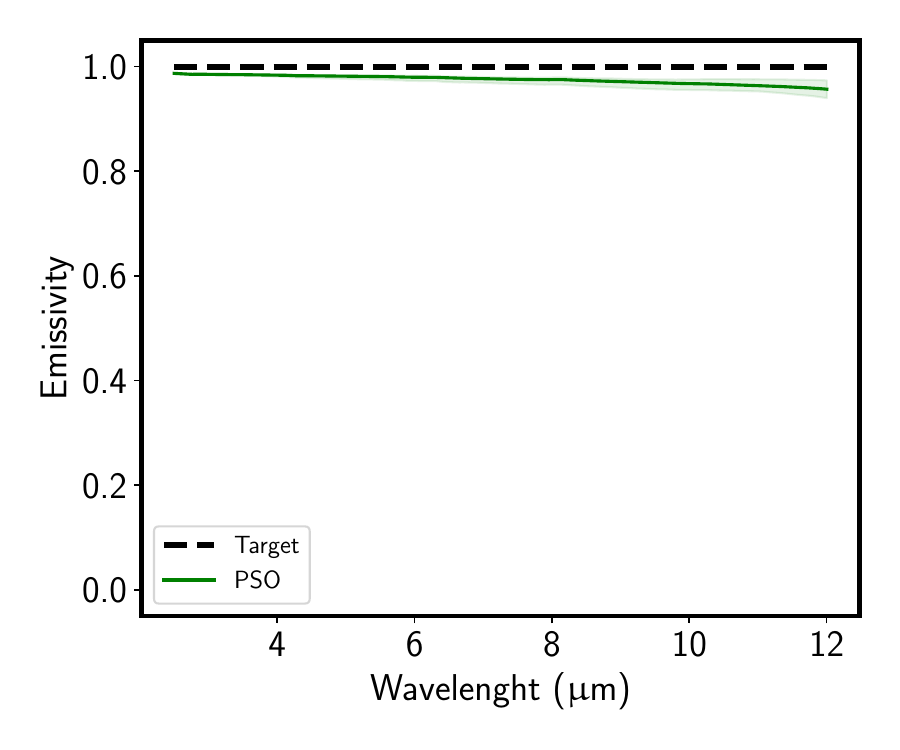}     
   \caption{}
   \label{fig:FigureC7f} 
\end{subfigure}
\hfill
\begin{subfigure}[b]{0.24\textwidth}
   \includegraphics[width=\linewidth]{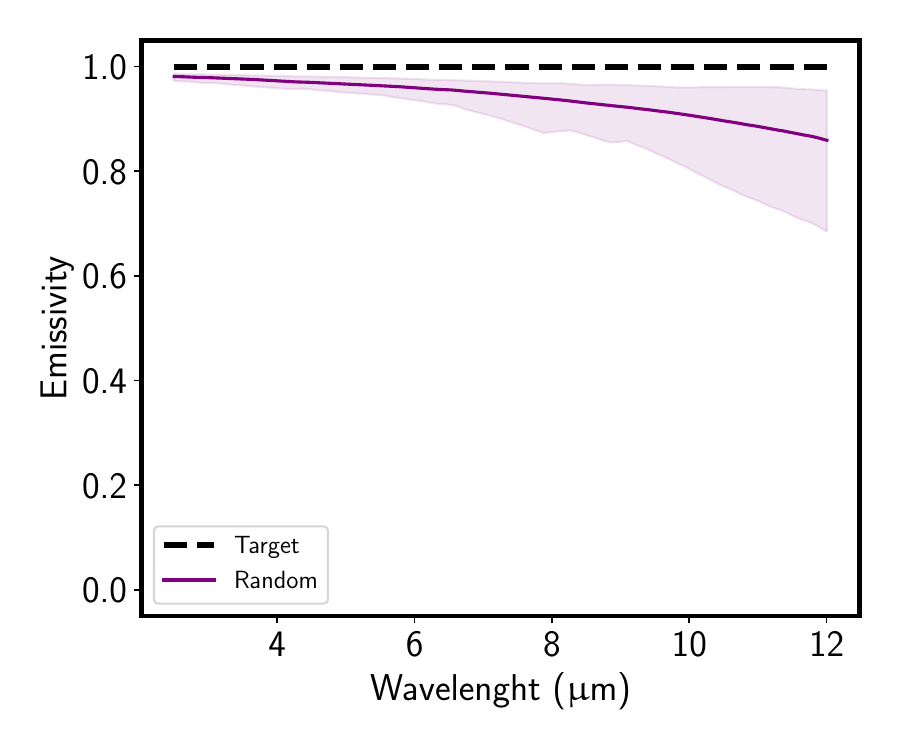}     
   \caption{}
   \label{fig:FigureC7g} 
\end{subfigure}
\hfill
\begin{subfigure}[b]{0.24\textwidth}
   \includegraphics[width=\linewidth]{Figures/Figure3g.pdf}     
   \caption{}
   \label{fig:FigureC7h} 
\end{subfigure}
\end{minipage}
}
\caption{The solution reconstruction graphs of the Stainless steel near-perfect emitter target for all optimization algorithms: (a) BO algorithm. (b) DE algorithm. (c) LBFGSB algorithm. (d) MADS algorithm. (e) NM algorithm. (f) PSO algorithm. (g) Random sampling. (h) ALPS algorithm.}
\label{fig:FigureC7}
\end{figure}

\begin{figure}[H]
\centering
\adjustbox{center}{
\begin{minipage}{1.5\textwidth}
\begin{subfigure}[b]{0.24\textwidth}
   \includegraphics[width=\linewidth]{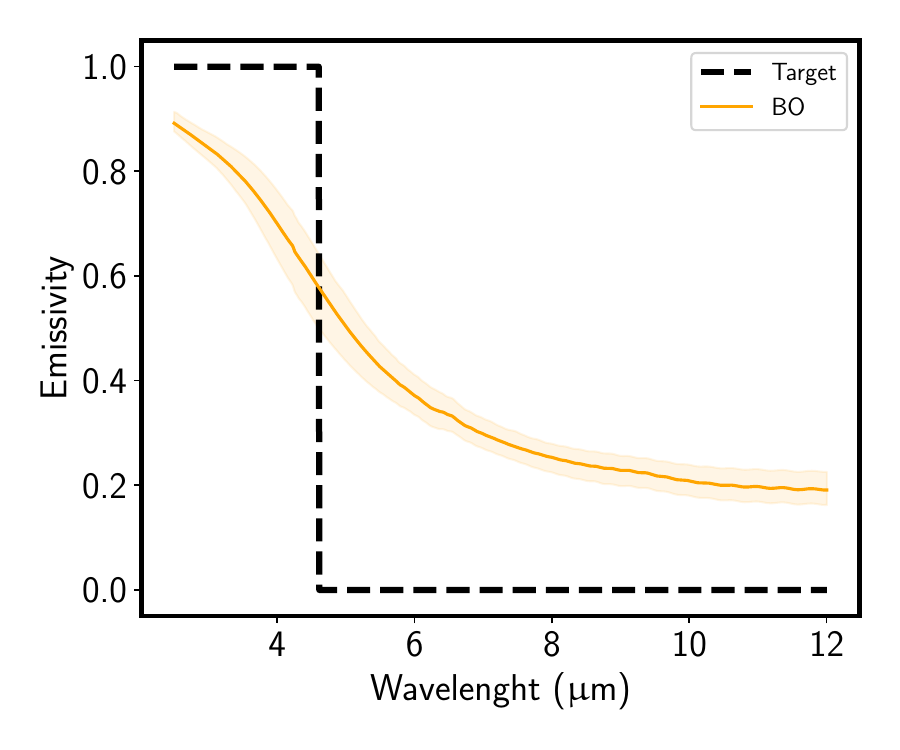}
   \caption{}
   \label{fig:FigureC8a}
\end{subfigure}
\hfill
\begin{subfigure}[b]{0.24\textwidth}
   \includegraphics[width=\linewidth]{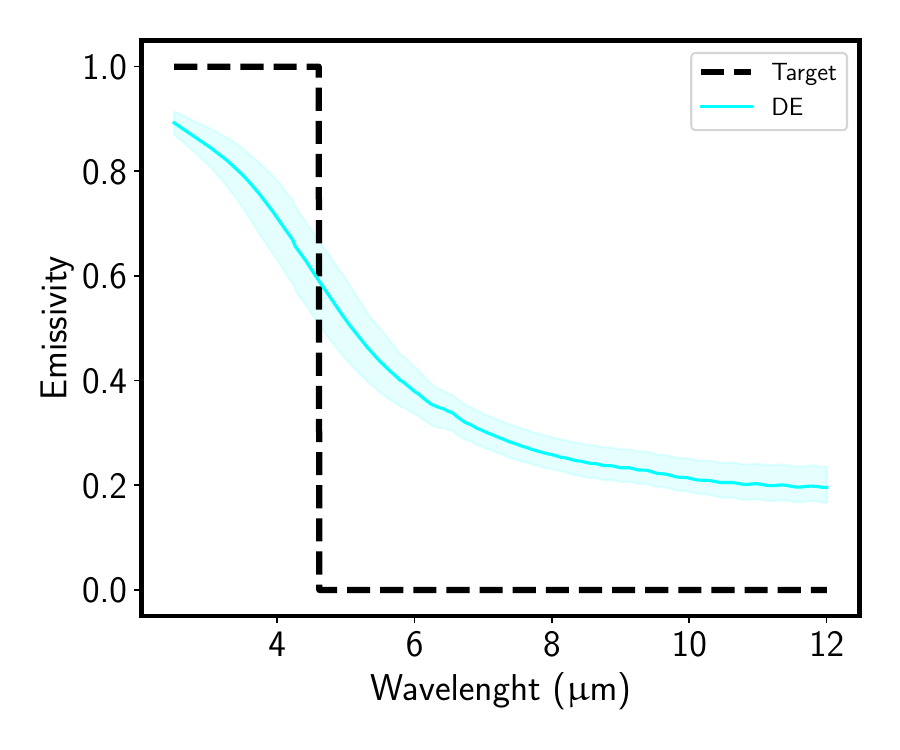}     
   \caption{}
   \label{fig:FigureC8b} 
\end{subfigure}
\hfill
\begin{subfigure}[b]{0.24\textwidth}
   \includegraphics[width=\linewidth]{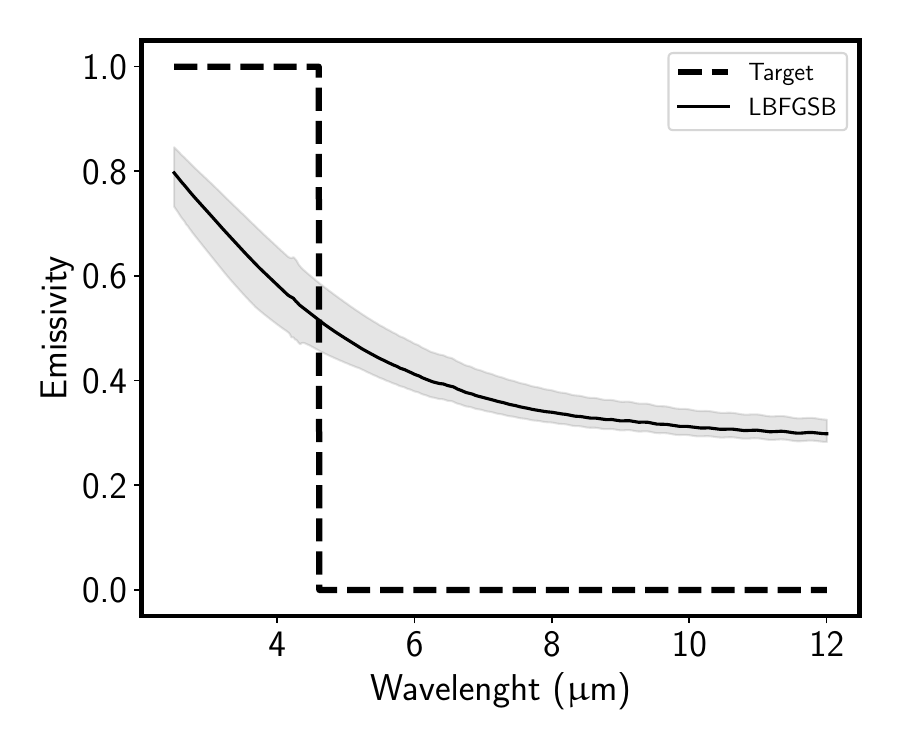}     
   \caption{}
   \label{fig:FigureC8c} 
\end{subfigure}
\hfill
\begin{subfigure}[b]{0.24\textwidth}
   \includegraphics[width=\linewidth]{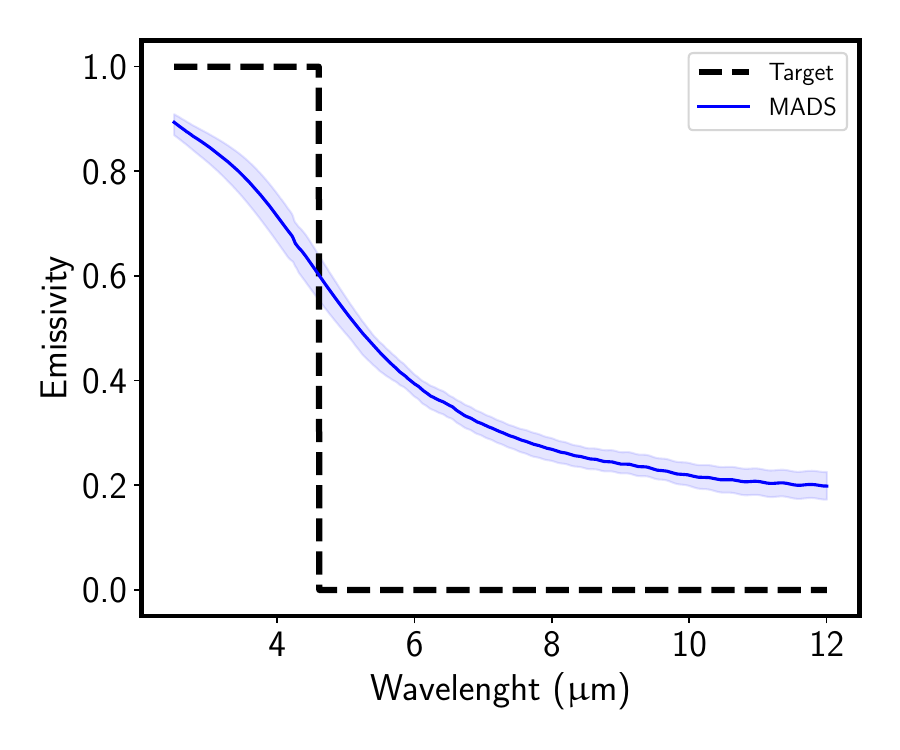}     
   \caption{}
   \label{fig:FigureC8d} 
\end{subfigure}
\end{minipage}
}

\adjustbox{center}{
\begin{minipage}{1.5\textwidth}
\begin{subfigure}[b]{0.24\textwidth}
   \includegraphics[width=\linewidth]{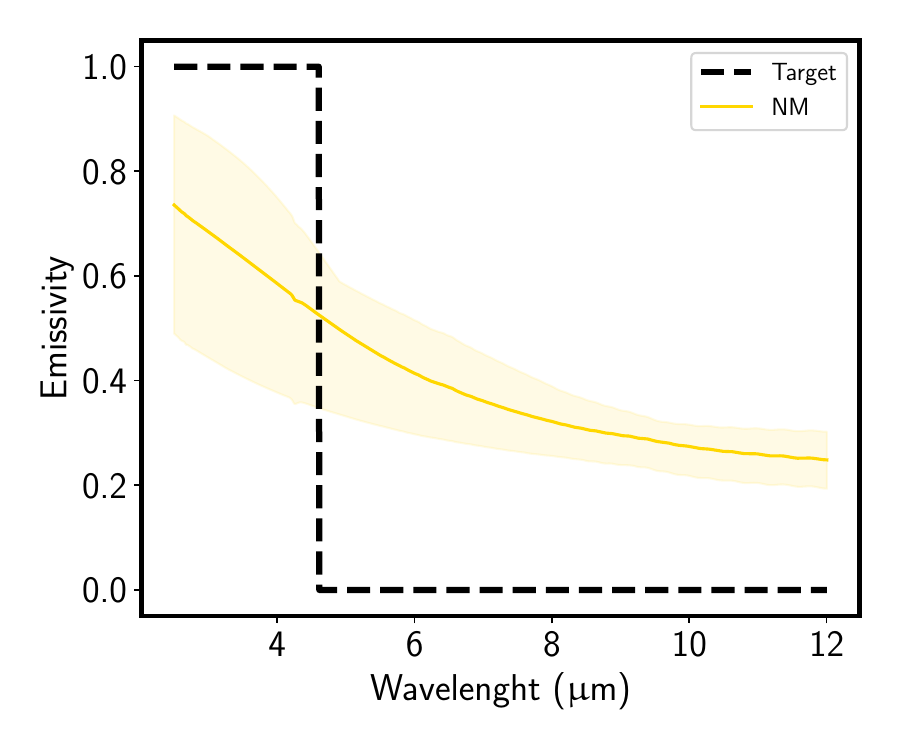}     
   \caption{}
   \label{fig:FigureC8e} 
\end{subfigure}
\hfill
\begin{subfigure}[b]{0.24\textwidth}
   \includegraphics[width=\linewidth]{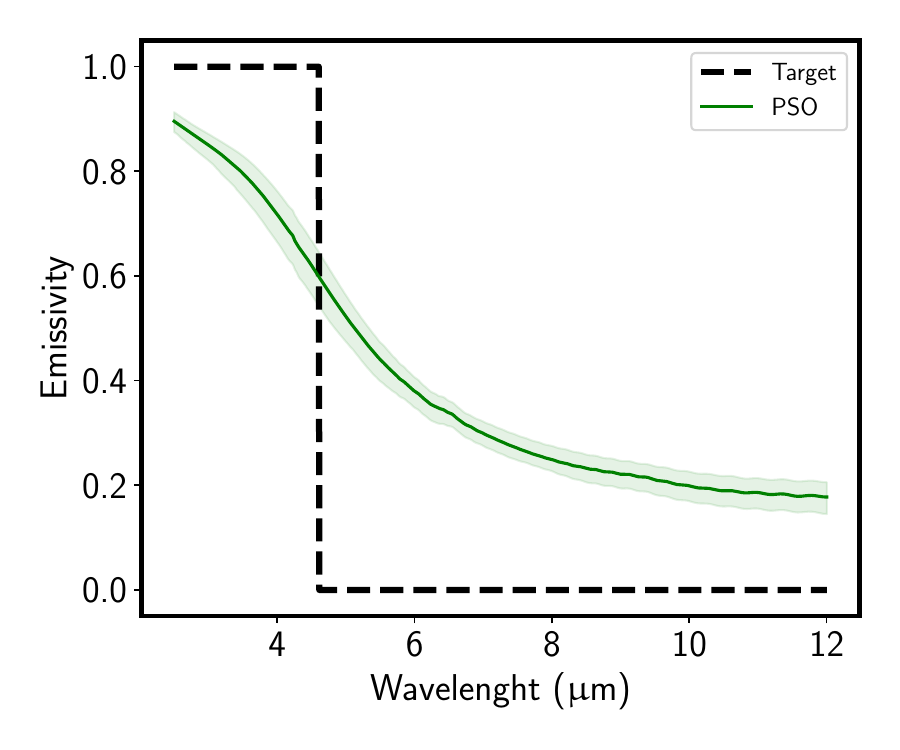}     
   \caption{}
   \label{fig:FigureC8f} 
\end{subfigure}
\hfill
\begin{subfigure}[b]{0.24\textwidth}
   \includegraphics[width=\linewidth]{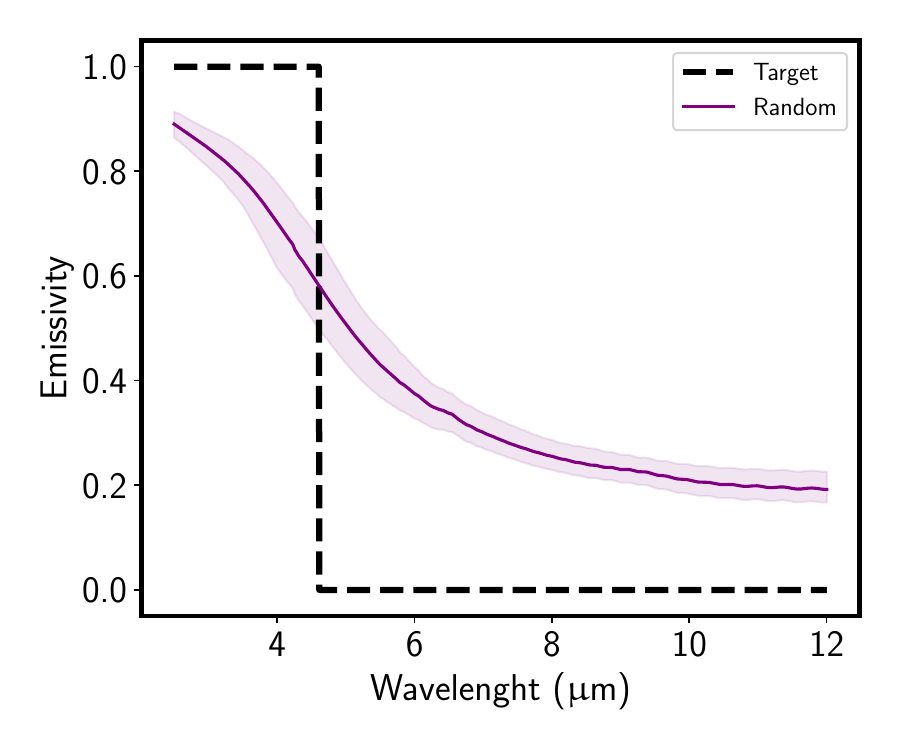}     
   \caption{}
   \label{fig:FigureC8g} 
\end{subfigure}
\hfill
\begin{subfigure}[b]{0.24\textwidth}
   \includegraphics[width=\linewidth]{Figures/Figure3h.pdf}     
   \caption{}
   \label{fig:FigureC8h} 
\end{subfigure}
\end{minipage}
}
\caption{The solution reconstruction graphs of the Stainless steel TPV emitter target for all optimization algorithms: (a) BO algorithm. (b) DE algorithm. (c) LBFGSB algorithm. (d) MADS algorithm. (e) NM algorithm. (f) PSO algorithm. (g) Random sampling. (h) ALPS algorithm.}
\label{fig:FigureC8}
\end{figure}

The detailed results of the photonic surface inverse design benchmarks are presented in this section. In Fig. \ref{fig:FigureC4}, the convergence graphs with included uncertainty for all benchmarks and algorithms are shown. Fig. \ref{fig:FigureC5}-\ref{fig:FigureC8} show the solution reconstruction graphs of all algorithms (except for ALPS) for the Inconel near-perfect emitter, Inconel TPV emitter, Stainless steel near-perfect emitter, and Stainless-steel TPV emitter, respectively. Finally, detailed convergence statistics are presented in Tab. \ref{tab:detailed_results_inconel_near_perfect}-\ref{tab:detailed_results_ss_TPV} following the same order as the solution reconstruction graphs.

\newpage

\begin{table}[!h]
\centering
\caption{Detailed convergence statistics of the $\epsilon$ value for the Inconel near-perfect emitter target, for all algorithms. The results are for 100 repeated runs. Smaller values are better.}
\begin{tabular}{ c | c | c | c | c}
 \textbf{Algorithm}  & \textbf{Mean $\epsilon$} &  \textbf{Std $\epsilon$} & \textbf{Maximum $\epsilon$} & \textbf{Minimum $\epsilon$}\\
\hline
\hline
ALPS & $\mathbf{0.020}$ & $\mathbf{0.003}$ & $\mathbf{0.020}$ & $\mathbf{0.020}$  \\
MADS & 0.050& 0.140 &0.620 & $\mathbf{0.020}$ \\
BO & 0.040 & 0.050  & 0.390 &  $\mathbf{0.020}$ \\
PSO & $\mathbf{0.020}$ &0.010&0.050& $\mathbf{0.020}$ \\
DE &0.150 & 0.120 & 0.500 &  $\mathbf{0.020}$\\
LBFGSB& 0.550&0.080 &0.690 &0.370  \\
Random &0.110&0.100 &   0.440 &$\mathbf{0.020}$  \\
NM & 0.570 &0.220 & 0.730& $\mathbf{0.020}$  \\
\end{tabular}

\label{tab:detailed_results_inconel_near_perfect}
\end{table}

\newpage
\begin{table}[!h]
\centering
\caption{Detailed convergence statistics of the $\epsilon$ value for the Inconel TPV emitter target, for all algorithms. The results are for 100 repeated runs. Smaller values are better.}
\begin{tabular}{ c | c | c | c | c}
 \textbf{Algorithm}  & \textbf{Mean $\epsilon$} &  \textbf{Std $\epsilon$} & \textbf{Maximum $\epsilon$} & \textbf{Minimum $\epsilon$}\\
\hline
\hline
ALPS & $\mathbf{0.300}$ & $\mathbf{0.003}$ & $\mathbf{0.320}$ & $\mathbf{0.300}$  \\
MADS & 0.320& 0.050 &0.510 & $\mathbf{0.300}$ \\
BO & 0.310 & 0.020  & 0.400 &  $\mathbf{0.300}$\\
PSO & 0.310 &0.010&0.340& $\mathbf{0.300}$\\
DE &0.320 & 0.010 & 0.350 &  $\mathbf{0.300}$\\
LBFGSB&0.390&0.060 &0.510 &0.320  \\
Random &0.320&0.010 & 0.340 &$\mathbf{0.300}$  \\
NM & 0.450 &0.080 & 0.530& $\mathbf{0.300}$  \\
\end{tabular}

\label{tab:detailed_results_inconel_TPV}
\end{table}

\begin{table}[!h]
\centering
\caption{Detailed convergence statistics of the $\epsilon$ value for the Stainless steel near-perfect emitter target, for all algorithms. The results are for 100 repeated runs. Smaller values are better.}
\begin{tabular}{ c | c | c | c | c}
 \textbf{Algorithm}  & \textbf{Mean $\epsilon$} &  \textbf{Std $\epsilon$} & \textbf{Maximum $\epsilon$} & \textbf{Minimum $\epsilon$}\\
\hline
\hline
ALPS & $\mathbf{0.021}$ & $\mathbf{0.003}$ & $\mathbf{0.042}$ & $\mathbf{0.017}$  \\
MADS & 0.057& 0.145 &0.632 &  $\mathbf{0.017}$ \\
BO & 0.024 & 0.006 & 0.053 &  0.018 \\
PSO & 0.022 & 0.008 &0.065&  $\mathbf{0.017}$ \\
DE &0.056 & 0.040 & 0.247 &  0.022\\
LBFGSB& 0.194 & 0.171 &0.530 &0.035  \\
Random &0.050 & 0.027 &  0.133 &0.019  \\
NM & 0.482 &0.220 & 0.650& 0.018  \\
\end{tabular}
\label{tab:detailed_results_ss_near_perfect}
\end{table}

\newpage
\begin{table}[!h]
\centering
\caption{Detailed convergence statistics of the $\epsilon$ value for the Stainless steel TPV emitter target, for all algorithms. The results are for 100 repeated runs. Smaller values are better.}
\begin{tabular}{ c | c | c | c | c}
 \textbf{Algorithm}  & \textbf{Mean $\epsilon$} &  \textbf{Std $\epsilon$} & \textbf{Maximum $\epsilon$} & \textbf{Minimum $\epsilon$}\\
\hline
\hline
ALPS & $\mathbf{0.290}$ & $\mathbf{0.010}$ & $\mathbf{0.300}$ & $\mathbf{0.280}$  \\
MADS & 0.300& 0.070 &0.320 & $\mathbf{0.280}$\\
BO & $\mathbf{0.290}$ & $\mathbf{0.010}$  & 0.310 &  $\mathbf{0.280}$ \\
PSO & $\mathbf{0.290}$ &$\mathbf{0.010}$ &0.320& $\mathbf{0.280}$\\
DE &$\mathbf{0.290}$ & $\mathbf{0.010}$ & 0.320 &  $\mathbf{0.280}$\\
LBFGSB& 0.360&0.020 &0.440 &0.350  \\
Random &$\mathbf{0.290}$&$\mathbf{0.010}$ &  0.330 &$\mathbf{0.280}$  \\
NM & 0.390 &0.070 & 0.600& $\mathbf{0.280}$  \\
\end{tabular}

\label{tab:detailed_results_ss_TPV}
\end{table}

\newpage
\subsection{ALPS Hyperparameter Analysis}\label{subapp:alps_hyperparameters}

In this section we present the analysis of how ALPS and the RF algorithm hyperparameters affect its performance. We separately analyze the RF hyperparameters and the ALPS hyperparamters on four different benchmarks (two synthetic benchmarks, and two photonic surface inverse design benchmarks, namely Inconel near-perfect emitter target, and the Stainless steel TPV emitter target). Each ALPS run is repeated  100 times for all benchmarks to account the stochasticity in the algorithm.

\subsubsection{Influence of RF Hyperparameters}

Firstly, in order to study the influence of the RF algorithm, we consider changing only two parameters, namely the number of estimators (or trees) and the maximum tree depth value. We further investigate these two hyperparameters because they are the only ones that differ from the default settings in scikit-learn 1.2.2 in our experimental model (see App. \ref{subapp:ml_exp_results}), which was optimized in our previous work (\citet{grbcic2024ensemble}). The ALPS hyperparameters used during the RF algorithm are the same as the ones presented in Sec. \ref{sec:results}, i.e., initial sample size and batch size are both set to 5, and the surrogate sample size is set to 600.  

\begin{figure}[!h]
\centering
\begin{subfigure}[b]{0.49\textwidth}
   \includegraphics[width=\linewidth]{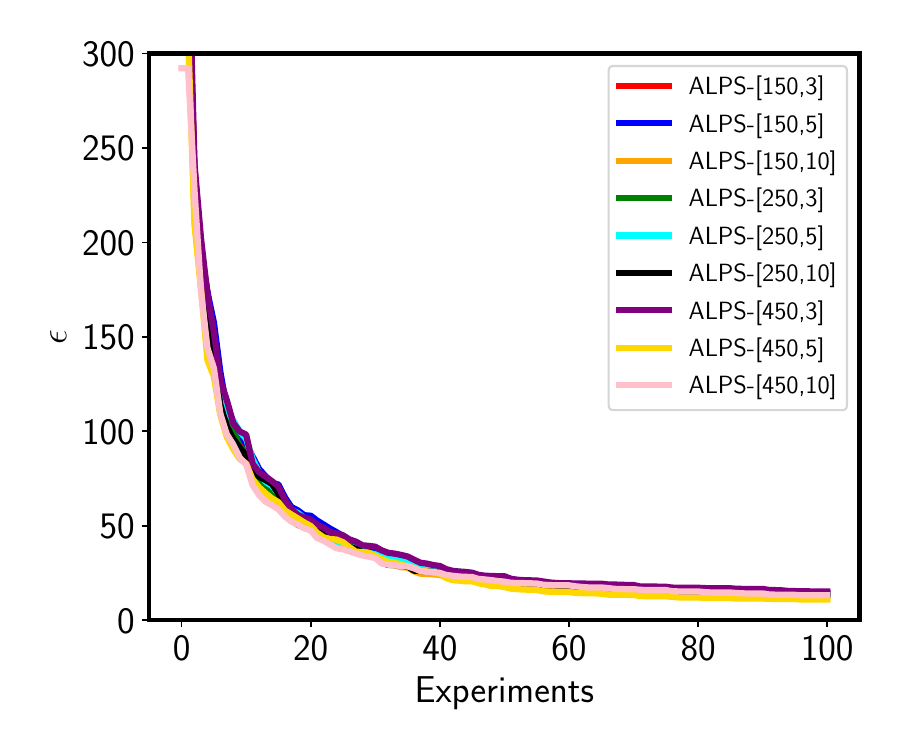}
   \caption{}
   \label{fig:FigureC9a}
\end{subfigure}
\begin{subfigure}[b]{0.49\textwidth}
   \includegraphics[width=\linewidth]{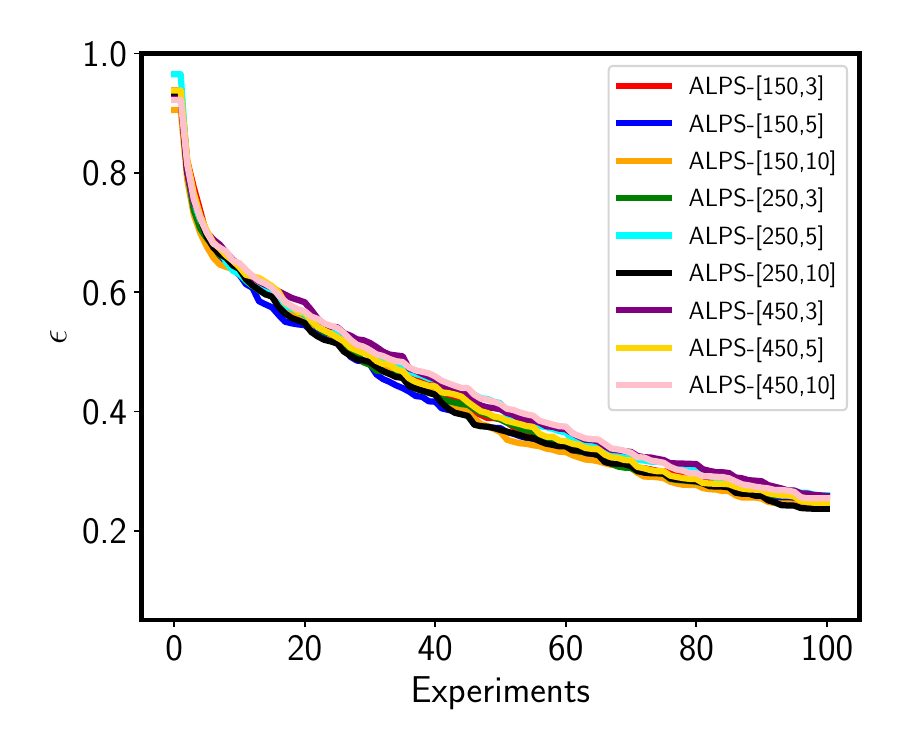}     
   \caption{}
   \label{fig:FigureC9b} 
\end{subfigure}
\begin{subfigure}[b]{0.49\textwidth}
   \includegraphics[width=\linewidth]{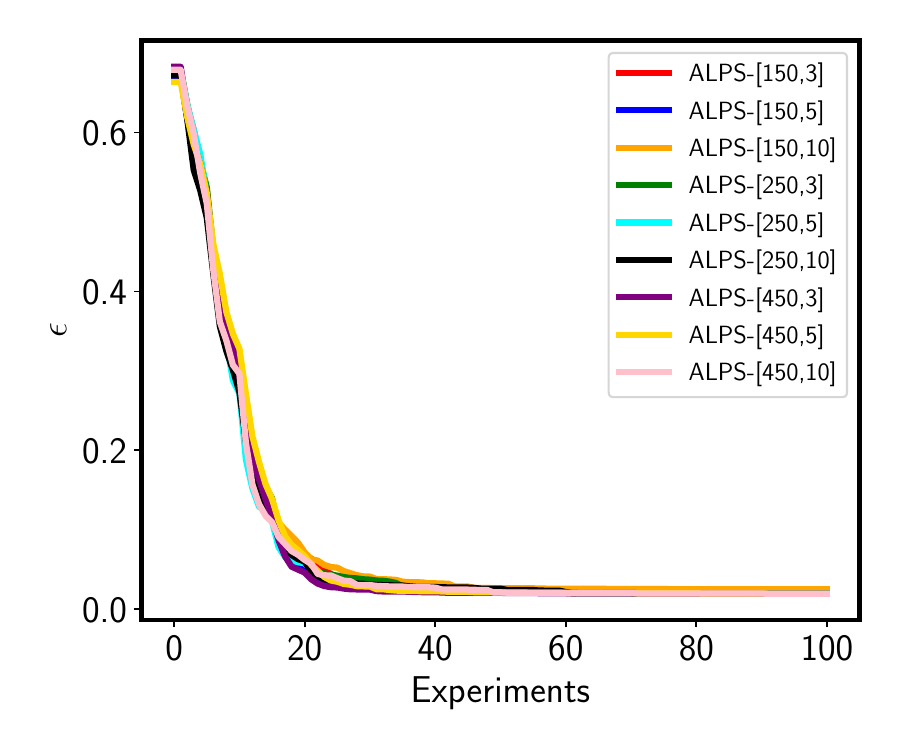}
   \caption{}
   \label{fig:FigureC9c}
\end{subfigure}
\begin{subfigure}[b]{0.49\textwidth}
   \includegraphics[width=\linewidth]{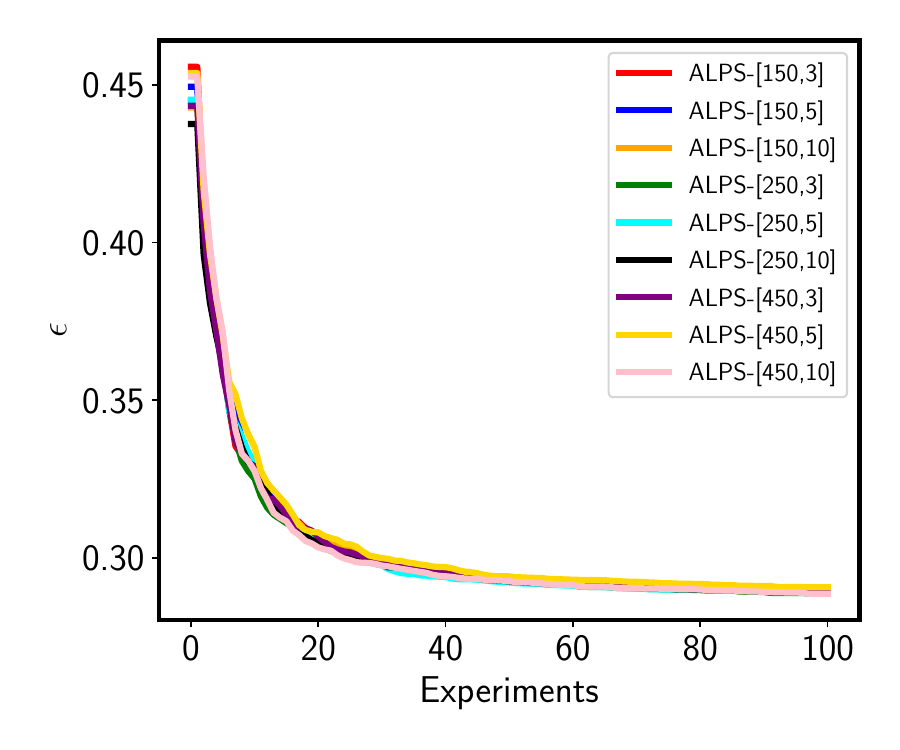}
   \caption{}
   \label{fig:FigureC9d}
\end{subfigure}
\caption[]{Convergence graphs for all investigated RF algorithm hyperparameters and a selection of benchmarks, the first number in each bracket represents the number of estimators (trees), while the second number is the maximum depth value: (a) Logistic growth benchmark. (b) Sinusoidal oscillation with damping benchmark. (c) Photonic surface benchmark: Inconel near-perfect emitter target. (d) Photonic surface benchmark: Stainless steel TPV emitter target.}
\label{fig:FigureC9}
\end{figure}

In Fig. \ref{fig:FigureC9}, the results of the RF algorithm used in ALPS are presented for the four different benchmarks. The number of estimators considered for this analysis is set to 150, 250, and 450, while the varied maximum depth values are 3, 5, and 10. It can be observed that changing the hyperparameters of the RF algorithm within ALPS does not affect the performance of ALPS with the exception the logistic growth benchmark case (Tab. \ref{tab:rf_hyperparameters_logistic_growth}). This means that using the default out-of-the-box parameters of the RF algorithm (in scikit-learn v1.2.2.) are suitable when we want to apply it within ALPS. This is beneficial as it shows the robustness of the approach, however, it also implies that there is not a lot of room for improvement of the approach that can be achieved by tuning the RF hyperparameters. Detailed convergence statistics for all benchmarks and all RF hyperparameter combinations are given in Tab. \ref{tab:rf_hyperparameters_logistic_growth}-\ref{tab:rf_hyperparameters_ss_TPV}.

\newpage

\newpage
\begin{table}[!h]
\centering
\caption{Detailed convergence statistics of the $\epsilon$ value for the logistic growth benchmark, for all investigated RF hyperparameter values.  Results are for 100 repeated runs.}
\begin{tabular}{ c | c | c | c }
 \textbf{Estimators}  &  \textbf{Max depth}  & \textbf{Mean $\epsilon$} &  \textbf{Std $\epsilon$} \\
\hline
\hline
150 & 3 & 14.37 & 11.97 \\
150 & 5 & 14.18 & 17.40 \\
150 & 10 & 11.45 & $\mathbf{6.66}$ \\
250 & 3 & 13.40 & 11.37 \\
250 & 5 & 13.95 & 13.12 \\
250 & 10 & 13.29 & 11.70 \\
450 & 3 & 15.14 & 13.40 \\
450 & 5 & $\mathbf{10.88}$ & 7.16 \\
450 & 10 & 13.13 & 8.29 \\
\end{tabular}

\label{tab:rf_hyperparameters_logistic_growth}
\end{table}

\begin{table}[!h]
\centering
\caption{Detailed convergence statistics of the $\epsilon$ value for the sinusoidal oscillation with damping benchmark, for all investigated RF hyperparameter values.  Results are for 100 repeated runs.}
\begin{tabular}{ c | c | c | c }
 \textbf{Estimators}  &  \textbf{Max depth}  & \textbf{Mean $\epsilon$} &  \textbf{Std $\epsilon$} \\
\hline
\hline
150 & 3 & 0.25 & 0.12 \\
150 & 5 & 0.25 & 0.12 \\
150 & 10 & $\mathbf{0.24}$& $\mathbf{0.10}$ \\
250 & 3 & 0.25 & 0.12 \\
250 & 5 & 0.26 & 0.12 \\
250 & 10 & $\mathbf{0.24}$ & $\mathbf{0.10}$ \\
450 & 3 & 0.26 & 0.12 \\
450 & 5 & 0.25 & 0.11 \\
450 & 10 & 0.25 & 0.11 \\
\end{tabular}

\label{tab:rf_hyperparameters_oscillation}
\end{table}

\newpage
\begin{table}[!h]
\centering
\caption{Detailed convergence statistics of the $\epsilon$ value for the photonic surface benchmark--Inconel near-perfect emitter target, for all investigated RF hyperparameter values. Results are for 100 repeated runs.}
\begin{tabular}{ c | c | c | c }
 \textbf{Estimators}  &  \textbf{Max depth}  & \textbf{Mean $\epsilon$} &  \textbf{Std $\epsilon$} \\
\hline
\hline
150 & 3 & $\mathbf{0.020}$ & $\mathbf{0.020}$ \\
150 & 5 & $\mathbf{0.020}$ & $\mathbf{0.020}$ \\
150 & 10 & 0.025 & 0.060 \\
250 & 3 & $\mathbf{0.020}$ & $\mathbf{0.020}$ \\
250 & 5 & $\mathbf{0.020}$ & $\mathbf{0.020}$ \\
250 & 10 & $\mathbf{0.020}$ & $\mathbf{0.020}$ \\
450 & 3 & $\mathbf{0.020}$ & $\mathbf{0.020}$ \\
450 & 5 & $\mathbf{0.020}$ & $\mathbf{0.020}$ \\
450 & 10 & $\mathbf{0.020}$ & $\mathbf{0.020}$ \\
\end{tabular}

\label{tab:rf_hyperparameters_inconel_near_perfect}
\end{table}

\begin{table}[!h]
\centering
\caption{Detailed convergence statistics of the $\epsilon$ value for the photonic surface benchmark--Stainless steel TPV emitter target, for all investigated RF hyperparameter values.  Results are for 100 repeated runs.}
\begin{tabular}{ c | c | c | c }
 \textbf{Estimators}  &  \textbf{Max depth}  & \textbf{Mean $\epsilon$} &  \textbf{Std $\epsilon$} \\
\hline
\hline
150 & 3 & 0.290 & 0.010 \\
150 & 5 & 0.290 & 0.010 \\
150 & 10 & 0.290 & 0.010 \\
250 & 3 & 0.290 & 0.010 \\
250 & 5 & 0.290 & 0.010 \\
250 & 10 & 0.290 & 0.010 \\
450 & 3 & 0.290 & 0.010 \\
450 & 5 & 0.290 & 0.010 \\
450 & 10 & 0.290 & 0.010 \\
\end{tabular}

\label{tab:rf_hyperparameters_ss_TPV}
\end{table}

\newpage
\subsubsection{Influence of ALPS Hyperparameters}

\begin{figure}[!h]
\centering
\begin{subfigure}[b]{0.49\textwidth}
   \includegraphics[width=\linewidth]{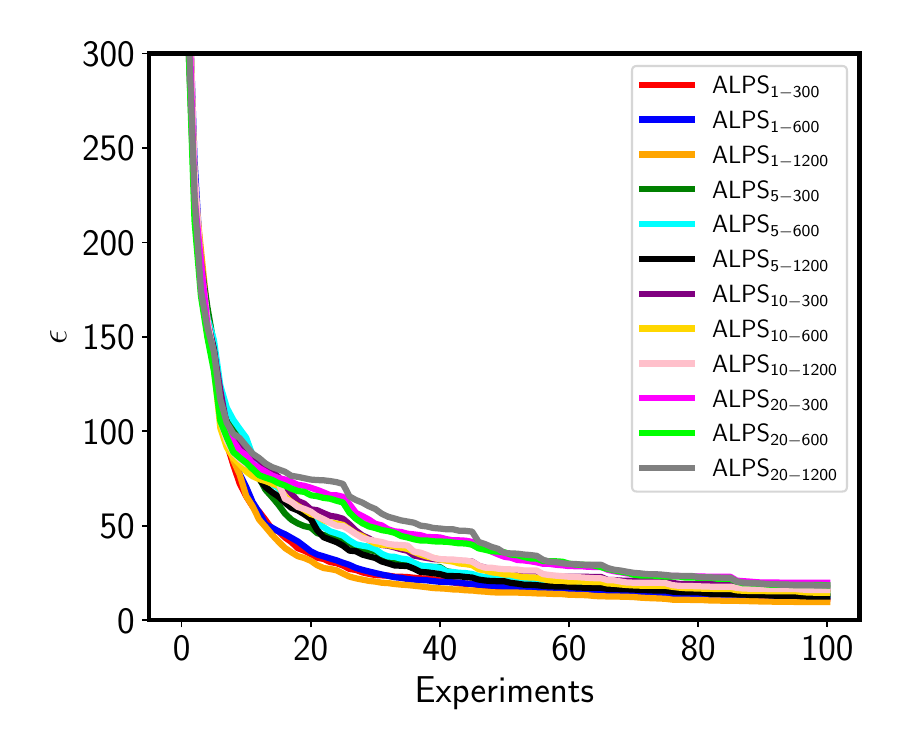}
   \caption{}
   \label{fig:FigureC10a}
\end{subfigure}
\begin{subfigure}[b]{0.49\textwidth}
   \includegraphics[width=\linewidth]{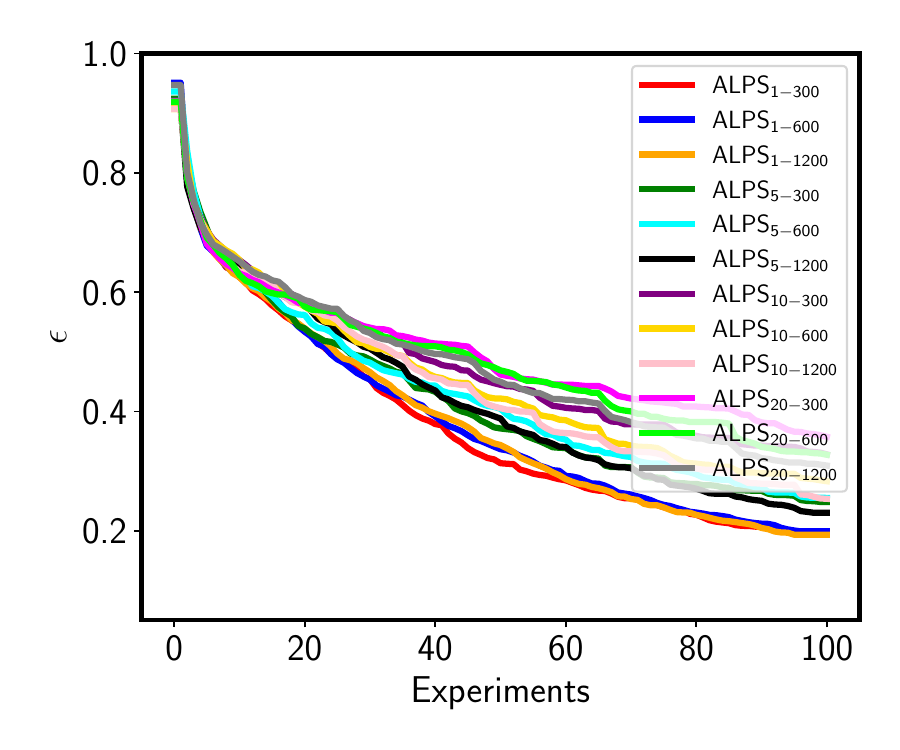}     
   \caption{}
   \label{fig:FigureC10b} 
\end{subfigure}
\begin{subfigure}[b]{0.49\textwidth}
   \includegraphics[width=\linewidth]{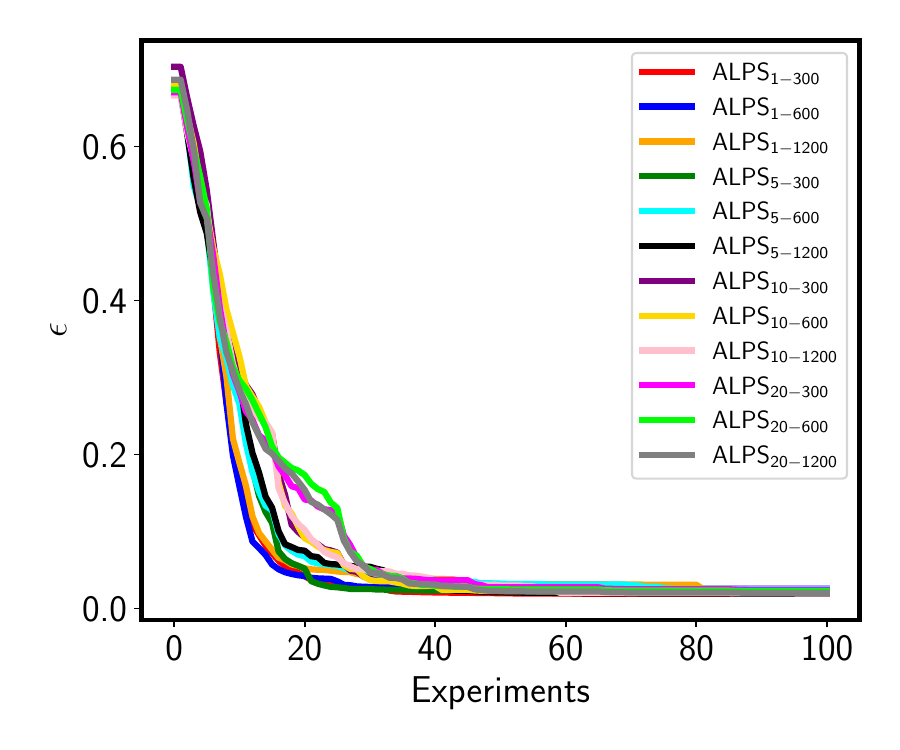}
   \caption{}
   \label{fig:FigureC10c}
\end{subfigure}
\begin{subfigure}[b]{0.49\textwidth}
   \includegraphics[width=\linewidth]{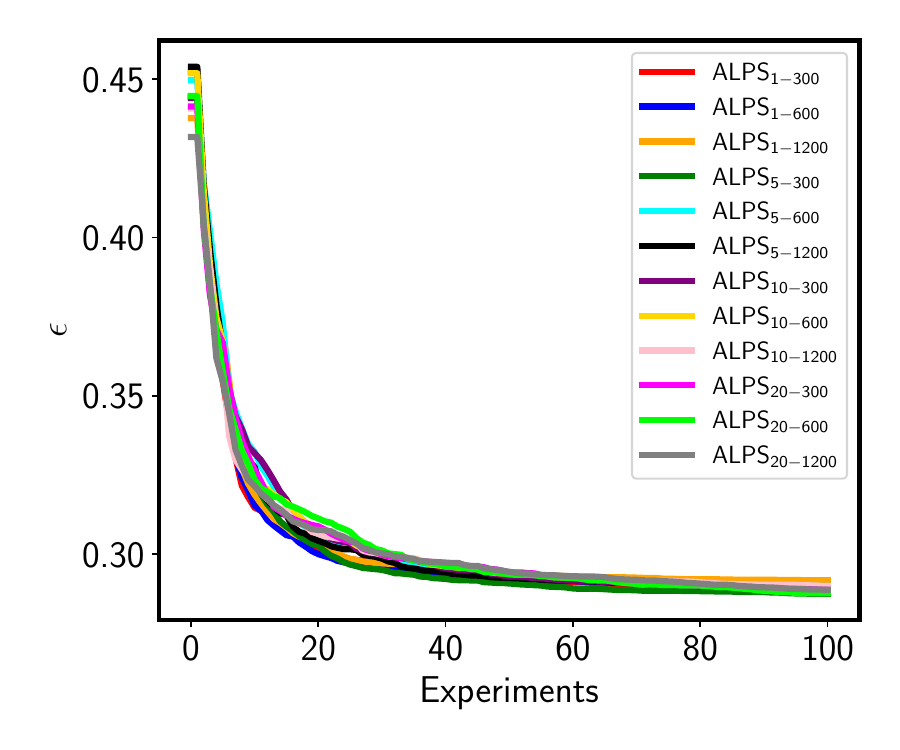}
   \caption{}
   \label{fig:FigureC10d}
\end{subfigure} 
\caption[]{Convergence graphs for all investigated ALPS hyperparameters and a selection of benchmarks, the first number in the index represents the batch size, while the second number is the surrogate sample size: (a) Logistic growth benchmark. (b) Sinusoidal oscillation with damping benchmark. (c) Photonic surface benchmark: Inconel near--perfect emitter target. (d) Photonic surface benchmark: Stainless steel TPV emitter target.}
\label{fig:FigureC10}
\end{figure}

In this section we show the results of how ALPS hyperparameters influence its performance. For our initial comparison with other optimization algorithms, we use the initial sample size, and batch size both set to 5, and the surrogate sample size set to 600. We keep the same initial size, however, we vary the batch size as 1, 5, 10, and 20, and we set the surrogate sample size as either 300, 600 or 1200. We repeat every run  100 times and obtain the convergence graphs for each ALPS hyperparameter variant. Fig. \ref{fig:FigureC10} shows the convergence graphs of all ALPS investigated hyperparameters, for all four benchmarks. ALPS hyperparameters have a larger impact on the synthetic benchmarks (Fig. \ref{fig:FigureC10a}-\ref{fig:FigureC10b}) than on the photonic surface inverse design benchmarks (Fig. \ref{fig:FigureC10c}-\ref{fig:FigureC10d}). Detailed convergence statistics for all benchmarks are given in Tab. \ref{tab:ALPS_hyperparameters_logistic_growth}-\ref{tab:ALPS_hyperparameters_ss_tpv}.

Increasing the surrogate sample size has an impact on the logistic growth benchmark (Tab \ref{tab:ALPS_hyperparameters_logistic_growth}). For both synthetic benchmarks, using smaller batch sizes (sequential sampling) enhances the performance of ALPS in terms of the mean $\epsilon$, even more so on the sinusoidal oscillation with damping benchmark (Tab \ref{tab:ALPS_hyperparameters_sinusoidal}). Utilizing the worst performing hyperparameter set ($n_{batch}$=20, $n_s$=300) on the logistic growth benchmark yields better results in the mean $\epsilon$ (19.79), than PSO, which is the second best performing optimization algorithm, that yields the mean  of$\epsilon$ of 26.20 (Tab. \ref{tab:detailed_results_logistic}). For the sinusoidal oscillation with damping benchmark, the worst performing hyperparameter set is the same, and it would rank as the third best if compared to the other investigated optimization algorithms (Tab. \ref{tab:detailed_results_sinusoidal}).

For the main problem of photonic surface inverse design, the influence of ALPS hyperparameters is not as pronounced. In Fig. \ref{fig:FigureC10c}, the convergence graph for the Inconel near-perfect emitter target is shown. In terms of the final mean value of $\epsilon$, there is not a large difference, with the exception of the case when $n_{batch}$ is 1, and the $n_s$ values are either 600 or 1200, the final mean of $\epsilon$ is higher compared to other hyperparameter sets (Tab. \ref{tab:ALPS_hyperparameters_inconel_near_perfect}). Both of these hyperparameter sets had exactly one run out of 100, where ALPS converged to a local optimum and obtained a maximum final $\epsilon$ value of 0.68 and 0.61 for $n_s$=600, and $n_s$=1200, respectively. For the Stainless steel TPV emitter target, all ALPS hyperparameter variants perform the same as presented in Tab. \ref{tab:ALPS_hyperparameters_ss_tpv}). The main benefit of using a smaller value of $n_{batch}$ can be observed in Fig. \ref{fig:FigureC10c} for the Inconel near-perfect emitter, where lower values have a tendency of converging to a satisfying design shortly after 20 experimental model evaluations, while a larger value delays this for additional experimental model evaluations. This performance is not apparent for the Stainless steel TPV emitter (Fig. \ref{fig:FigureC10d}).
\newpage

\newpage
\begin{table}[!h]
\centering
\caption{Detailed convergence statistics of the $\epsilon$ value for the logistic growth benchmark, for all investigated ALPS hyperparameter values.  Results are for 100 repeated runs.}
\begin{tabular}{ c | c | c | c }
 $\mathbf{n_{batch}}$  &   $\mathbf{n_s}$   & \textbf{Mean $\epsilon$} &  \textbf{Std $\epsilon$} \\
\hline
\hline
1 & 300 & 11.50 & 12.15 \\
1 & 600 & 12.58 & 13.40 \\
1 & 1200 & $\mathbf{9.67}$ & 10.51 \\
5 & 300 & 14.10 & $\mathbf{7.42}$ \\
5 & 600 & 13.43 & 9.68 \\
5 & 1200 & 12.56 & 12.25 \\
10 & 300 & 16.75 & 9.64 \\
10 & 600 & 14.57 & 10.33 \\
10 & 1200 & 16.02 & 13.62 \\
20 & 300 & 19.79 & 8.21 \\
20 & 600 & 18.40 & 10.44 \\
20 & 1200 & 18.33 & 11.96 \\
\end{tabular}

\label{tab:ALPS_hyperparameters_logistic_growth}
\end{table}

\newpage
\begin{table}[!h]
\centering
\caption{Detailed convergence statistics of the $\epsilon$ value for the sinusoidal oscillation with damping benchmark, for all investigated ALPS hyperparameter values.  Results are for 100 repeated runs.}
\begin{tabular}{ c | c | c | c }
 $\mathbf{n_{batch}}$  &   $\mathbf{n_s}$   & \textbf{Mean $\epsilon$} &  \textbf{Std $\epsilon$} \\
\hline
\hline
1 & 300 & $\mathbf{0.19}$ & $\mathbf{0.10}$ \\
1 & 600 & 0.20 & 0.12 \\
1 & 1200 &$\mathbf{0.19}$ & 0.14 \\
5 & 300 & 0.25 & 0.10 \\
5 & 600 & 0.25 & 0.12 \\
5 & 1200 & 0.23 & 0.11 \\
10 & 300 & 0.33 & 0.12 \\
10 & 600 & 0.28 & 0.12 \\
10 & 1200 & 0.25 & 0.13 \\
20 & 300 & 0.36 & 0.12 \\
20 & 600 & 0.33 & 0.12 \\
20 & 1200 & 0.31 & 0.14 \\
\end{tabular}

\label{tab:ALPS_hyperparameters_sinusoidal}
\end{table}

\newpage
\begin{table}[!h]
\centering
\caption{Detailed convergence statistics of the $\epsilon$ value for the photonic surface benchmark--Inconel near-perfect emitter target, for all investigated ALPS hyperparameter values.  Results are for 100 repeated runs.}
\begin{tabular}{ c | c | c | c }
 $\mathbf{n_{batch}}$  &   $\mathbf{n_s}$   & \textbf{Mean $\epsilon$} &  \textbf{Std $\epsilon$} \\
\hline
\hline
1 & 300 & $\mathbf{0.020}$ &  $\mathbf{0.001}$ \\
1 & 600 & 0.030 & 0.070 \\
1 & 1200 & 0.030 & 0.060 \\
5 & 300 &  $\mathbf{0.020}$ & $\mathbf{0.001}$  \\
5 & 600 &  $\mathbf{0.020}$ & $\mathbf{0.001}$  \\
5 & 1200 &  $\mathbf{0.020}$ & $\mathbf{0.001}$  \\
10 & 300 &  $\mathbf{0.020}$ & $\mathbf{0.001}$  \\
10 & 600 &  $\mathbf{0.020}$ & $\mathbf{0.001}$  \\
10 & 1200 &  $\mathbf{0.020}$ & $\mathbf{0.001}$  \\
20 & 300 &  $\mathbf{0.020}$ & $\mathbf{0.001}$  \\
20 & 600 &  $\mathbf{0.020}$ & $\mathbf{0.001}$  \\
20 & 1200 &  $\mathbf{0.020}$ & $\mathbf{0.001}$  \\
\end{tabular}

\label{tab:ALPS_hyperparameters_inconel_near_perfect}
\end{table}

\newpage
\begin{table}[!h]
\centering
\caption{Detailed convergence statistics of the $\epsilon$ value for the photonic surface benchmark--Stainless steel TPV emitter target, for all investigated ALPS hyperparameter values.  Results are for 100 repeated runs.}
\begin{tabular}{ c | c | c | c }
 $\mathbf{n_{batch}}$  &   $\mathbf{n_s}$   & \textbf{Mean $\epsilon$} &  \textbf{Std $\epsilon$} \\
\hline
\hline
1 & 300 & $\mathbf{0.290}$ & $\mathbf{0.010}$ \\
1 & 600 & $\mathbf{0.290}$ & $\mathbf{0.010}$ \\
1 & 1200 & $\mathbf{0.290}$ & $\mathbf{0.010}$ \\
5 & 300 & $\mathbf{0.290}$ & $\mathbf{0.010}$ \\
5 & 600 & $\mathbf{0.290}$ & $\mathbf{0.010}$ \\
5 & 1200 & $\mathbf{0.290}$ & $\mathbf{0.010}$ \\
10 & 300 & $\mathbf{0.290}$ & $\mathbf{0.010}$ \\
10 & 600 & $\mathbf{0.290}$ & $\mathbf{0.010}$ \\
10 & 1200 & $\mathbf{0.290}$ & $\mathbf{0.010}$\\
20 & 300 & $\mathbf{0.290}$ & $\mathbf{0.010}$ \\
20 & 600 & $\mathbf{0.290}$ & $\mathbf{0.010}$ \\
20 & 1200 & $\mathbf{0.290}$ & $\mathbf{0.010}$ \\
\end{tabular}

\label{tab:ALPS_hyperparameters_ss_tpv}
\end{table}

\subsection{ALPS with Warm Starting Inverse Design Convergence}\label{subapp:alpswarmstarting}

In this section the ALPS convergence results with and without warm starting   are presented. The convergence details for the cross-target warm starting cases are presented in Tab. \ref{tab:detailed_results_cross_target_inconel_100}-\ref{tab:detailed_results_cross_target_ss_step}, while Tab. \ref{tab:detailed_results_cross_material_inconel_100}-\ref{tab:detailed_results_cross_material_ss_step} show the results of the cross-material warm starting cases.

\begin{table}[!h]
\centering
\caption{Detailed convergence statistics of the $\epsilon$ value for ALPS with and without cross-target warm starting. The target is the Inconel near-perfect emitter warm started by a prior inverse design model for an Inconel TPV emitter. The results are for 100 repeated runs. }
\begin{tabular}{ c | c | c | c | c}
 \textbf{Algorithm}  & \textbf{Mean $\epsilon$} &  \textbf{Std $\epsilon$} & \textbf{Maximum $\epsilon$} & \textbf{Minimum $\epsilon$}\\
\hline
\hline
ALPS$_{ws=True}$ & $\mathbf{0.020}$ &  $\mathbf{0.001}$  &  $\mathbf{0.020}$  &  $\mathbf{0.020}$  \\
ALPS$_{ws=False}$ & $\mathbf{0.020}$  & $\mathbf{0.001}$ &  $\mathbf{0.020}$  &  $\mathbf{0.020}$ \\
\end{tabular}

\label{tab:detailed_results_cross_target_inconel_100}
\end{table}

\begin{table}[!h]
\centering
\caption{Detailed convergence statistics of the $\epsilon$ value for ALPS with and without cross-target warm starting. The target is the Inconel TPV emitter warm started by a prior inverse design model for an Inconel near-perfect emitter. The results are for 100 repeated runs.}
\begin{tabular}{ c | c | c | c | c}
 \textbf{Algorithm}  & \textbf{Mean $\epsilon$} &  \textbf{Std $\epsilon$} & \textbf{Maximum $\epsilon$} & \textbf{Minimum $\epsilon$}\\
\hline
\hline
ALPS$_{ws=True}$ & $\mathbf{0.300}$ &$\mathbf{0.001}$ & $\mathbf{0.310}$ & $\mathbf{0.300}$ \\
ALPS$_{ws=False}$ & $\mathbf{0.300}$ & $\mathbf{0.001}$ & $\mathbf{0.310}$ & $\mathbf{0.300}$ \\
\end{tabular}

\label{tab:detailed_results_cross_target_inconel_step}
\end{table}

\newpage
\begin{table}[!h]
\centering
\caption{Detailed convergence statistics of the $\epsilon$ value for ALPS with and without cross-target warm starting. The target is the Stainless steel near-perfect emitter warm started by a prior inverse design model for a Stainless steel TPV emitter. The results are for 100 repeated runs. }
\begin{tabular}{ c | c | c | c | c}
 \textbf{Algorithm}  & \textbf{Mean $\epsilon$} &  \textbf{Std $\epsilon$} & \textbf{Maximum $\epsilon$} & \textbf{Minimum $\epsilon$}\\
\hline
\hline
ALPS$_{ws=True}$ & $\mathbf{0.020}$ & $\mathbf{0.001}$ & $\mathbf{0.040}$& $\mathbf{0.020}$\\
ALPS$_{ws=False}$ & $\mathbf{0.020}$ & $\mathbf{0.001}$ & $\mathbf{0.040}$ & $mathbf{0.020}$ \\
\end{tabular}

\label{tab:detailed_results_cross_target_ss_100}
\end{table}

\begin{table}[!h]
\centering
\caption{Detailed convergence statistics of the $\epsilon$ value for ALPS with and without cross-target warm starting. The target is the Stainless steel TPV emitter warm started by a prior inverse design model for a Stainless steel near-perfect emitter. The results are for 100 repeated runs.}
\begin{tabular}{ c | c | c | c | c}
 \textbf{Algorithm}  & \textbf{Mean $\epsilon$} &  \textbf{Std $\epsilon$} & \textbf{Maximum $\epsilon$} & \textbf{Minimum $\epsilon$}\\
\hline
\hline
ALPS$_{ws=True}$ & $\mathbf{0.290}$&$\mathbf{0.010}$& $\mathbf{0.300}$ & $\mathbf{0.280}$ \\
ALPS$_{ws=False}$ & $\mathbf{0.290}$ & $\mathbf{0.010}$ & 0.310 & $\mathbf{0.280}$\\
\end{tabular}

\label{tab:detailed_results_cross_target_ss_step}
\end{table}

\begin{table}[!h]
\centering
\caption{Detailed convergence statistics of the $\epsilon$ value for ALPS with and without cross-material warm starting. The target is the Inconel near-perfect emitter warm started by a prior inverse design model for an Stainless steel TPV emitter. The results are for 100 repeated runs. }
\begin{tabular}{ c | c | c | c | c}
 \textbf{Algorithm}  & \textbf{Mean $\epsilon$} &  \textbf{Std $\epsilon$} & \textbf{Maximum $\epsilon$} & \textbf{Minimum $\epsilon$}\\
\hline
\hline
ALPS$_{ws=True}$ &  $\mathbf{0.020}$ & $\mathbf{0.001}$ & $\mathbf{0.020}$ & $\mathbf{0.020}$ \\
ALPS$_{ws=False}$ & $\mathbf{0.020}$ &$\mathbf{0.001}$ & $\mathbf{0.020}$ & $\mathbf{0.020}$ \\
\end{tabular}

\label{tab:detailed_results_cross_material_inconel_100}
\end{table}

\begin{table}[!h]
\centering
\caption{Detailed convergence statistics of the $\epsilon$ value for ALPS with and without cross-material warm starting. The target is the Inconel TPV emitter warm started by a prior inverse design model for an Stainless steel near-perfect emitter. The results are for 100 repeated runs. }
\begin{tabular}{ c | c | c | c | c}
 \textbf{Algorithm}  & \textbf{Mean $\epsilon$} &  \textbf{Std $\epsilon$} & \textbf{Maximum $\epsilon$} & \textbf{Minimum $\epsilon$}\\
\hline
\hline
ALPS$_{ws=True}$ & 0.300 & $\mathbf{0.010}$ & 0.310 & $\mathbf{0.280}$\\
ALPS$_{ws=False}$ & $\mathbf{0.290}$ & $\mathbf{0.010}$ & $\mathbf{0.300}$& $\mathbf{0.280}$ \\
\end{tabular}

\label{tab:detailed_results_cross_material_inconel_step}
\end{table}

\newpage
\begin{table}[!h]
\centering
\caption{Detailed convergence statistics of the $\epsilon$ value for ALPS with and without cross-material warm starting. The target is the Stainless steel near-perfect emitter warm started by a prior inverse design model for a Inconel TPV emitter. The results are for 100 repeated runs. }
\begin{tabular}{ c | c | c | c | c}
 \textbf{Algorithm}  & \textbf{Mean $\epsilon$} &  \textbf{Std $\epsilon$} & \textbf{Maximum $\epsilon$} & \textbf{Minimum $\epsilon$}\\
\hline
\hline
ALPS$_{ws=True}$ & $\mathbf{0.020}$ & $\mathbf{0.001}$ & $\mathbf{0.040}$ & $\mathbf{0.020}$ \\
ALPS$_{ws=False}$ & $\mathbf{0.020}$ & 0$\mathbf{0.001}$ &$\mathbf{0.040}$ & $\mathbf{0.020}$\\
\end{tabular}

\label{tab:detailed_results_cross_material_ss_100}
\end{table}

\begin{table}[!h]
\centering
\caption{Detailed convergence statistics of the $\epsilon$ value for ALPS with and without cross-material warm starting. The target is the Stainless steel TPV emitter warm started by a prior inverse design model for a Inconel near-perfect emitter. The results are for 100 repeated runs. }
\begin{tabular}{ c | c | c | c | c}
 \textbf{Algorithm}  & \textbf{Mean $\epsilon$} &  \textbf{Std $\epsilon$} & \textbf{Maximum $\epsilon$} & \textbf{Minimum $\epsilon$}\\
\hline
\hline
ALPS$_{ws=True}$ & $\mathbf{0.300}$ & $\mathbf{0.001}$&$\mathbf{0.310}$ & $\mathbf{0.300}$ \\
ALPS$_{ws=False}$ & $\mathbf{0.300}$ & $\mathbf{0.001}$& 0.320 & $\mathbf{0.300}$ \\
\end{tabular}

\label{tab:detailed_results_cross_material_ss_step}
\end{table}

\end{appendices}

\newpage
\bibliography{sn-bibliography}

\begin{thebibliography}{10}
\providecommand{\doi}[1]{\url{https://doi.org/#1}}
\bibcommenthead

\bibitem[\protect\citeauthoryear{Brewster}{1992}]{brewster1992thermal}
Brewster MQ.
\newblock Thermal radiative transfer and properties.
\newblock John Wiley \& Sons; 1992.

\bibitem[\protect\citeauthoryear{Howell et~al.}{2020}]{howell2020thermal}
Howell JR, Meng{\"u}{\c{c}} MP, Daun K, Siegel R.
\newblock Thermal radiation heat transfer.
\newblock CRC press; 2020.

\bibitem[\protect\citeauthoryear{Fan et~al.}{2020}]{fan2020near}
Fan D, Burger T, McSherry S, Lee B, Lenert A, Forrest SR.
\newblock Near-perfect photon utilization in an air-bridge thermophotovoltaic
  cell.
\newblock Nature. 2020;586(7828):237--241.

\bibitem[\protect\citeauthoryear{LaPotin
  et~al.}{2022}]{lapotin2022thermophotovoltaic}
LaPotin A, Schulte KL, Steiner MA, Buznitsky K, Kelsall CC, Friedman DJ, et~al.
\newblock Thermophotovoltaic efficiency of 40\%.
\newblock Nature. 2022;604(7905):287--291.

\bibitem[\protect\citeauthoryear{Raman et~al.}{2014}]{raman2014passive}
Raman AP, Anoma MA, Zhu L, Rephaeli E, Fan S.
\newblock Passive radiative cooling below ambient air temperature under direct
  sunlight.
\newblock Nature. 2014;515(7528):540--544.

\bibitem[\protect\citeauthoryear{Heo et~al.}{2020}]{heo2020janus}
Heo SY, Lee GJ, Kim DH, Kim YJ, Ishii S, Kim MS, et~al.
\newblock A Janus emitter for passive heat release from enclosures.
\newblock Science advances. 2020;6(36):eabb1906.

\bibitem[\protect\citeauthoryear{Menon et~al.}{2020}]{menon2020enhanced}
Menon AK, Haechler I, Kaur S, Lubner S, Prasher RS.
\newblock Enhanced solar evaporation using a photo-thermal umbrella for
  wastewater management.
\newblock Nature Sustainability. 2020;3(2):144--151.

\bibitem[\protect\citeauthoryear{Ni et~al.}{2016}]{ni2016steam}
Ni G, Li G, Boriskina SV, Li H, Yang W, Zhang T, et~al.
\newblock Steam generation under one sun enabled by a floating structure with
  thermal concentration.
\newblock Nature Energy. 2016;1(9):1--7.

\bibitem[\protect\citeauthoryear{Weinstein
  et~al.}{2015}]{weinstein2015concentrating}
Weinstein LA, Loomis J, Bhatia B, Bierman DM, Wang EN, Chen G.
\newblock Concentrating solar power.
\newblock Chemical reviews. 2015;115(23):12797--12838.

\bibitem[\protect\citeauthoryear{He et~al.}{2020}]{he2020perspective}
He YL, Qiu Y, Wang K, Yuan F, Wang WQ, Li MJ, et~al.
\newblock Perspective of concentrating solar power.
\newblock Energy. 2020;198:117373.

\bibitem[\protect\citeauthoryear{Wiecha et~al.}{2021}]{wiecha2021deep}
Wiecha PR, Arbouet A, Girard C, Muskens OL.
\newblock Deep learning in nano-photonics: inverse design and beyond.
\newblock Photonics Research. 2021;9(5):B182--B200.

\bibitem[\protect\citeauthoryear{Ma et~al.}{2021}]{ma2021deep}
Ma W, Liu Z, Kudyshev ZA, Boltasseva A, Cai W, Liu Y.
\newblock Deep learning for the design of photonic structures.
\newblock Nature Photonics. 2021;15(2):77--90.

\bibitem[\protect\citeauthoryear{Mao et~al.}{2021}]{mao2021inverse}
Mao S, Cheng L, Zhao C, Khan FN, Li Q, Fu H.
\newblock Inverse design for silicon photonics: from iterative optimization
  algorithms to deep neural networks.
\newblock Applied Sciences. 2021;11(9):3822.

\bibitem[\protect\citeauthoryear{Wang et~al.}{2021}]{wang2021intelligent}
Wang N, Yan W, Qu Y, Ma S, Li SZ, Qiu M.
\newblock Intelligent designs in nanophotonics: from optimization towards
  inverse creation.
\newblock PhotoniX. 2021;2:1--35.

\bibitem[\protect\citeauthoryear{Xu et~al.}{2021}]{xu2021improved}
Xu X, Sun C, Li Y, Zhao J, Han J, Huang W.
\newblock An improved tandem neural network for the inverse design of
  nanophotonics devices.
\newblock Optics Communications. 2021;481:126513.

\bibitem[\protect\citeauthoryear{Park et~al.}{2024}]{park2024tnn}
Park M, Grb{\v{c}}i{\'c} L, Motameni P, Song S, Singh A, Malagrino D, et~al.
\newblock Inverse Design of Photonic Surfaces via High throughput Femtosecond
  Laser Processing and Tandem Neural Networks.
\newblock Advanced Science. 2024;p. 2401951.

\bibitem[\protect\citeauthoryear{Ma et~al.}{2022}]{ma2022benchmarking}
Ma T, Tobah M, Wang H, Guo LJ.
\newblock Benchmarking deep learning-based models on nanophotonic inverse
  design problems.
\newblock Opto-Electronic Science. 2022;1(1):210012--1.

\bibitem[\protect\citeauthoryear{Jiang and Fan}{2020}]{jiang2020simulator}
Jiang J, Fan JA.
\newblock Simulator-based training of generative neural networks for the
  inverse design of metasurfaces.
\newblock Nanophotonics. 2020;9(5):1059--1069.

\bibitem[\protect\citeauthoryear{Kudyshev et~al.}{2020}]{kudyshev2020machinea}
Kudyshev ZA, Kildishev AV, Shalaev VM, Boltasseva A.
\newblock Machine-learning-assisted metasurface design for high-efficiency
  thermal emitter optimization.
\newblock Applied Physics Reviews. 2020;7(2).

\bibitem[\protect\citeauthoryear{Liu et~al.}{2020}]{liu2020topological}
Liu Z, Zhu Z, Cai W.
\newblock Topological encoding method for data-driven photonics inverse design.
\newblock Optics express. 2020;28(4):4825--4835.

\bibitem[\protect\citeauthoryear{Liu et~al.}{2021}]{liu2021tackling}
Liu Z, Zhu D, Raju L, Cai W.
\newblock Tackling photonic inverse design with machine learning.
\newblock Advanced Science. 2021;8(5):2002923.

\bibitem[\protect\citeauthoryear{Habibi et~al.}{2023}]{habibi2023actually}
Habibi M, Wang J, Fuge M.
\newblock When Is it Actually Worth Learning Inverse Design?
\newblock In: International Design Engineering Technical Conferences and
  Computers and Information in Engineering Conference. vol. 87301. American
  Society of Mechanical Engineers; 2023. p. V03AT03A025.

\bibitem[\protect\citeauthoryear{Zhu et~al.}{2023}]{zhu2023inverse}
Zhu Y, Chen Y, Gorsky S, Shubitidze T, Dal~Negro L.
\newblock Inverse design of functional photonic patches by adjoint optimization
  coupled to the generalized Mie theory.
\newblock JOSA B. 2023;40(7):1857--1874.

\bibitem[\protect\citeauthoryear{Gershnabel
  et~al.}{2022}]{gershnabel2022reparameterization}
Gershnabel E, Chen M, Mao C, Wang EW, Lalanne P, Fan JA.
\newblock Reparameterization Approach to Gradient-Based Inverse Design of
  Three-Dimensional Nanophotonic Devices.
\newblock ACS Photonics. 2022;10(4):815--823.

\bibitem[\protect\citeauthoryear{Hughes et~al.}{2018}]{hughes2018adjoint}
Hughes TW, Minkov M, Williamson IA, Fan S.
\newblock Adjoint method and inverse design for nonlinear nanophotonic devices.
\newblock ACS Photonics. 2018;5(12):4781--4787.

\bibitem[\protect\citeauthoryear{Minkov et~al.}{2020}]{minkov2020inverse}
Minkov M, Williamson IA, Andreani LC, Gerace D, Lou B, Song AY, et~al.
\newblock Inverse design of photonic crystals through automatic
  differentiation.
\newblock Acs Photonics. 2020;7(7):1729--1741.

\bibitem[\protect\citeauthoryear{Lalau-Keraly et~al.}{2013}]{lalau2013adjoint}
Lalau-Keraly CM, Bhargava S, Miller OD, Yablonovitch E.
\newblock Adjoint shape optimization applied to electromagnetic design.
\newblock Optics express. 2013;21(18):21693--21701.

\bibitem[\protect\citeauthoryear{Wang et~al.}{2020}]{wang2020inverse}
Wang K, Ren X, Chang W, Lu L, Liu D, Zhang M.
\newblock Inverse design of digital nanophotonic devices using the adjoint
  method.
\newblock Photonics Research. 2020;8(4):528--533.

\bibitem[\protect\citeauthoryear{Molesky et~al.}{2018}]{molesky2018inverse}
Molesky S, Lin Z, Piggott AY, Jin W, Vuckovi{\'c} J, Rodriguez AW.
\newblock Inverse design in nanophotonics.
\newblock Nature Photonics. 2018;12(11):659--670.

\bibitem[\protect\citeauthoryear{Deng et~al.}{2022}]{deng2022hybrid}
Deng L, Xu Y, Liu Y.
\newblock Hybrid inverse design of photonic structures by combining
  optimization methods with neural networks.
\newblock Photonics and Nanostructures-Fundamentals and Applications.
  2022;52:101073.

\bibitem[\protect\citeauthoryear{Kudyshev et~al.}{2020}]{kudyshev2020machineb}
Kudyshev ZA, Kildishev AV, Shalaev VM, Boltasseva A.
\newblock Machine learning--assisted global optimization of photonic devices.
\newblock Nanophotonics. 2020;10(1):371--383.

\bibitem[\protect\citeauthoryear{Hegde}{2019}]{hegde2019photonics}
Hegde RS.
\newblock Photonics inverse design: pairing deep neural networks with
  evolutionary algorithms.
\newblock IEEE Journal of Selected Topics in Quantum Electronics.
  2019;26(1):1--8.

\bibitem[\protect\citeauthoryear{Yao et~al.}{2023}]{yao2023inverse}
Yao T, Huang T, Yan B, Ge M, Yin J, Peng C, et~al.
\newblock Inverse design of dispersion for photonic devices based on LSTM and
  gradient-free optimization algorithms hybridization.
\newblock JOSA B. 2023;40(6):1525--1532.

\bibitem[\protect\citeauthoryear{Kudyshev
  et~al.}{2021}]{kudyshev2021optimizing}
Kudyshev ZA, Kildishev AV, Shalaev VM, Boltasseva A.
\newblock Optimizing startshot lightsail design: A generative network-based
  approach.
\newblock ACS Photonics. 2021;9(1):190--196.

\bibitem[\protect\citeauthoryear{Yeung et~al.}{2022}]{yeung2022deepadjoint}
Yeung C, Pham B, Tsai R, Fountaine KT, Raman AP.
\newblock DeepAdjoint: an all-in-one photonic inverse design framework
  integrating data-driven machine learning with optimization algorithms.
\newblock ACS Photonics. 2022;10(4):884--891.

\bibitem[\protect\citeauthoryear{Yeung et~al.}{2022}]{yeung2022enhancing}
Yeung C, Ho D, Pham B, Fountaine KT, Zhang Z, Levy K, et~al.
\newblock Enhancing adjoint optimization-based photonic inverse design with
  explainable machine learning.
\newblock ACS Photonics. 2022;9(5):1577--1585.

\bibitem[\protect\citeauthoryear{Grbcic et~al.}{2024}]{grbcic2024ensemble}
Grbcic L, Park M, Elzouka M, Prasher R, Müller J, Grigoropoulos CP, et~al.:
  Inverse design of photonic surfaces on Inconel via multi-fidelity machine
  learning ensemble framework and high throughput femtosecond laser processing.

\bibitem[\protect\citeauthoryear{Noack and Ushizima}{2023}]{noack2023methods}
Noack M, Ushizima D.
\newblock Methods and Applications of Autonomous Experimentation.
\newblock CRC Press (Unlimited); 2023.

\bibitem[\protect\citeauthoryear{H{\"a}se et~al.}{2019}]{hase2019next}
H{\"a}se F, Roch LM, Aspuru-Guzik A.
\newblock Next-generation experimentation with self-driving laboratories.
\newblock Trends in Chemistry. 2019;1(3):282--291.

\bibitem[\protect\citeauthoryear{Kochenderfer and
  Wheeler}{2019}]{kochenderfer2019algorithms}
Kochenderfer MJ, Wheeler TA.
\newblock Algorithms for optimization.
\newblock Mit Press; 2019.

\bibitem[\protect\citeauthoryear{Paria et~al.}{2022}]{paria2022greedy}
Paria B, P{\`o}czos B, Ravikumar P, Schneider J, Suggala AS.
\newblock Be greedy--a simple algorithm for blackbox optimization using neural
  networks.
\newblock In: ICML2022 Workshop on Adaptive Experimental Design and Active
  Learning in the Real World; 2022. .

\bibitem[\protect\citeauthoryear{Wang et~al.}{2014}]{wang2014evaluation}
Wang C, Duan Q, Gong W, Ye A, Di Z, Miao C.
\newblock An evaluation of adaptive surrogate modeling based optimization with
  two benchmark problems.
\newblock Environmental Modelling \& Software. 2014;60:167--179.

\bibitem[\protect\citeauthoryear{Elzouka
  et~al.}{2020}]{elzouka2020interpretable}
Elzouka M, Yang C, Albert A, Prasher RS, Lubner SD.
\newblock Interpretable forward and inverse design of particle spectral
  emissivity using common machine-learning models.
\newblock Cell Reports Physical Science. 2020;1(12).

\bibitem[\protect\citeauthoryear{Liu et~al.}{2018}]{liu2018remarks}
Liu H, Cai J, Ong YS.
\newblock Remarks on multi-output Gaussian process regression.
\newblock Knowledge-Based Systems. 2018;144:102--121.

\bibitem[\protect\citeauthoryear{Pedregosa et~al.}{2011}]{pedregosa2011scikit}
Pedregosa F, Varoquaux G, Gramfort A, Michel V, Thirion B, Grisel O, et~al.
\newblock Scikit-learn: Machine learning in Python.
\newblock the Journal of machine Learning research. 2011;12:2825--2830.

\bibitem[\protect\citeauthoryear{Jones}{2001}]{jones2001taxonomy}
Jones DR.
\newblock A taxonomy of global optimization methods based on response surfaces.
\newblock Journal of global optimization. 2001;21:345--383.

\bibitem[\protect\citeauthoryear{Breiman}{2001}]{breiman2001random}
Breiman L.
\newblock Random forests.
\newblock Machine learning. 2001;45:5--32.

\bibitem[\protect\citeauthoryear{Kennedy and
  Eberhart}{1995}]{kennedy1995particle}
Kennedy J, Eberhart R.
\newblock Particle swarm optimization.
\newblock In: Proceedings of ICNN'95-international conference on neural
  networks. vol.~4. ieee; 1995. p. 1942--1948.

\bibitem[\protect\citeauthoryear{Storn and Price}{1997}]{storn1997differential}
Storn R, Price K.
\newblock Differential evolution--a simple and efficient heuristic for global
  optimization over continuous spaces.
\newblock Journal of global optimization. 1997;11:341--359.

\bibitem[\protect\citeauthoryear{Audet and Dennis~Jr}{2006}]{audet2006mesh}
Audet C, Dennis~Jr JE.
\newblock Mesh adaptive direct search algorithms for constrained optimization.
\newblock SIAM Journal on optimization. 2006;17(1):188--217.

\bibitem[\protect\citeauthoryear{Lagarias
  et~al.}{1998}]{lagarias1998convergence}
Lagarias JC, Reeds JA, Wright MH, Wright PE.
\newblock Convergence properties of the Nelder--Mead simplex method in low
  dimensions.
\newblock SIAM Journal on optimization. 1998;9(1):112--147.

\bibitem[\protect\citeauthoryear{Shahriari et~al.}{2015}]{shahriari2015taking}
Shahriari B, Swersky K, Wang Z, Adams RP, De~Freitas N.
\newblock Taking the human out of the loop: A review of Bayesian optimization.
\newblock Proceedings of the IEEE. 2015;104(1):148--175.

\bibitem[\protect\citeauthoryear{Ivic et~al.}{2024}]{Indago}
Ivic S, Druzeta S, Grbcic L.: Indago v0.5.0.
\newblock PyPI.
\newblock Available from: \url{https://pypi.org/project/Indago/}.

\bibitem[\protect\citeauthoryear{Bayoumy}{2022}]{OMADS_AB}
Bayoumy A.: OMADS.
\newblock Github.
\newblock Available from: \url{https://github.com/Ahmed-Bayoumy/OMADS}.

\bibitem[\protect\citeauthoryear{Virtanen et~al.}{2020}]{virtanen2020scipy}
Virtanen P, Gommers R, Oliphant TE, Haberland M, Reddy T, Cournapeau D, et~al.
\newblock SciPy 1.0: fundamental algorithms for scientific computing in Python.
\newblock Nature methods. 2020;17(3):261--272.

\bibitem[\protect\citeauthoryear{Zhan and Xing}{2020}]{zhan2020expected}
Zhan D, Xing H.
\newblock Expected improvement for expensive optimization: a review.
\newblock Journal of Global Optimization. 2020;78(3):507--544.

\bibitem[\protect\citeauthoryear{Head et~al.}{2020}]{head_2020_4014775}
Head T, Kumar M, Nahrstaedt H, Louppe G, Shcherbatyi I.:
  scikit-optimize/scikit-optimize.
\newblock Zenodo.
\newblock Available from: \url{https://doi.org/10.5281/zenodo.4014775}.

\end{thebibliography}

\end{document}